\documentclass[12pt]{article}

\usepackage{amsmath}
\usepackage{graphicx,psfrag,epsf}
\usepackage{enumerate}

\usepackage[sort,numbers]{natbib}

\usepackage{bbm}
\usepackage{amsfonts}
\usepackage{algpseudocode}
\usepackage{algorithm}
\usepackage{algorithmicx}
\usepackage{caption}
\usepackage{subfigure}
\usepackage{soul}
\usepackage{tikz}

\usepackage{tabularx,ragged2e,booktabs,caption}

\usepackage[english]{babel}
\usepackage{amsthm}
\usepackage{amsmath}
\usepackage{amsfonts}
\usepackage{amssymb}
\usepackage{makeidx}

\usepackage{graphicx}
\usepackage[normalem]{ulem}
\usepackage{color}
\usepackage{array}
\usepackage[colorlinks=true,linkcolor=blue,citecolor=blue,urlcolor=black]{hyperref}
\usepackage{pdflscape}

\usepackage[T1]{fontenc}
\usepackage[utf8]{inputenc}
\usepackage{tabularx,ragged2e,booktabs,caption}
\newcolumntype{C}[1]{>{\Centering}m{#1}}

\providecommand{\keywords}[1]{\textit{Keywords.} #1}
\newcommand{\blind}{0}

\def\mR{\mathbb{R}}

\def\real{\mR}
\def\bc{\begin{center}}
\def\ec{\end{center}}

\newcommand{\beq}{\begin{eqnarray}}
\newcommand{\eeq}{\end{eqnarray}}
\newcommand{\beqq}{\begin{eqnarray*}}
\newcommand{\eeqq}{\end{eqnarray*}}

\begin{document}

\def\spacingset#1{\renewcommand{\baselinestretch}%
{#1}\small\normalsize} \spacingset{1}

\if0\blind
{

  \title{\bf Vietoris-Rips Persistent Homology}
  \author{Sarit Agami\hspace{.2cm}\\
    Andrew and Erna Viterbi Faculty of  Electrical Engineering,\\
     Technion -- Israel Institute of Technology}
  \maketitle
} \fi

%

\begin{abstract}
Persistence diagrams are useful displays that give a summary information regarding the topological features of some phenomenon. Usually, only one persistence diagram is available, and replicated persistence diagrams are needed for statistical inference. One option for generating these replications is to fit a distribution for the points on the persistence diagram. The type of the relevant distribution depends on the way the persistence diagram is builded. There are two approaches for building the persistence diagram, one is based on the Vietoris-Rips complex, and the second is based on some fitted function such as the kernel density estimator. The two approaches yield a two dimensional persistence diagram, where the coordinates of each point are the 'birth' and 'death' times. For the first approach, however, the 'birth' time is zero for all the points that present the connected components of the phenomenon. In this paper we examine the distribution of the connected components when the persistence diagram is based on Vietoris-Rips complex. In addition, we study the behaviour of the connected components when the phenomenon is measured with noise.
\end{abstract}


\section{Introduction}
\label{sec:intro}
The aim of topological data analysis (TDA) is to provide methods for analyzing the complex topological and geometric structures underlying data. The data is assumed to be a finite set of points coming with a notion of distance, or similarity, between them. Given data, the topological features of the underlying space can be quantified via persistent homology (see, for example, \cite{carlsson2009, carlsson2014,Edelsbrunner2008, Edelsbrunner2010, Edelsbrunner2014,Zomorodian2005,Oudot2015,Ghrist2014}). One way to calculate the persistent homology is to consider some function $f$ on the space, for example the distance function or the kernel density estimator, and then to calculate the sub-level (or the super-level) sets $f^{-1}(-\infty,x]$ ($f^{-1}[x,\infty)$) of the function $f$. Using the sub-level (or the super-level) sets, the homology changes: new connected components can appear, existing components can be merged, loops can appear or be filled, etc. Persistent homology tracks these changes, identifies the appearing features and associated a lifetime to them. The resulting information is encoded as a set of intervals called a barcode, and each interval is called a bar. Equivalently, the barcode is a multi-set of points in $\mathbb{R}^2$ where the coordinates of each point are the starting and the end of the corresponding interval, which are also named the 'birth time' and the 'death time', respectively.  The obtained diagram is called the 'persistence diagram', where the connected components are the zero-th homology, or the $H_0$ points, the loops are the one-th homology, or the $H_1$ points, etc. A 'length' or a 'lifetime' of a given bar is defined as the death time minus the birth time. Another way to calculate the persistent homology is based on the Vietoris-Rips filtration, which is given by a union of growing balls centered at each of the original data. Starting with radius $r=0$, the union of the balls is the original data, where each point is considered as a connected component. This is translated into the persistence diagram by creating an interval for the birth for each of these features. As the value of $r$ increases, some of the balls start to overlap, that is, some of the connected components that get merged together are 'died'. The persistence diagram keeps track of these deaths, putting an end point to the corresponding intervals as they disappear. In other words, the persistence diagram can be seen as a multiscale topological signature encoding the homology of the union of the balls for all radii as well as its evolution across the values of $r$. By this building, the birth time is zero for all the points. Denote by $n$ the number of the data points. The number of the $H_0$ points is $n$ as well, when always these $n$ points include one point with infinity lifetime.
The building of the persistence diagram depends on 'maxscale parameter' which describes the maximum value of the Vietoris-Rips filtration. Which value of the maxscale is the best choosing depends on the behaviour of the data that describes the specific phenomenon. When the maxscale is not large enough, different $H_0$ persistence diagrams are obtained for the various values of maxscale, where the difference is among the longer bars. But as the maxscale becomes large enough, its value or larger values obtain the same $H_0$ persistence diagrams. Therefore, a possible way to determine the 'optimal' maxscale, is to compare the maximal length of the $H_0$ points with the given maxscale: if the maximal length is strictly smaller than the given maxscale, then this value is the optimal maxscale.
\section{Distribution of the Connected Components}
\subsection{General}
For studying the distribution of the points on the persistence diagram, \cite{adler2017, adler2018} suggested a parametric model when the persistence homology is based on the kernel density estimator. Once the persistence diagram is based on the Vietoris-Rips filtration, we can use this parametric model to describe the behaviour of the $H_1$ points. However, the $H_0$ points in that case cannot be described by the same model due to the zero birth times for all the points. From the other hand, the last property leads, concentrating on the death times only, to formulate the relevant distribution as one dimensional distribution. We expect to see, at least for low dimensional data, and while ignoring the $H_0$ point at infinity, a right-tailed distribution fitting.
This fitting can be done using some built-in procedures that exist in the softwares. We use the procedure \emph{'allfitdist'} in Matlab; this procedure suggests some distributions for each data, and the best fitted distribution can be chosen via the AIC and BIC criterions \ (cf.\cite{burnham}). Some of the suggested distributions have infinity support, while the support of the $H_0$ points (ignoring the point at infinity) is finite. Therefore the considered family of the appropriate distributions should be limited to those with finite support. When the most appropriate distribution is fitted, it is easy to obtain replicated $H_0$ persistence diagrams by generating random numbers from this fitted distribution. The next step will be to use these replications for statistical inference, such as identification of the topological signals that belongs to the original $H_0$ persistence diagram. Particularly, the replications enable to examine the significance of the points that are suspect as topological signals. For this purpose we use, similarly to \cite{adler2017,adler2018}, and \cite{agami2017}, the order statistics of the death times $T_j$, $j=1,...,(n-1)$. For each $T_j$, we can calculate its confidence interval and its p-value in a similar way of the 'percentile bootstrap' method. We studying the performance of these statistics by the examples in Section 3.
\\
\indent In many situations, the data is observed with some noise. The noise, as we show in the examples in Section 4, adds longer bars relative to the maximal length of the bars when the data is observed without noise.
Since a long bar presents a topological signal, it is important to recognize if each long bar comes from noise, or if it is a real topological signal. We study the behaviour of the $H_0$ bars in the setting of noisy data by the examples in Section 4.

\subsection{Clean Data}
We refer the data that is measured without any noise as clean data. For such data set, as we show in the considered examples in Section 3, usually the fitted distribution for the $H_0$ points is the beta distribution. This is true independently on the sample size $n$. The influence of $n$ on the $H_0$ points in the setting of clean data, is only on their lengths: as $n$ increases, the lengths decrease. This is reasonable since larger $n$ means more closer points to each other, and therefore existing components can be merged faster. Once we have the fitted distribution, we can use it to generate replications of the $H_0$ persistence diagram we have in hand. For all the considered examples in Section 3, we get the correct number of significant connected components using the procedure with the order statistics $T_j$ that was mentioned above. That is, these statistics perform well for the identification of the $H_0$ topological signals.
As for the goodness of the distribution fitting, we compare real $H_0$ persistence diagrams with their corresponded simulated $H_0$ persistence diagrams (where again, the simulated persistence diagrams are random numbers from the fitted distribution). This comparison can be done by using some summary statistics for each persistence diagram. We use the kurtosis and the skewness statistics; we calculate these statistics for each real and each simulated persistence diagram, and then compare their distributions over the real and the simulated persistence diagrams. We can see close behaviours of these distributions.

\subsection{Noisy Data}
Adding noise to the data requires a larger value of maxscale in order to capture the shape of the data. If the maxscale is small and not large enough, then the maximal $H_0$ death time will be the value of the maxscale. The 'optimal' value of the maxscale depends on the amount of the noisy points.  Let $M$ be the percent of the $n$ data points that are measured with noise. Particulary, in our considered examples in Section 4, the noise is additive and generated from the bivariate normal distribution with zero mean and the identity matrix divided by 9 for the covariance. Denote by $c_{max}$ the maximal $H_0$ death time under the clean data. We can define a 'long bar' as the bar that his death time is greater than $c_{max}$, and similarly, a 'short bar' as the bar that his death time is smaller or equals $c_{max}$. Based on the results of the examples in Section 4, we can indicate on three properties:
\begin{itemize}
\item[(i)] Given $n$, the ratio of the long bars relative to the short bars increases as $M$ increases. This is somehow with the contrary to the fact that larger $n$ decreases the death times. But, looking inside the $H_0$ points distribution, we can see that larger $n$ decreases the death times in the major part of the distribution, when this is also depends on the maxscale value. Therefore, the maximal length does not necessarily become larger as $M$ increases.
\item[(ii)] For a given $M$, the ratio of the long bars relative to the short bars increases as $n$ increases. Therefore we can see in the figures in Section 4 that for a given $M$, the tail of the fitted distribution becomes thinner as $n$ increases, this is due to some point which is isolated from the other points.
\item[(iii)] For a given $M$ and $n$, the ratio of the long bars relative to the short bars is the same for the various values of the maxscale. Although, there is some difference among the long bars for the different values of the maxscale. Usually, the difference is in their maximal death times. But, when the maxscale is small relative to the optimal maxscale, the difference can start at the $95$-th percentile of the $H_0$ death times, and even at a lower percentile, depends on how smaller is the maxscale relative to the optimal one.

\end{itemize}
In this setting of noisy data, the parametric beta distribution (or other distribution with a finite support) is no longer behave as the best distribution to describe the behaviour of the connected components. We demonstrate it in Section 4.

\section{Examples of Clean Data}
We present now five examples, which are different in their topological structure.
The aim is to explore the behaviour of the $H_0$ points in a given single persistence diagram generated by the Vietoris-Rips complex.
The examples that we shall treat are one circle, a collection of two concentric circles in the plane, two distinct circles, 2-sphere, and 3-torus. The first three examples are two-dimensional objects, whereas the two later examples are higher-dimension objects.
In this section we concentrate on clean data for each example, that is, the examples are measured without noise. On the next section we examine three of these five examples but with additive noise.

\subsection{One Circle}
\subsubsection{Description and the Distribution}
\begin{figure}[h!]
\bc
\text{One circle with radius 1}\par\medskip
\includegraphics[width=1.45in, height=1.45in]{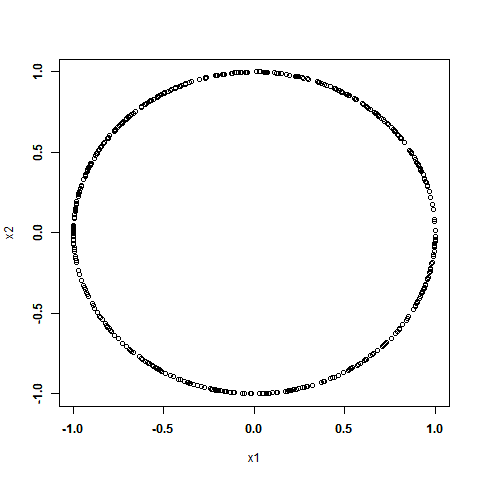} \hskip0.01truein
\includegraphics[width=1.45in, height=1.45in]{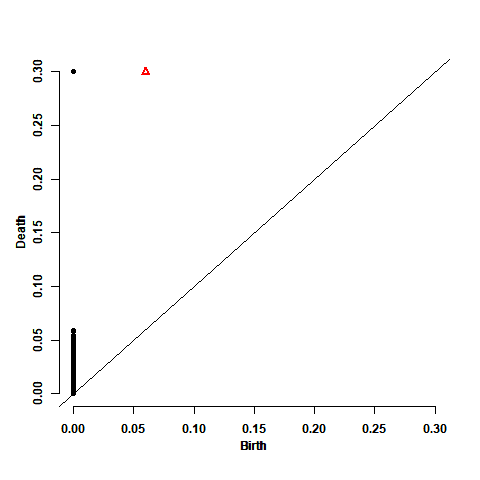} \hskip0.01truein
\\
\includegraphics[width=1.45in, height=1.45in]{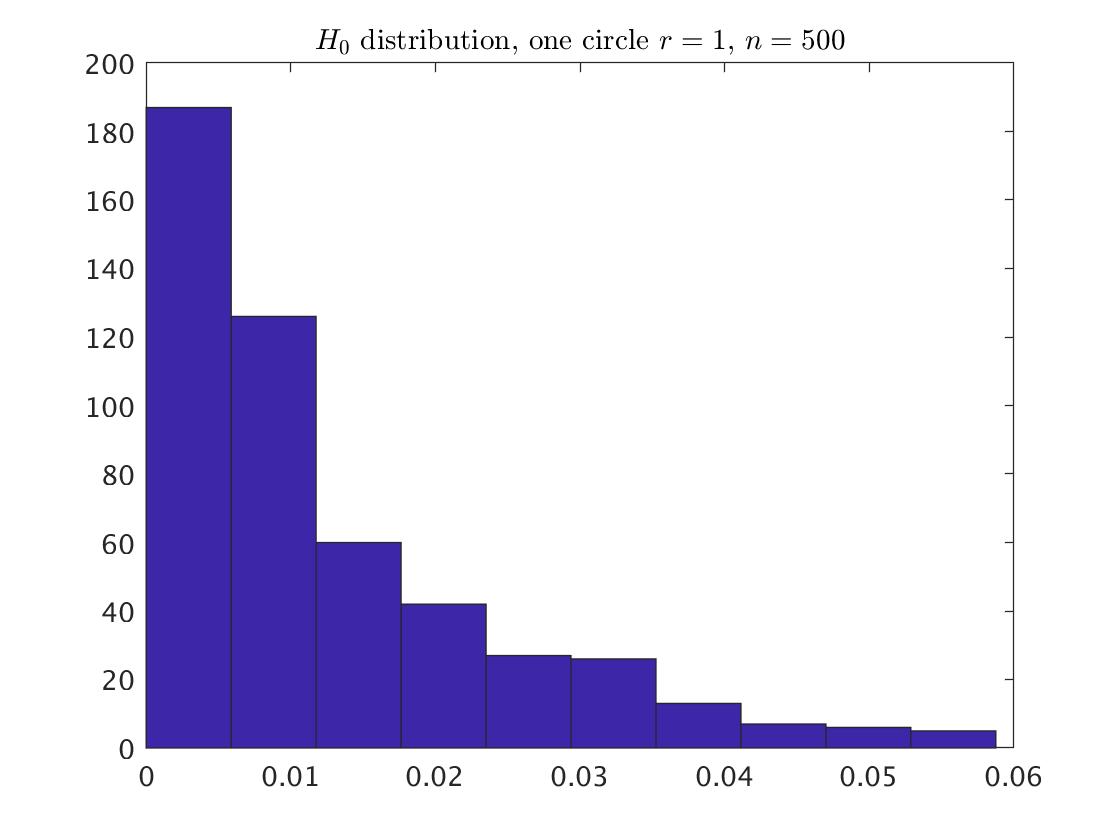} \hskip0.01truein
\includegraphics[width=1.45in, height=1.45in]{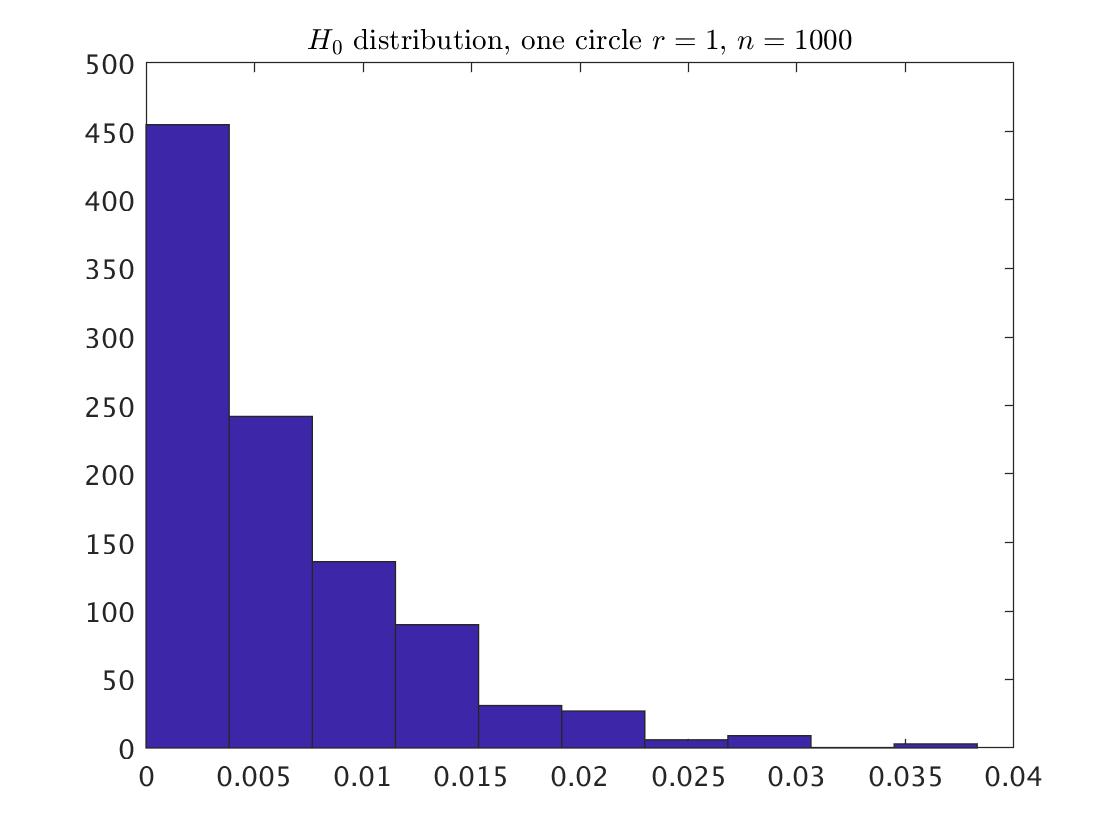} \hskip0.01truein
\includegraphics[width=1.45in, height=1.45in]{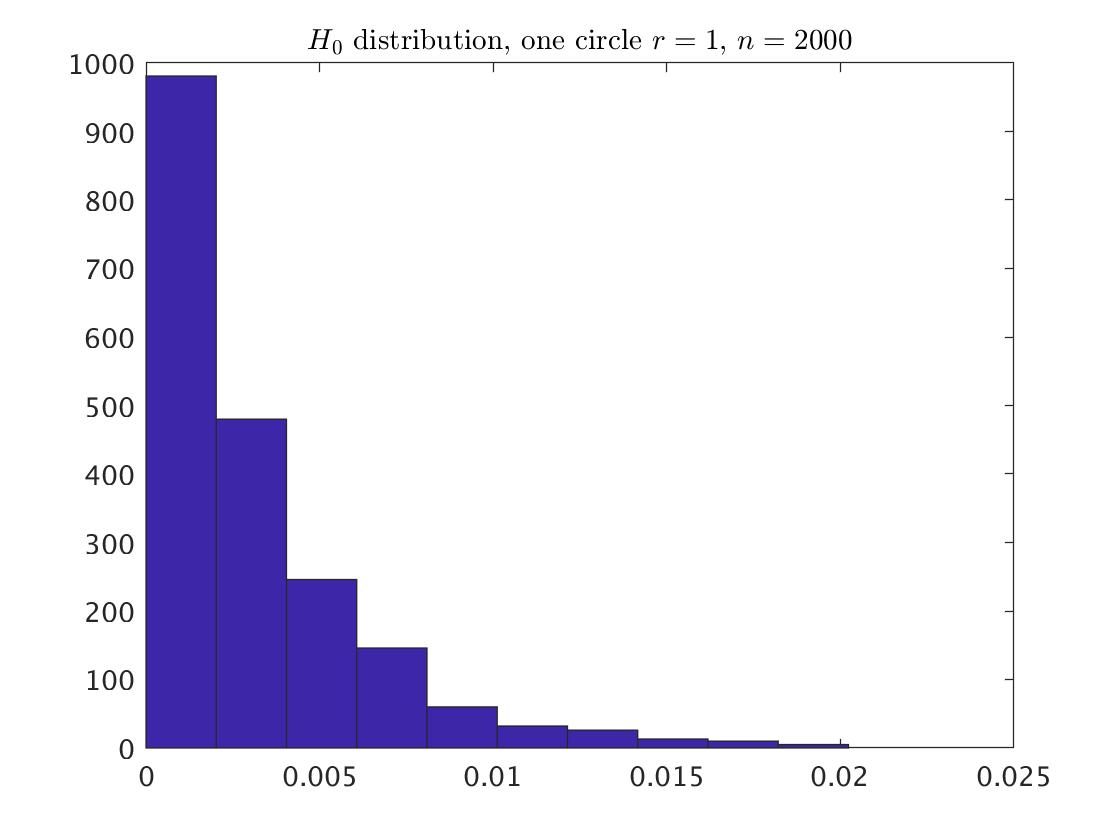} \hskip0.01truein
\\
\text{One circle with radius 3}\par\medskip
\includegraphics[width=1.45in, height=1.45in]{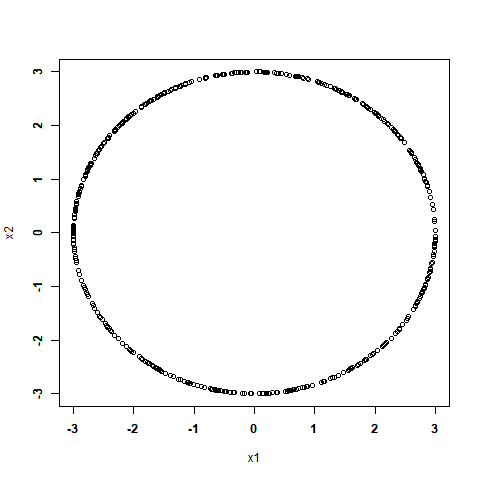} \hskip0.01truein
\includegraphics[width=1.45in, height=1.45in]{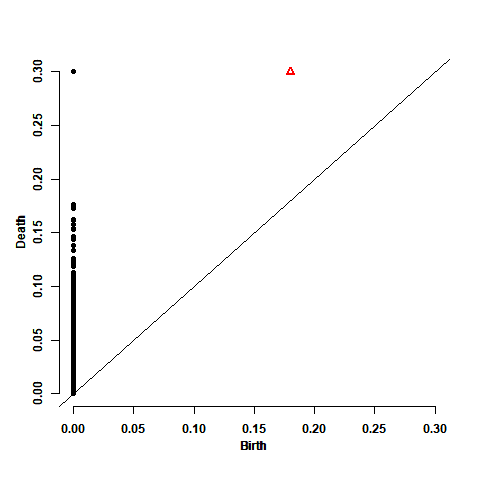} \hskip0.01truein
\\
\includegraphics[width=1.45in, height=1.45in]{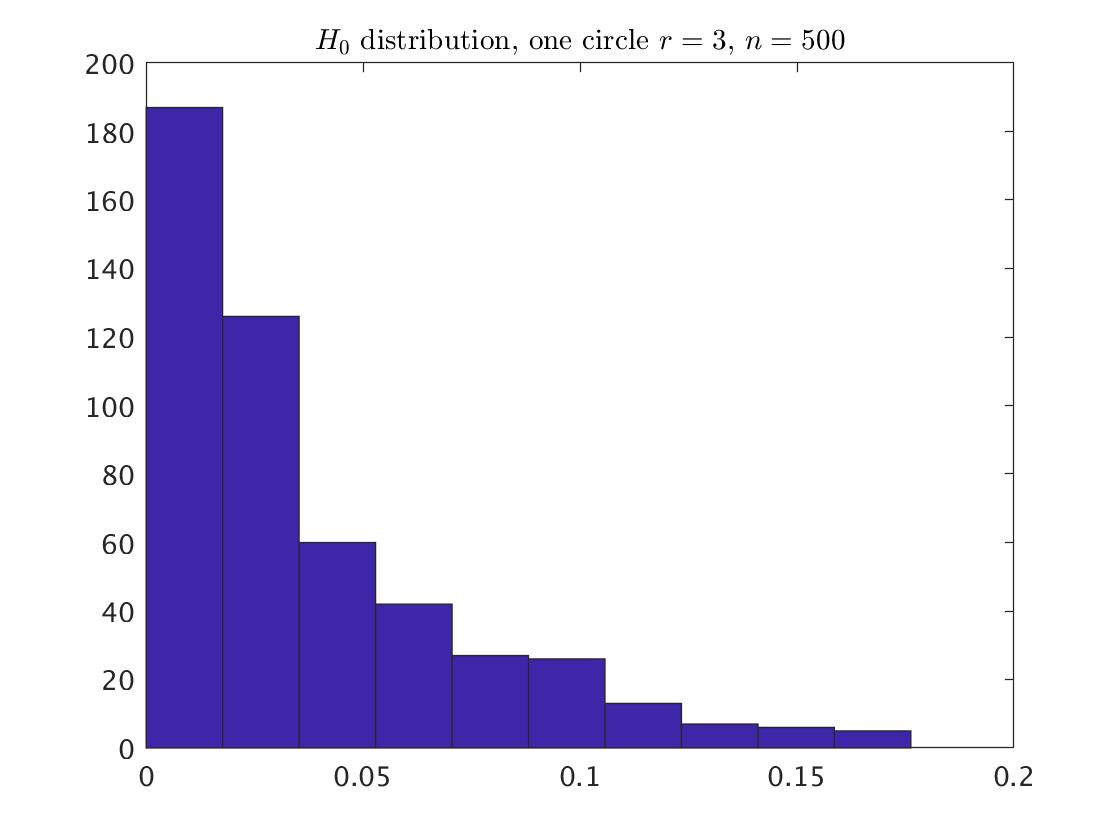} \hskip0.01truein
\includegraphics[width=1.45in, height=1.45in]{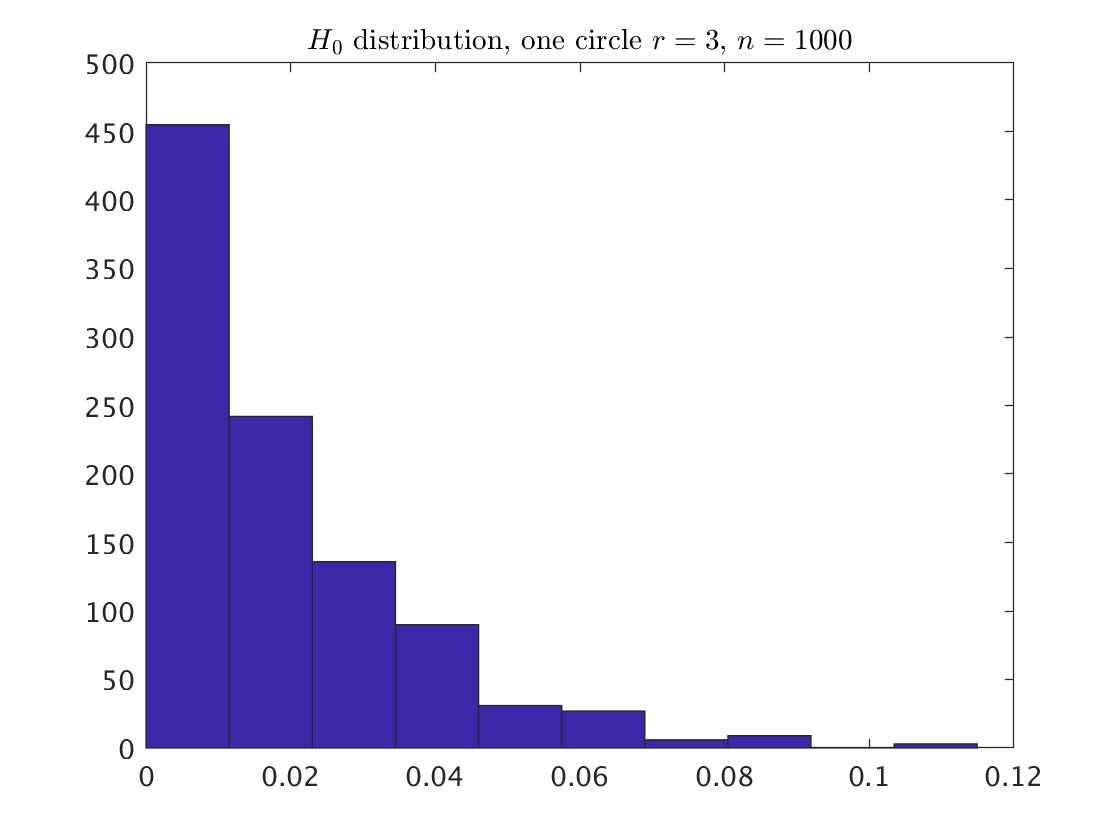} \hskip0.01truein
\includegraphics[width=1.45in, height=1.45in]{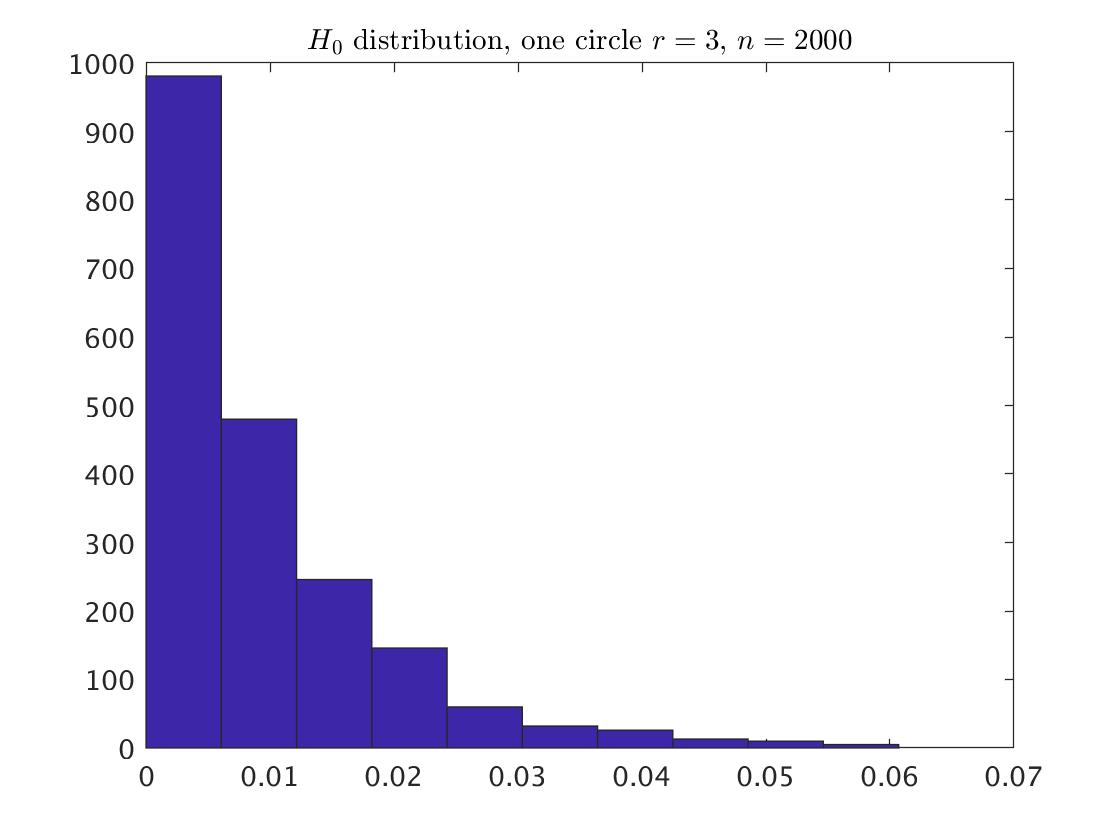} \hskip0.01truein
\ec
\caption{\footnotesize
 The plots describe one circle with radius 1 and radius 3. For each case, the plots from left to right are: the plot of the data, the corresponded persistence diagram based on the Vietoris-Rips filtration, the corresponded $H_0$ distribution without the point at infinity, and two more histograms that describe the $H_0$ distribution (without the point at infinity) for larger samples.
In the plot of the persistence diagram, the black circles are the $H_0$ persistence points, while the red triangles are the $H_1$ points.
}
\label{fig:onecircler1_M0}
\end{figure}

The following example includes a random sample of $n=500$ points from one circle with radius $r=1$, comparing with $r=3$. These circles are described in Fig.\ \ref{fig:onecircler1_M0}. For each circle we see to its right the corresponding persistence diagram of the Vietoris-Rips filtration using maxscale=0.3. This diagram contains 500 points of $H_0$, with the black circles indicating the $H_0$ points and the red triangles corresponding to $H_1$, in both cases trying to capture the underlying homology of the one circle. As described above, each point in the diagram is a `birth-death' pair. We expect to see one black circle and one red triangle somewhat isolated from the other points in the diagram, and this is in fact the case. Note that visually, due to page constrain, the location of the point at infinity of $H_0$ is at (0, maxscale) (the persistence diagram was calculated and plotted by using the package 'TDA' in $R$ software). Next to that plot we have the corresponding histogram of $H_0$ death times without including the point at infinity. The next histograms for each circle describe the $H_0$ points distributions that are corresponded to the larger samples $n=1,000$ and $n=2,000$, respectively.
At this example it is enough to use maxscale equals to 0.3. Maxscale equals to or greater than 0.3 obtain the same $H_0$ persistence diagram. This relative small value of maxscale is enough since all the points are sitting exactly on the circle, and there are no points inside or outside the circle line. However, this is not the case for the $H_1$ one point, which has a different value when using maxscale=0.3, and maxscale=5, for example.

\normalsize

\newpage
Looking at the $H_0$ persistence diagram and the histograms that describe the $H_0$ distribution for the both circles, it seems that for a given $n$, the shape of the distributions for $r=1$ and $r=3$ are the same. The difference between them is in the values of the death times: the death times under $r=3$ are three times larger than those under $r=1$, and therefore yield different distributions. For fitting the parametric distribution, we first remove the point at infinity, and then use the \emph{'allfitdist'}  procedure in Matlab. For example for $n=500$ and $r=1$, the suggested finite support distributions are the beta distribution with support=[0,1], and the generalized Pareto with support=[0,0.256]. Taking the minimal BIC and AIC measures yield that the best fitting is the beta distribution with the parameters $[a,b]= [0.976,77.302]$. Similarly, the best fitting for the $H_0$ points corresponding to the circle with $r=3$ and $n=500$ is the beta distribution with $[a,b]= [0.957,24.658]$. The fitted distributions for the larger samples $n=1,000$ and $n=2,000$ are the beta distribution as well, where the relevant parameters are summarized in Appendix A.1.

\subsubsection{Goodness of Fit}
\begin{figure}[h!]
\bc
\includegraphics[width=1.45in, height=1.45in]{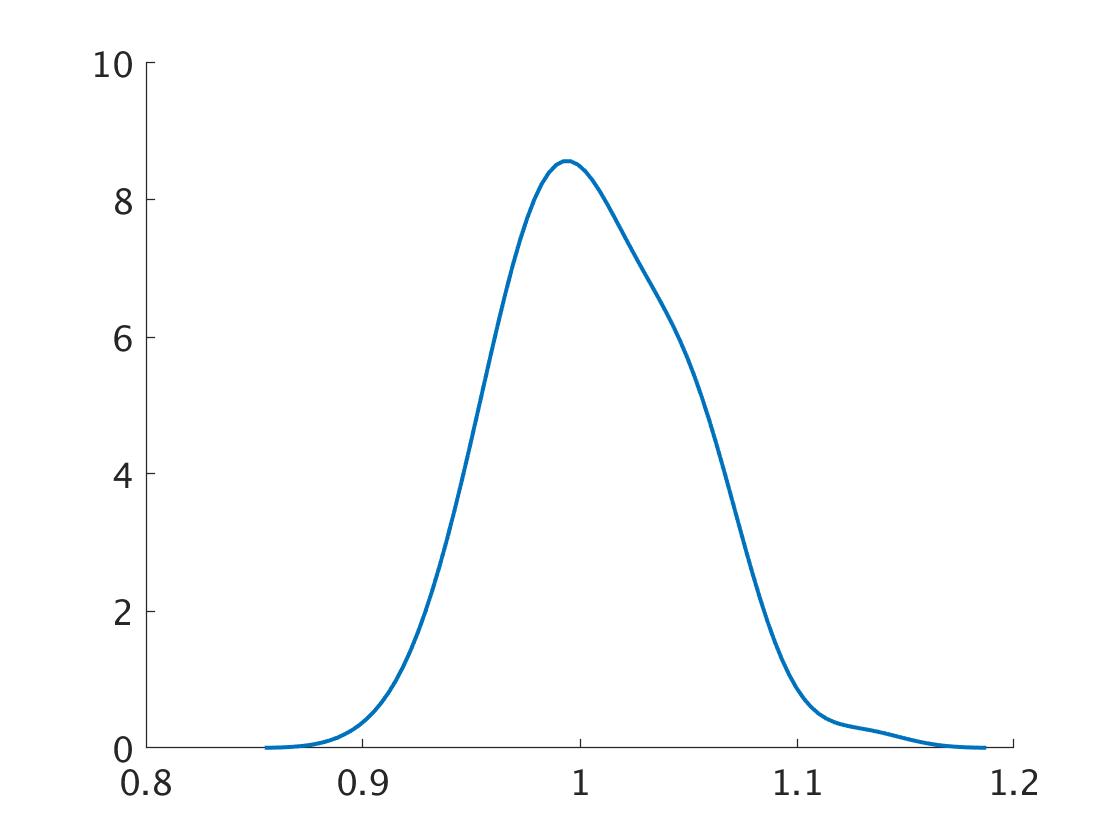} \hskip0.01truein
\includegraphics[width=1.45in, height=1.45in]{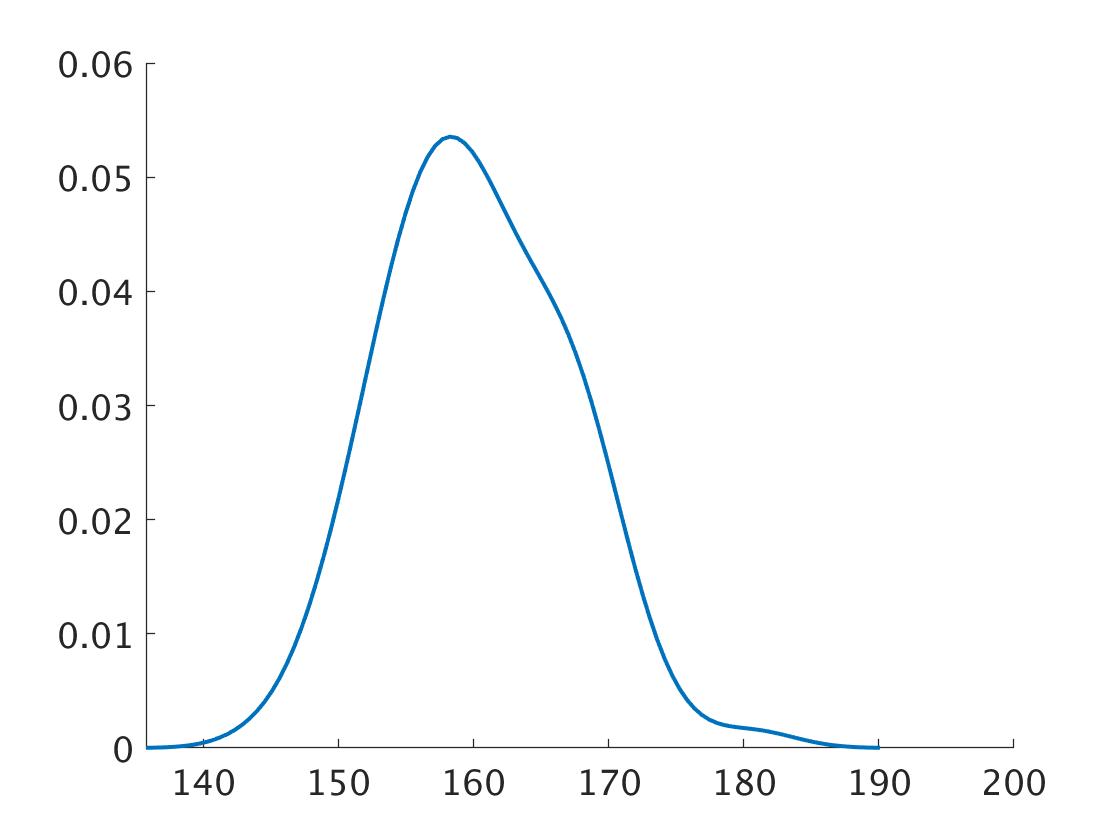} \hskip0.01truein
\includegraphics[width=1.45in, height=1.45in]{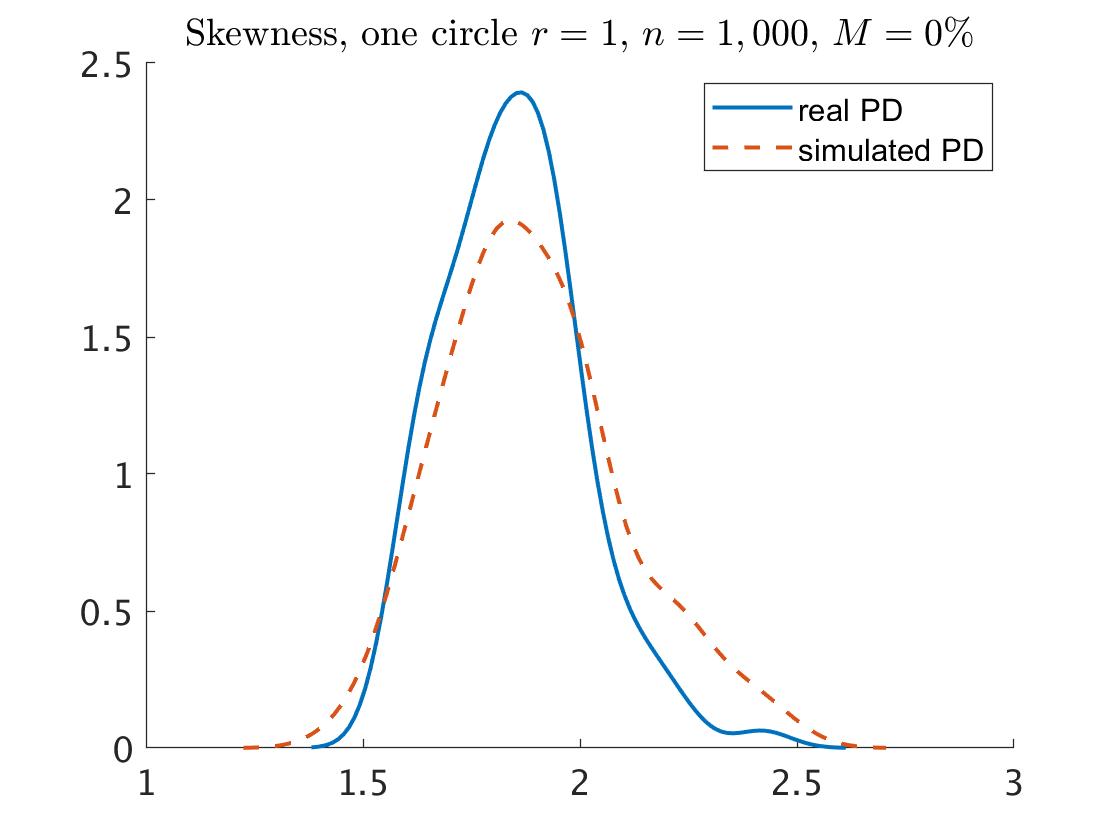} \hskip0.01truein
\includegraphics[width=1.45in, height=1.45in]{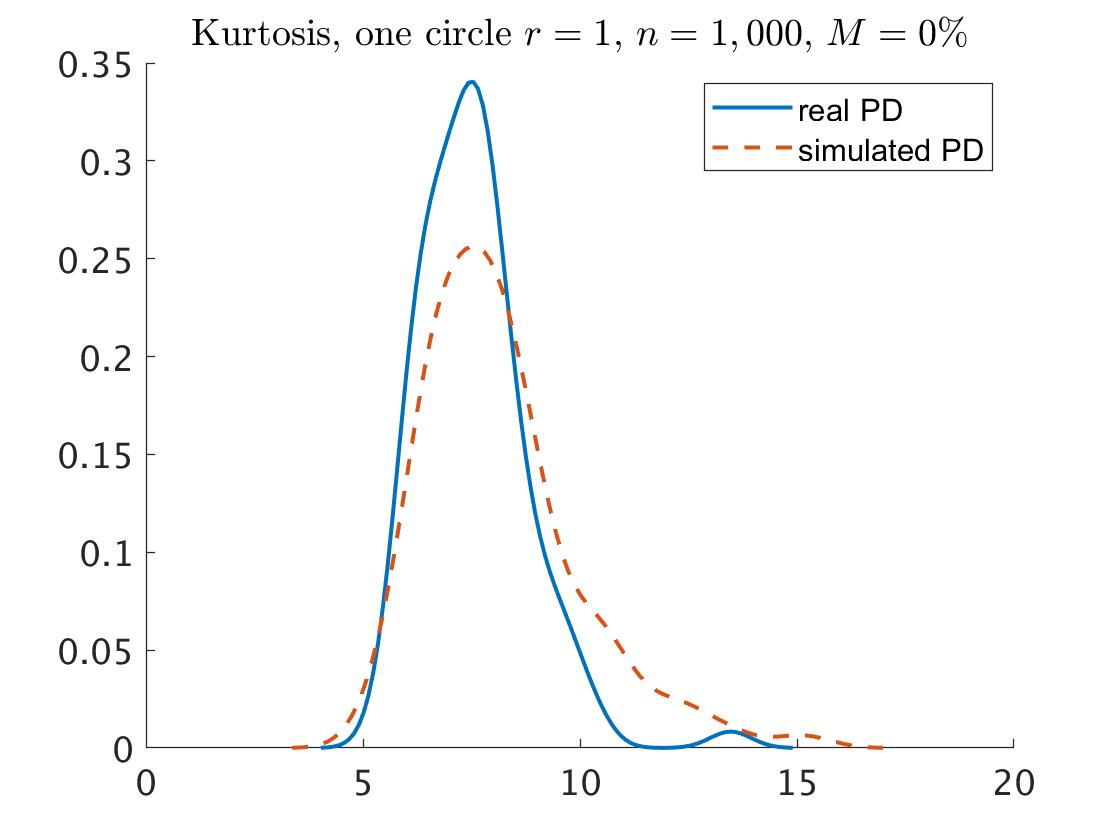} \hskip0.01truein
\ec
\caption{\footnotesize
 The first two plots are smoothed empirical densities for the parameters $a$ and $b$ of the beta distribution. Each beta distribution (over 86 cases) describes the $H_0$ points coming from one circle with $n=1,000$ and $r=1$. The next two plots are smoothed empirical densities of the skewness and kurtosis distributions over 100 real and fitted distributions for the $H_0$ points coming from one circle with $n=1,000$ and $r=1$.
}
\label{fig:onecircle_par_dist}
\end{figure}

In order to test how well the fitted distribution matches the $H_0$ persistence diagram in hand, we generated 100 collections of samples from one circle according to the same procedure that generated the original data, with $n=1,000$ and $r=1$, and for each one we fitted the best distribution to the $H_0$ points. The fitted distribution for 86 collections was the beta distribution, while for the rest 14 collections the fitted distribution was the generalized Pareto distribution. The first two plots of Fig. \ref{fig:onecircle_par_dist} show the (smoothed) empirical densities of the resulting  parameters estimates for the beta distribution (over the 86 cases). Overall, the results indicate that the distribution fitting is stable, with what seems to be  an acceptable  spread in the distribution of the estimates.
In addition we considered summary statistics of the skewness and kurtosis of the 100 $H_0$ persistence diagrams, to see how well the simulations replicate the statistical properties of the original $H_0$ persistence diagrams.
The results are presented in the two right plots of Fig. \ref{fig:onecircle_par_dist}.
The blue (full line) curves show the empirical probability densities of the skewness and kurtosis for each of the original 100 $H_0$ persistence diagrams. The red (dashed line) curves show the same phenomenon, but for the simulated  100 diagrams based on the fitted distributions. The curves for the each statistic are close, which indicate a good fitting of the distributions for the original $H_0$ persistence diagrams.
\normalsize
\subsubsection{Statistical Inference}
For the identification of topological signals of the one circle with $r=1$ and a sample size $n$, we generated 1,000 replicated $H_0$ persistence diagrams by taking random numbers from the above fitted beta distribution, each one contained $(n-1)$ points. We calculated the maximum statistic of the death times $T_1$, its confidence interval and its p-value. The results are summarized in Table \ref{table:OneCirc_stat}. We get that $T_1$ is insignificant, that is, no significant $H_0$ points among the set of the $H_0$ points without the point at infinity. Hence, the point at infinity is the only topological signal among the set of the $H_0$ points, as we hoped to find. The same result is obtained for one circle with $r=3$, as presented in Table \ref{table:OneCirc_stat}. Note that this result is independent on the value of $n$.

\begin{center}
\fontsize{8.5}{0.9}\selectfont
\captionof{table}{One circle - confidence interval and p-value}
\begin{tabular}{l|lcc|ccc|cccccc}
\\
&& $n=500$ &&& $n=1,000$&&& $n=2,000$\\\hline
\\
\\
\\
\\
Noise=$0\%$&$T_1$ real PD & CI & $p$-value& $T_1$ real PD & CI & $p$-value&  $T_1$ real PD & CI & $p$-value\\\hline
\\
\\
\\
\\
\textbf{$r=1$}& 0.059 &[0,0.112]&0.990  &0.038&[0, 0.058]&0.813&0.020&[0, 0.033]&0.940 \\

\\
\\
\\
\\

\textbf{$r=3$}& 0.176 &[0,0.303]&0.972  &0.115&[0,0.173]&0.780&0.061&[0,0.096]&0.926 \\
\\
\\
\\
\\
\\
\\
\\
\label{table:OneCirc_stat}
\end{tabular}
\end{center}
\footnotesize{Maximum statistic $T_1$ for the real $H_0$ persistence diagram (PD) and the simulated $H_0$ persistence diagrams of a sample $n$ of one circle with radius $r$. The confidence interval (CI) is a one-side confidence interval with $95\%$ confidence level. The $p$-value is also a one-side. Both the CI and the $p$-value are based on 1,000 simulated persistence diagrams.
}\\
\normalsize

\subsubsection{Comparison of persistence diagrams}
Given the fitted distribution, we can compare two persistence diagrams and conclude if they corresponded to the same data or not. Consider, for example, the above persistence diagram of one circle with $r=1$ based on a sample of $n=1,000$ points, and another, different, sample from the same object. As was mentioned above, the distribution of the $H_0$ points of the first sample is beta with [a,b]=[1.032,164.246]. Similarly, the distribution of the second sample is beta with [a,b]=[1.056,167.978]. The Kolmogorov-Smirnov statistic does not reject the null hypothesis of equal distributions, with p-value = 0.105. That is, these two persistence diagrams come form the same data, as was desired to get.

\subsection{Two Concentric Circles}
\subsubsection{Description and the distribution}
\begin{figure}[h!]
\bc
\includegraphics[width=1.45in, height=1.45in]{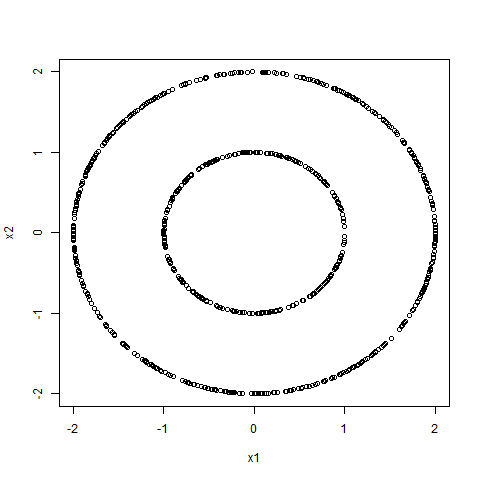} \hskip0.01truein
\includegraphics[width=1.45in, height=1.45in]{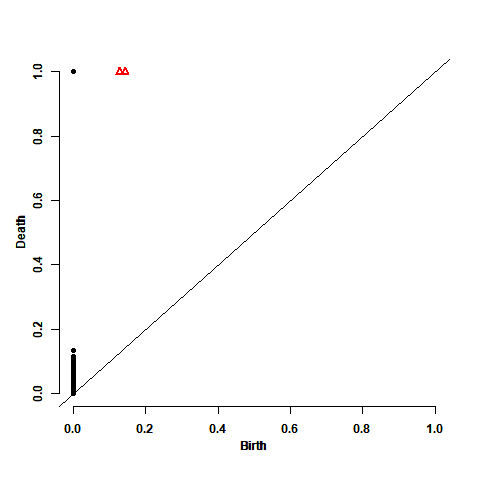} \hskip0.01truein
\\

\includegraphics[width=1.45in, height=1.45in]{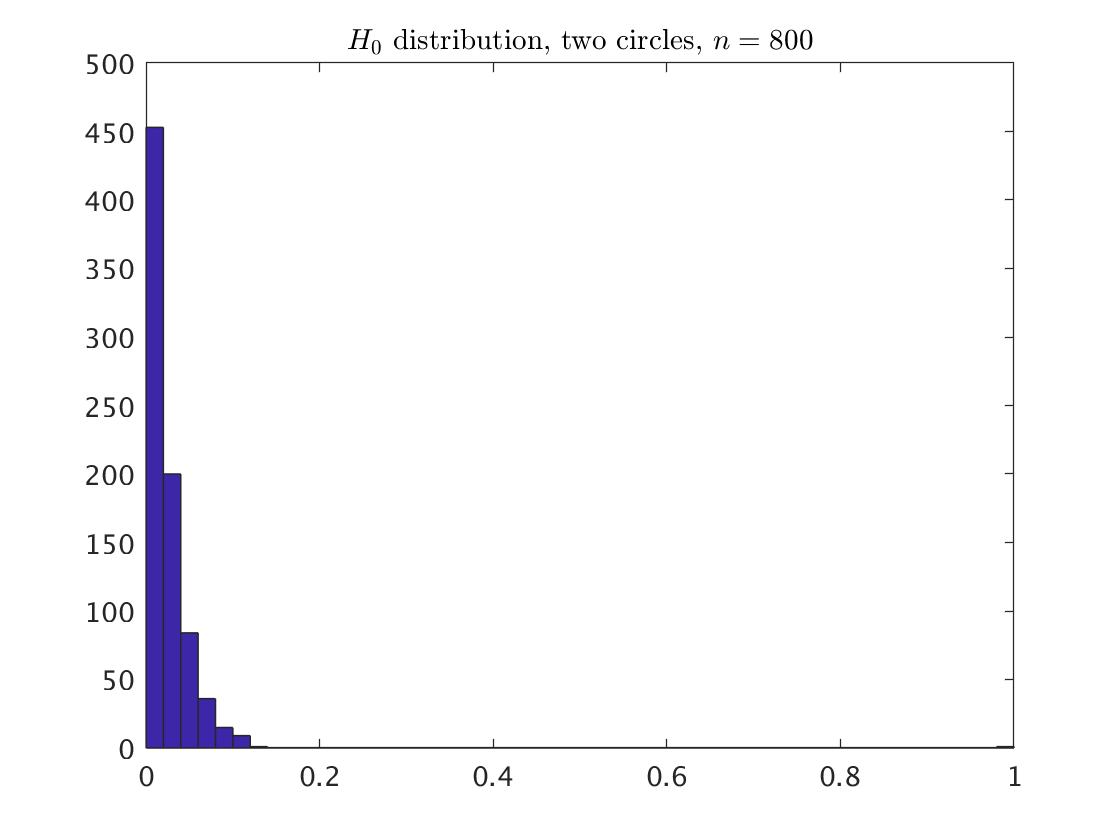} \hskip0.01truein
\includegraphics[width=1.45in, height=1.45in]{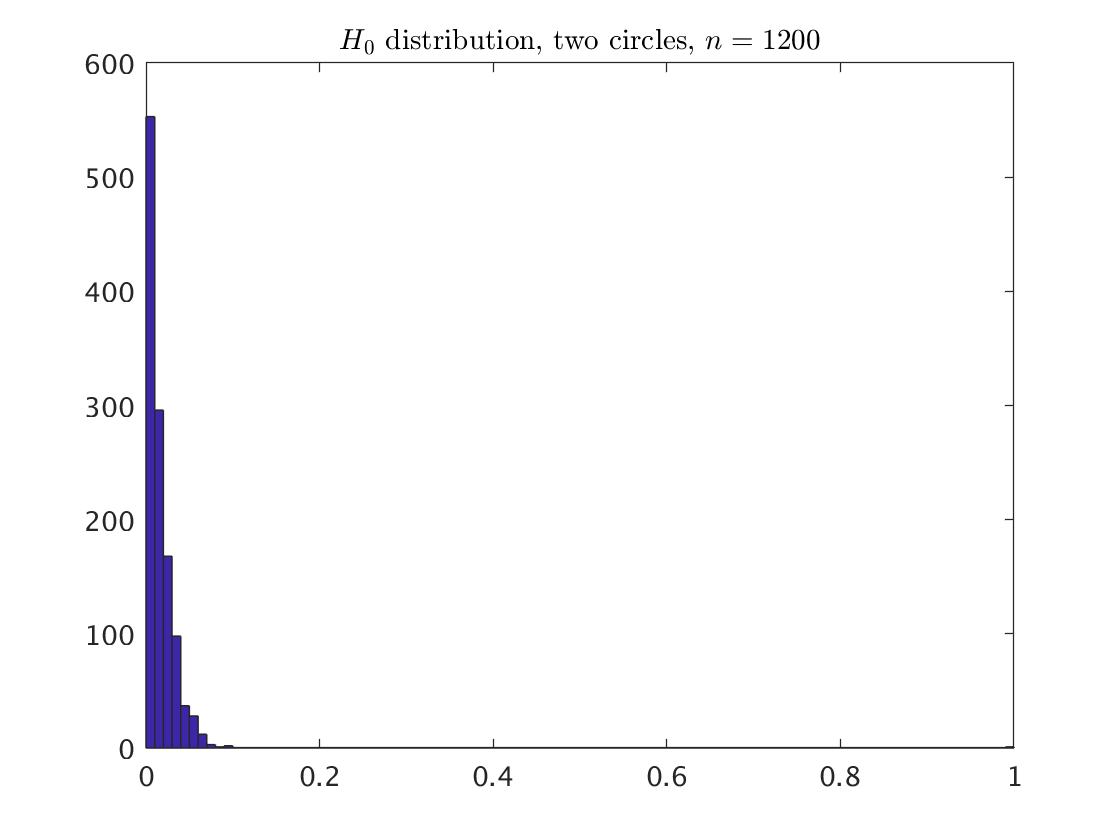} \hskip0.01truein
\includegraphics[width=1.45in, height=1.45in]{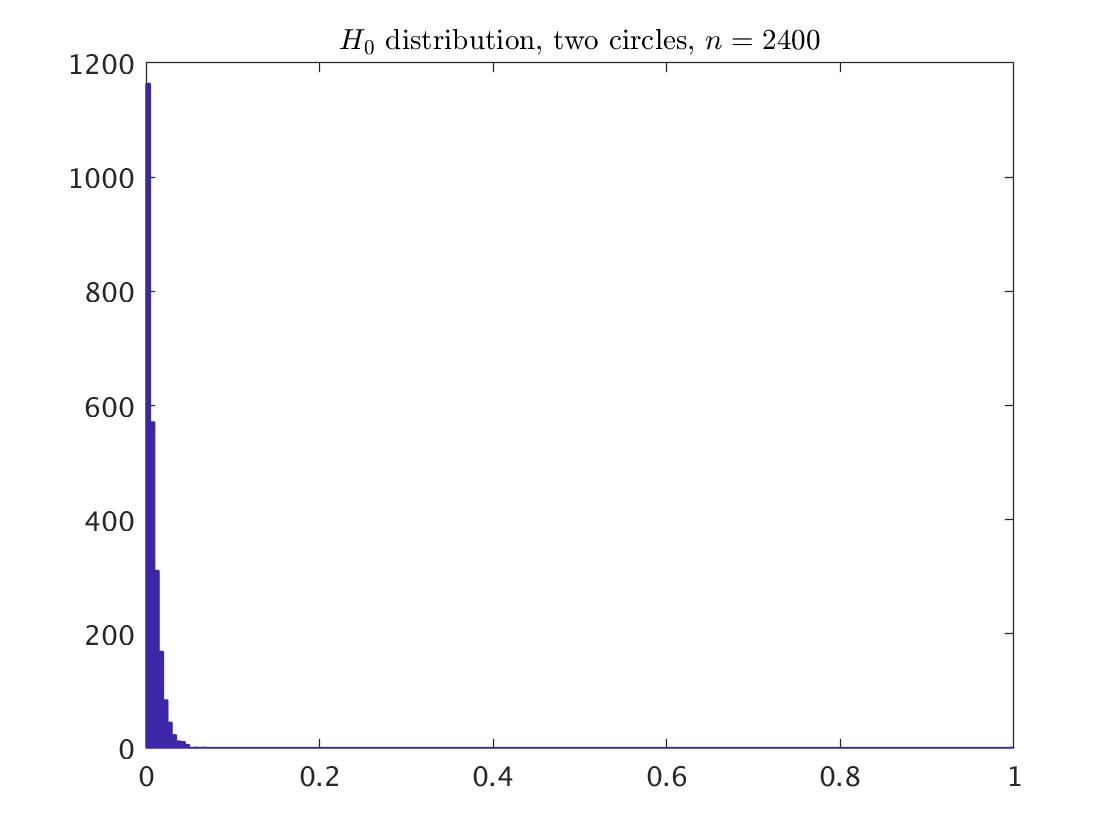} \hskip0.01truein

\ec
\caption{\footnotesize
 The plots describe two concentric circles. From left ro right: plot of the data, the corresponded persistence diagram based on the Vietoris-Rips filtration, the corresponded $H_0$ distribution without the point at infinity, and two more histograms that describe the $H_0$ distribution (without the point at infinity) for larger samples.
In the plot of the persistence diagram, the black circles are the $H_0$ persistence points, while the red triangles are the $H_1$ points.
}
\label{fig:Twocircle_M0}
\end{figure}

The second example describes two concentric circles. The sample includes $n=800$ points from two circles of diameters 4 and 2, where 500 points were chosen from the larger circle, and 300 from the smaller one. This sample is described in the left of Fig.\ \ref{fig:Twocircle_M0}.
Maxscale of 1 is the best value to cover the two concentric circles of this example. Smaller values than 1 obtain different $H_0$ persistence diagrams, whereas larger values than 1 obtain the same $H_0$ persistence diagram. The second plot of Fig.\ \ref{fig:Twocircle_M0} describes the corresponding persistence diagram of the Vietoris-Rips filtration using maxscale=1. We expect to see two black circles and two red triangles somewhat isolated from the other points in the diagram, and this is in fact the case. Next to that plot we have the corresponding histogram of the death times of the $H_0$ points without the point at infinity. The next histograms describe the $H_0$ points distributions that are corresponded to the larger samples $n=1,200$ and $n=2,400$ (using again maxscale=1), respectively, when the ratio of the points number from the smaller circle relative to that of the larger circle is continue to be 0.6.

\normalsize
\subsubsection{Statistical Inference}
In the following analysis, we distinguish between persistence diagrams that are based on maxsccale equals to 0.3, 0.5, and 1. Fitting the one dimensional distribution for the $H_0$ points, the best fitting is the beta distribution, where the relevant parameters are summarized in Appendix A.1.

For the identification of the topological signals, we again as in the previous example, generated 1,000 replicated $H_0$ persistence diagrams. Each of these persistence diagrams contains $(n-1)$ points corresponding to data with $n$ points. We calculated the two first maximum statistics of the death times $T_1$, and $T_2$, their confidence intervals and their p-values. Table \ref{table:TwoCircles_pv} summarizes the results. We get that $T_1$ is significant, and $T_2$ is insignificant. That is, there is one significant $H_0$ point among the set of the $H_0$ points without the point at infinity, and with adding back the point at infinity yield that there are two $H_0$ topological signals, as we hoped to find. \noindent As a result, also here where there are two connected components (comparing to the previous example that included one connected component), we have that the fitted distribution and the statistics perform well. Note that the result is independent on the maxscale value.
\begin{center}
\fontsize{8.5}{0.9}\selectfont
\captionof{table}{Two concentric circles - confidence interval and p-value}
\begin{tabular}{l|lccc|ccc}
\\
\\
\\
\\
Noise=0\%&$n$&$T_1$ real PD & CI & p-value  & $T_2$ real PD & CI & p-value\\\hline
\\
\\
\\
\\
maxscale=0.3&800& 0.300 &[0, 0.205]&0  &0.134&	[0, 0.170]&	0.612\\
\\
&1,200& 0.300&	[0, 0.153]&	0&	0.097&	[0, 0.125]&	0.661\\
\\
&2,400& 0.300&	[0, 0.084]&	0&	0.067&	[0, 0.069]&	0.085\\
\\
\\
\\
\\
\\
\\
maxscale=0.5&800& 0.500&	[0, 0.211]&	0&	0.134&	[0, 0.175]&	0.721\\
\\
&1,200& 0.500&	[0, 0.157]	&0&	0.097&	[0, 0.126]&	0.758\\
\\
&2,400& 0.500&	[0, 0.087]&	0&	0.067	&[0, 0.072]&	0.146\\
\\
\\
\\
\\
\\
\\
maxscale=1&800& 1&	[0, 0.727]&	0&	0.134&	[0, 0.631]&	1\\
\\
&1,200& 1&	[0, 0.594]&	0&	0.097&	[0, 0.505]&	1\\
\\
&2,400& 1&	[0, 0.390]&	0&	0.067&	[0, 0.326]&	1\\
\\
\\
\\
\\
\\
\\

\label{table:TwoCircles_pv}
\end{tabular}
\end{center}
\footnotesize{Maximum statistics $T_1$ and $T_2$, for the real $H_0$ persistence diagram and the simulated $H_0$ persistence diagrams of the two concentric circles example. The CI is a one-side confidence interval with $95\%$ confidence level. The p-value is also a one-side. Both the CI and the $p$-value are based on 1,000 simulated persistence diagrams.}\\
\normalsize

\subsection{Two Distinct Circles}
\subsubsection{Description and the Distribution}
\begin{figure}[h!]
\bc
\includegraphics[width=1.45in, height=1.45in]{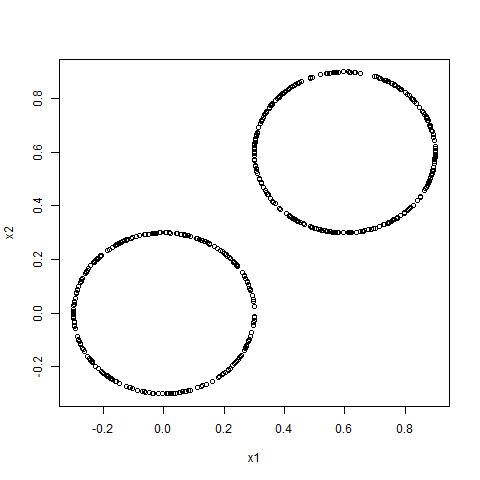} \hskip0.01truein
\includegraphics[width=1.45in, height=1.45in]{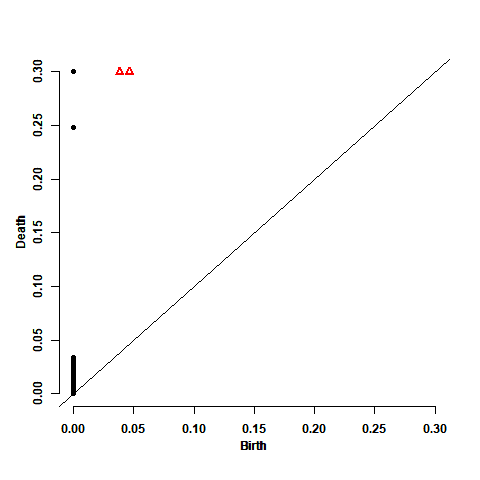} \hskip0.01truein
\includegraphics[width=1.45in, height=1.45in]{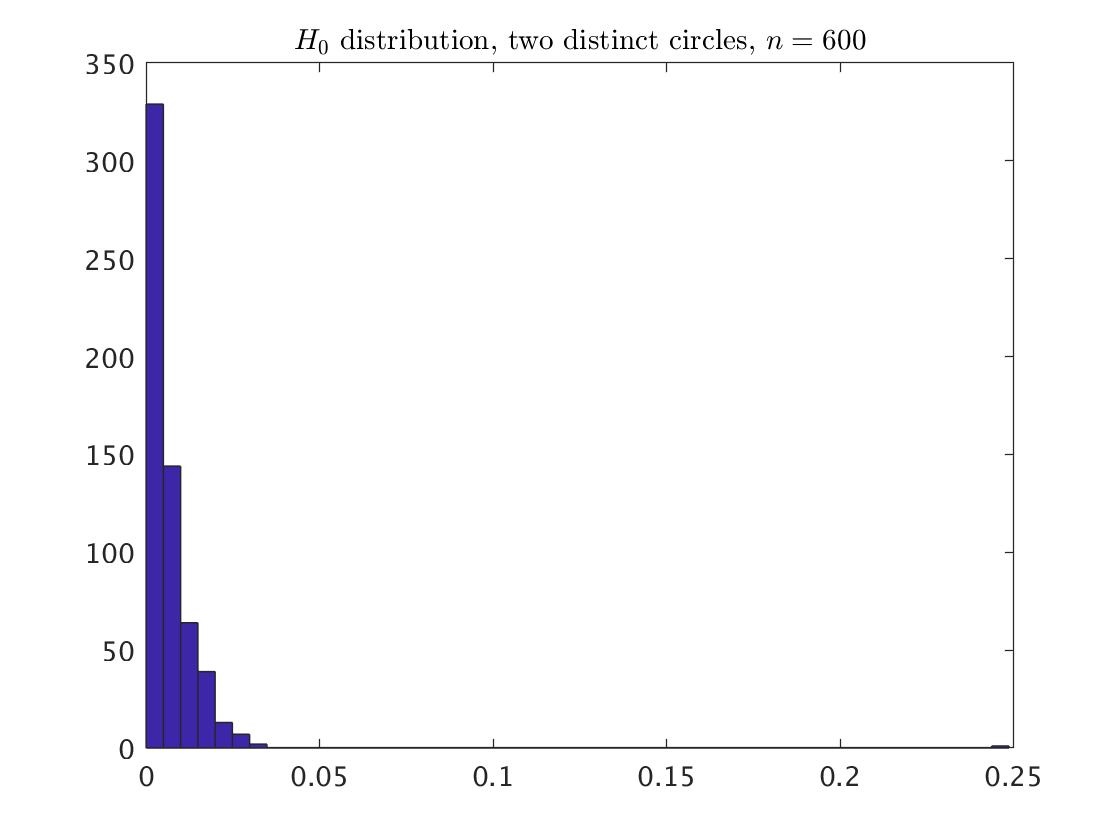} \hskip0.01truein
\includegraphics[width=1.45in, height=1.45in]{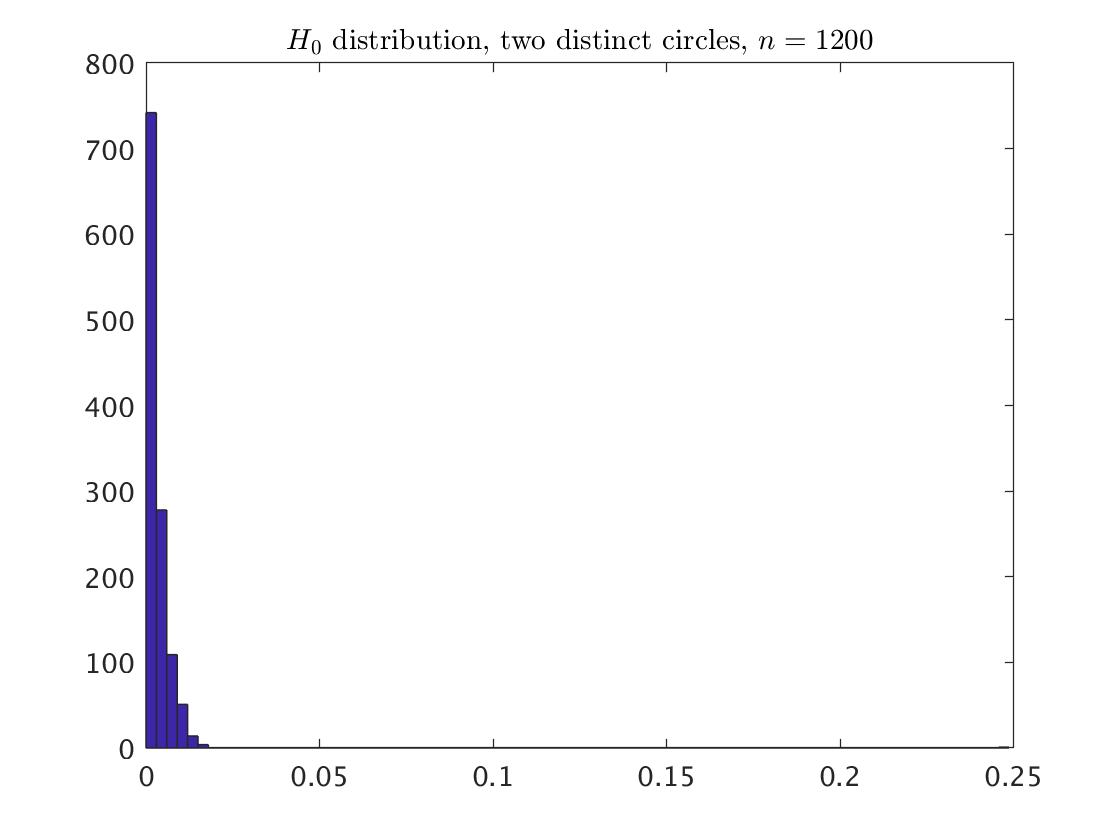} \hskip0.01truein
\ec
\caption{\footnotesize
 The plots describe two distinct circles. From left ro right: plot of the data, the corresponded persistence diagram based on the Vietoris-Rips filtration, the corresponded $H_0$ distribution without the point at infinity, and one more histogram that describes the $H_0$ distribution (without the point at infinity) for larger sample.
In the plot of the persistence diagram, the black circles are the $H_0$ persistence points, while the red triangles are the $H_1$ points.
}
\label{fig:Distcircle_M0}
\end{figure}

The third example includes two distinct circles. The sample has $n=600$ points from two distinct circles, each of them has radius $r=0.3$ and contains 300 points. The distance between the two circles is 0.6 for each point. This sample is described in the left of Fig.\ \ref{fig:Distcircle_M0}.
\normalsize

To its right, we see the corresponding persistence diagram of the Vietoris-Rips filtration using maxscale=0.3. We expect to see two black circles and two red triangles somewhat isolated from the other points in the diagram, and this is in fact the case. The next plot is the histogram of the corresponding $H_0$ points (without including the $H_0$ point at infinity). The last plot is the histogram of the $H_0$ points under larger sample of $n=1,200$ points, when the ratio between the points number of the two circles is continue to be 0.5, keeping the distance of 0.6 between the circles.
The maxscale of 0.3 is enough to cover the two distinct circles of this example, we get the same $H_0$ persistence diagrams for maxscale equals to or greater than 0.3.

Again the best fitting of the $H_0$ points is the beta distribution, where the relevant parameters are summarized in Appendix A.1.
\subsubsection{Statistical Inference}
Also here we generated 1,000 replicated $H_0$ persistence diagrams, and calculated the two first maximum statistics of the death times $T_1$, and $T_2$, their confidence intervals and their p-values. Table \ref{table:TwoFarCircles_pv} summarizes the results. We get that $T_1$ is significant, and $T_2$ is insignificant. That is, there is one significant $H_0$ point among the set of the $H_0$ points without the point at infinity, and with adding back the point at infinity yield that there are two $H_0$ topological signals, as we hoped to find.

\begin{center}
\fontsize{8.5}{0.9}\selectfont
\captionof{table}{Two distinct circles - confidence interval and p-value}
\begin{tabular}{l|lccc|ccc}
\\
\\
\\
\\
Noise=0\%&$n$&$T_1$ real PD & CI & p-value  & $T_2$ real PD & CI & p-value\\\hline
\\
\\
\\
\\
maxscale=0.3&600& 0.249&	[0, 0.063]&	0&	0.034&	[0, 0.051]&	0.939\\
\\
&1,200& 0.249&	[0, 0.034]	&0	&0.017&	[0, 0.028]&	0.998\\
\\
\\
\\
\\
\\
\\
\label{table:TwoFarCircles_pv}
\end{tabular}
\end{center}
\footnotesize{Maximum statistics $T_1$ and $T_2$, for the real $H_0$ persistence diagram and the simulated $H_0$ persistence diagrams of the two distinct circles example. The CI is a one-side confidence interval with $95\%$ confidence level. The p-value is also a one-side. Both the CI and the $p$-value are based on 1,000 simulated persistence diagrams.}\\
\normalsize

\subsection{The two dimensional sphere}
\subsubsection{Description and the Distribution}
\begin{figure}[h!]
\bc
\includegraphics[width=1.45in, height=1.45in]{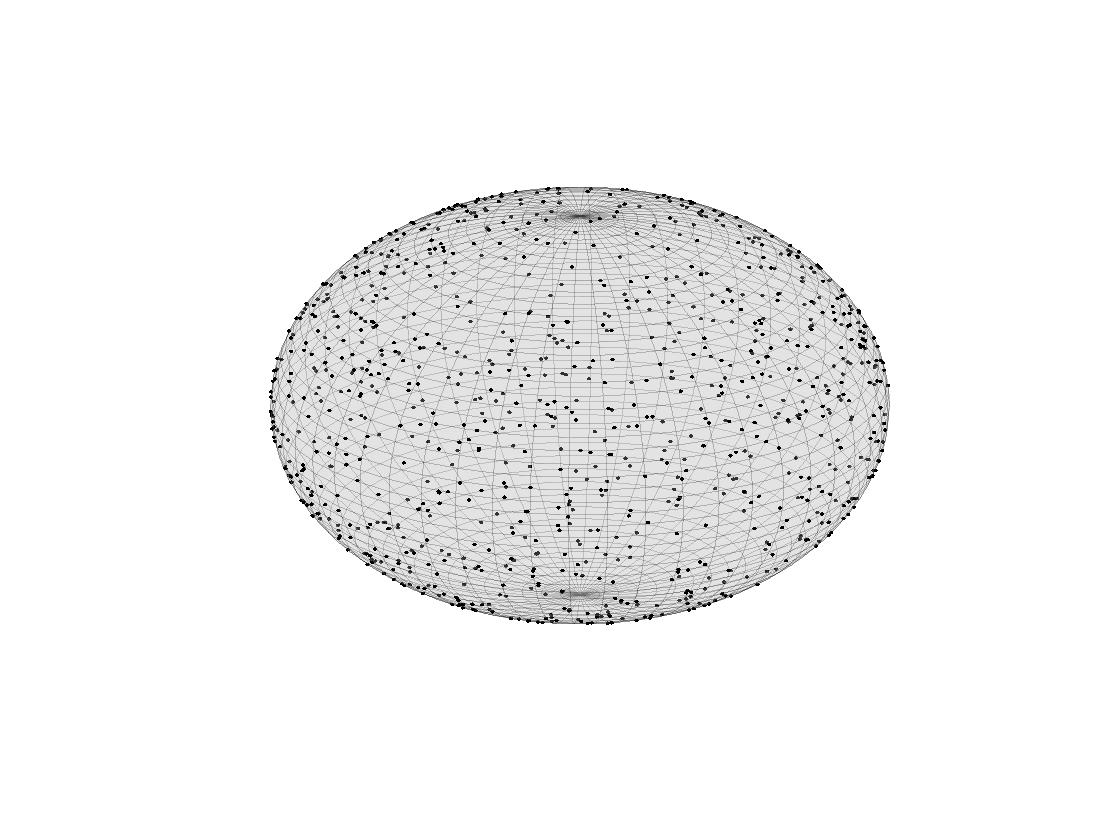} \hskip0.01truein
\includegraphics[width=1.45in, height=1.45in]{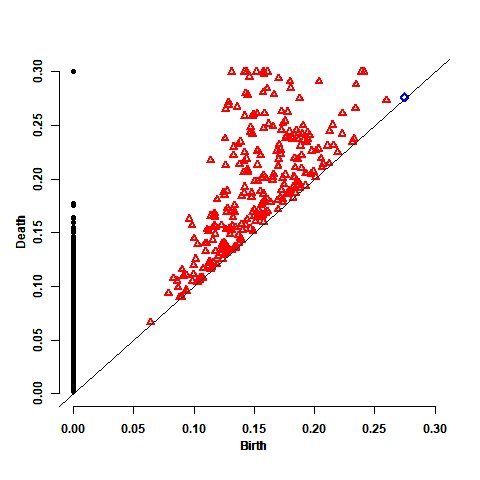} \hskip0.01truein
\includegraphics[width=1.45in, height=1.45in]{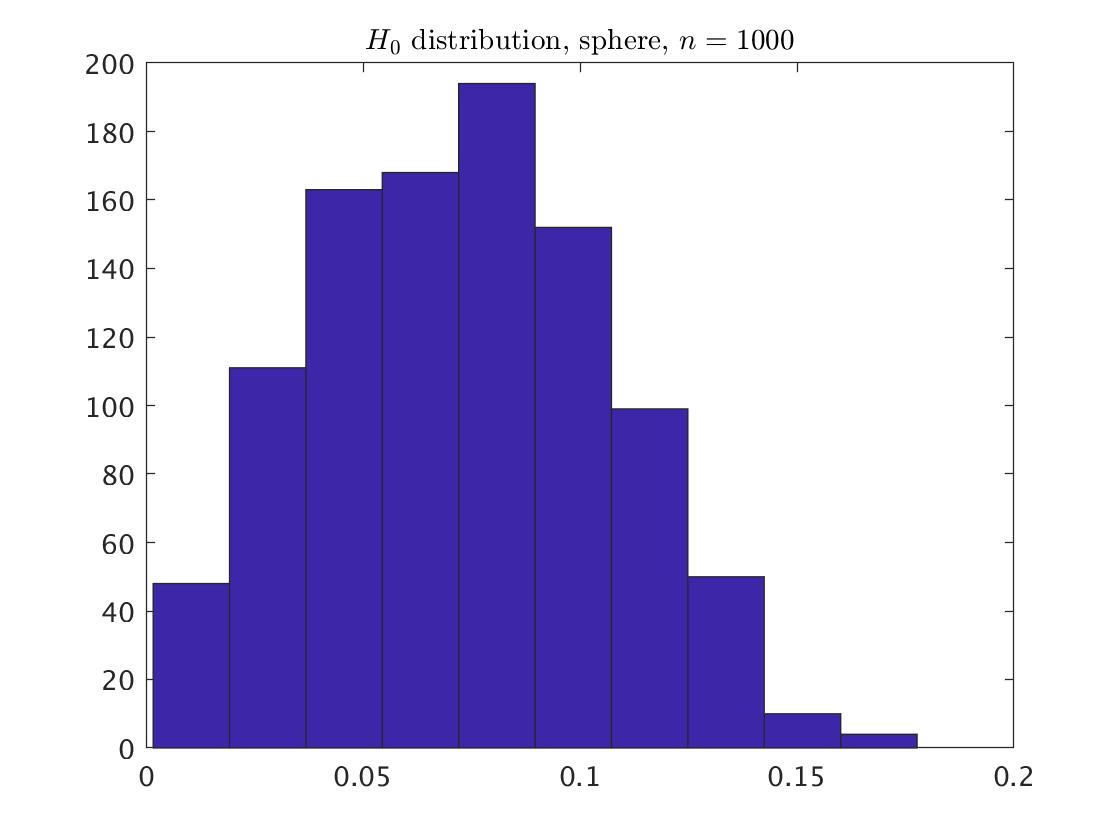} \hskip0.01truein
\includegraphics[width=1.45in, height=1.45in]{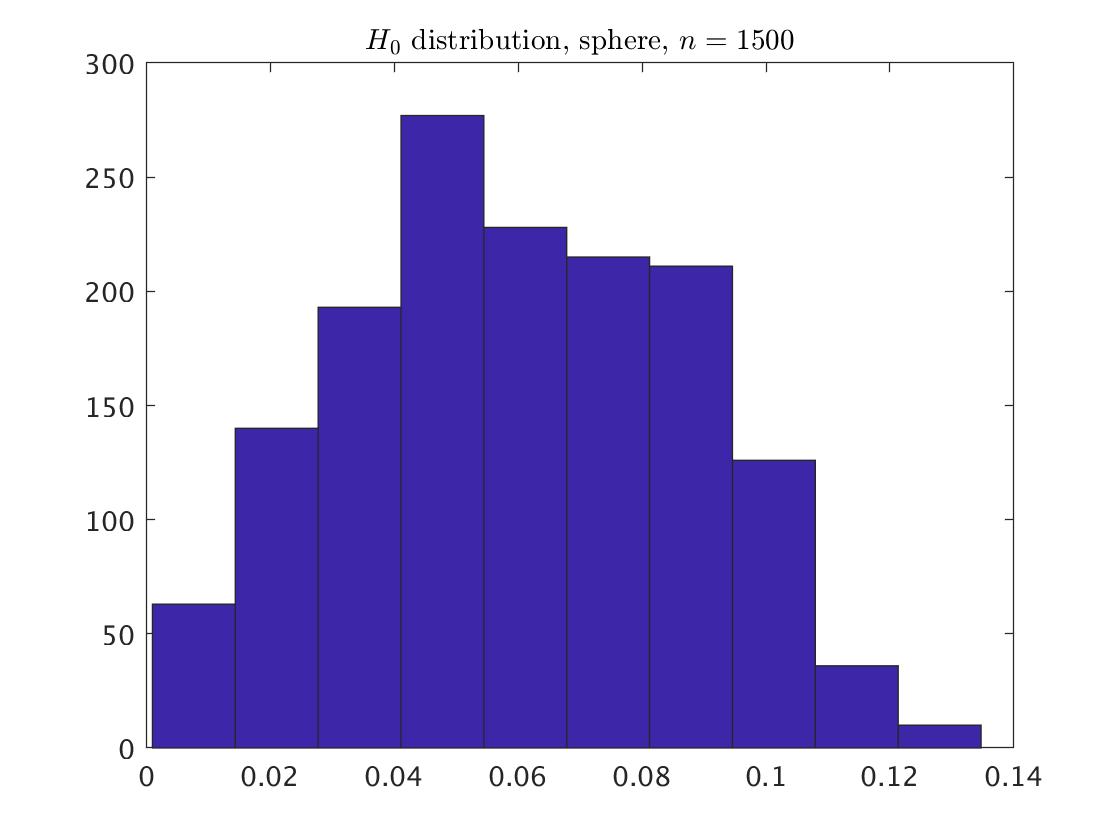} \hskip0.01truein
\ec
\caption{\footnotesize
The plots describe the 2-sphere. From left to right: plot of the data, the corresponded persistence diagram based on the Vietoris-Rips filtration, the corresponded $H_0$ distribution without the point at infinity, and one more histogram that describes the $H_0$ distribution (without the point at infinity) for larger sample.
In the plot of the persistence diagram, the black circles are the $H_0$ persistence points, the red triangles are the $H_1$ points, and the blue diamond is the $H_2$ point.
}
\label{fig:sphere_M0}
\end{figure}

Whereas the three previous examples included two dimensional objects, the current and the next example include higher dimensional objects. The current example is  a random sample of $n=1,000$ points from the uniform distribution on the sphere $S^2$ in $R^3$ with radius $r=1$. This sample is described in the left of Fig.\ \ref{fig:sphere_M0}.
\normalsize

 The maxscale of 0.3 is enough to cover the 2-sphere of this example, and the persistence diagram based on the Vietoris-Rips filtration with this maxscale is presented in the second plot of Fig.\ \ref{fig:sphere_M0}. To its right presented are the histograms of the corresponding $H_0$ points, and again, are not including the $H_0$ point at infinity. They are corresponded to random samples of $n=1,000$ and $n=1,500$ points.
The 2-sphere is characterized by having a single connected component and a single void. Therefore we expect to have one black circle somewhat isolated from the other points in the diagram and one blue diamond. The void does not have to be isolated from the other points due to the short lifetimes of high dimensional homologies. This is in fact the case.

The best fitting of the $H_0$ points is the beta distribution. The relevant parameters are summarized in Appendix A.1.
\subsubsection{Statistical Inference}
For the identification of the topological signals for each of the considered values of $n$, we generated 1,000 replicated $H_0$ persistence diagrams, each one contained $(n-1)$ points. We calculated the maximum statistic of the death times $T_1$, it confidence interval and its p-value. Table \ref{table:sphere_pv} summarizes the results. We get that $T_1$ is insignificant, and together with adding back the point at infinity yield that there is one $H_0$ topological signal, as we hoped to find.

\begin{center}
\fontsize{8.5}{0.9}\selectfont
\captionof{table}{Two dimensional sphere - confidence interval and p-value}
\begin{tabular}{l|lccc}
\\
\\
\\
\\
Noise=0\%&$n$&$T_1$ real PD & CI & p-value\\\hline
\\
\\
\\
\\
maxscale=0.3&1,000& 0.178&	[0, 0.300]&	1\\
\\
&1,500& 0.135&	[0, 0.250]&	1\\

\label{table:sphere_pv}
\end{tabular}
\end{center}
\footnotesize{Maximum statistic $T_1$ for the real $H_0$ persistence diagram and the simulated $H_0$ persistence diagrams of the 2-sphere. The CI is a one-side confidence interval with $95\%$ confidence level. The p-value is also a one-side. Both the CI and the $p$-value are based on 1,000 simulated persistence diagrams. }\\
\normalsize

\subsection{3-Torus}
\subsubsection{Description and the Distribution}
\begin{figure}[h!]
\bc
\includegraphics[width=1.45in, height=1.45in]{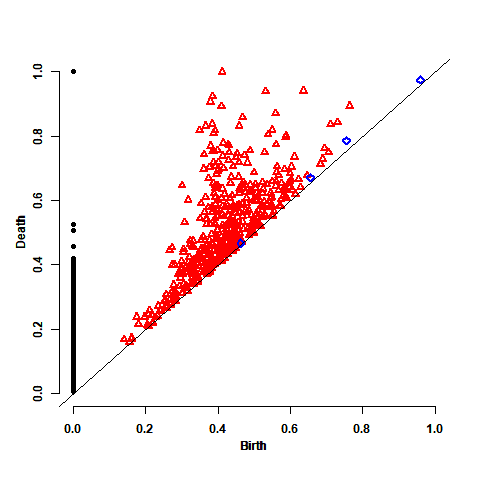} \hskip0.01truein
\includegraphics[width=1.45in, height=1.45in]{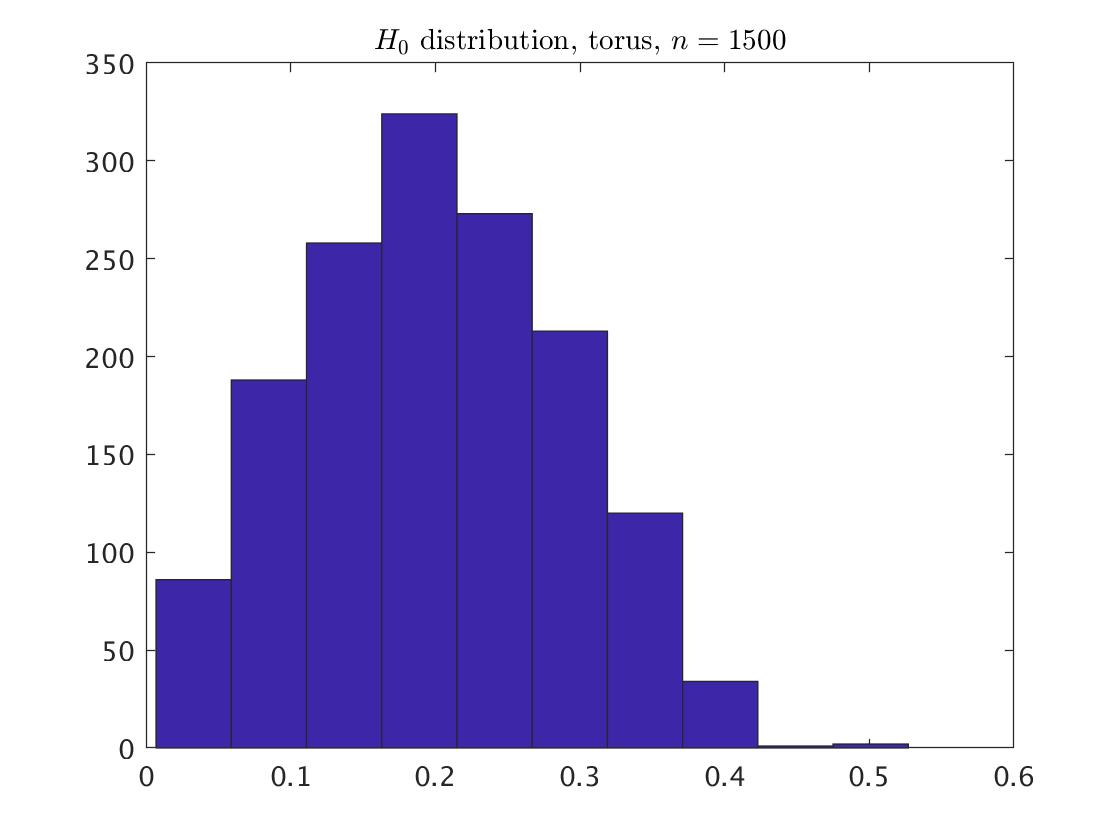} \hskip0.01truein
\includegraphics[width=1.45in, height=1.45in]{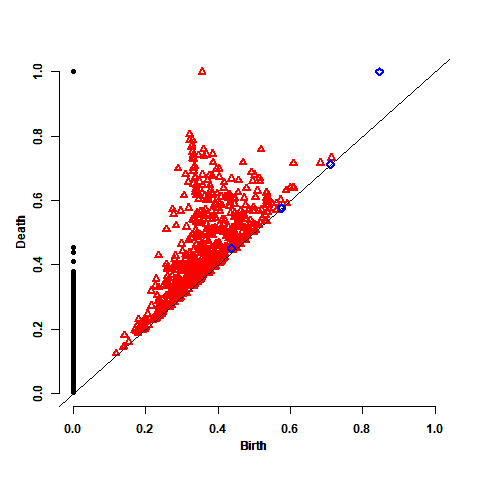} \hskip0.01truein
\includegraphics[width=1.45in, height=1.45in]{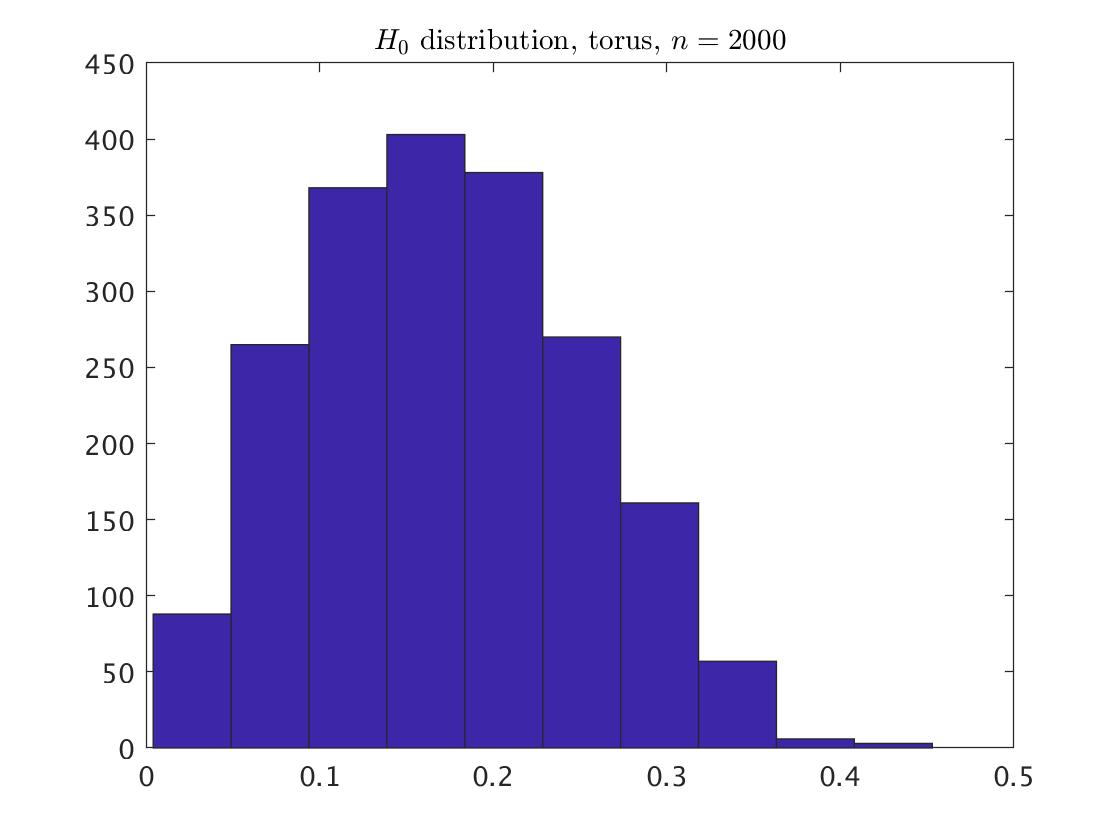} \hskip0.01truein
\ec
\caption{\footnotesize
The plots describe the 3-torus. From left to right: a persistence diagram based on the Vietoris-Rips filtration for a sample of $n=1,500$, the corresponded $H_0$ distribution without the point at infinity, a persistence diagram based on the Vietoris-Rips filtration for a sample of $n=2,000$, and the corresponded $H_0$ distribution without the point at infinity. In the plots of the persistence diagram, the black circles are the $H_0$ persistence points, the red triangles are the $H_1$ points, and the blue diamond are the $H_2$ points.
}
\label{fig:Torus_M0}
\end{figure}

Here we take a sample of $n=1,500$ points from the 3-torus $T^3$, chosen uniformly with respect to the natural Riemannian metric induced on it as a subset on $\mathbb R^4$. In this example, maxscale of 1 is enough to cover the whole 3-torus.
Since $T^3$ lives in $\real^4$, we cannot show the picture of the sample. However, the persistence diagram and the $H_0$ persistence diagram distribution are just as easy to see here as they were before, and they are shown in the two first plots Fig.\ \ref{fig:Torus_M0}, based on maxscale=1. In the persistence diagram we expect to see one black circle, three red triangles, and one blue diamond, somewhat isolated from the other points in the diagram. We can see the one black circle but for the three triangles we need a larger sample to recognize them. For example, the persistence diagram that is based on a sample of $n=2,000$ (based on maxscale=0.5), as described in the third plot of Fig.\ \ref{fig:Torus_M0}, can now recognize one of the three red isolated triangles. To its right we have the histogram corresponded to the $H_0$ points under the larger sample of $n=2,000$.

\normalsize
Fitting the one dimensional distribution for the $H_0$ points, the best fitting for $n=1,500$ is the beta distribution under maxscale=0.5, 1, but the generalized Pareto distribution for maxscale=0.3. For $n=2,000$, the best fitting is the generalized Pareto distribution for maxscale=0.3 and the beta distribution for maxscale=0.5. The relevant parameters are summarized in Appendix A.1.
%
\subsubsection{Statistical Inference}
For the identification of the topological signals for each of the considered values of $n$, we generated 1,000 replicated $H_0$ persistence diagrams, each one contained $(n-1)$ points. We calculated the maximum statistic of the death times $T_1$, it confidence interval and its p-value. Table \ref{table:Torus_pv} summarizes the results. We get that $T_1$ is insignificant among the $H_0$ without the point at infinity. Adding back the point at infinity yields that there is one $H_0$ topological signal, as we hoped to find. Note that this result is independent on the maxscale value.

\begin{center}
\fontsize{8.5}{0.9}\selectfont
\captionof{table}{Torus $T^3$ - confidence interval and p-value}
\begin{tabular}{l|lccc}
\\
\\
\\
\\
Noise=0\%&$n$&$T_1$ real PD & CI & p-value\\\hline
\\
\\
\\
\\
maxscale=0.3&1,500& 0.300&	[0,0.300]&	0.926\\
\\
&2,000& 0.300&	[0,0.300]&	0\\
\\
\\
\\
\\
\\
\\
\\
maxscale=0.5&1,500& 0.500&	[0,0.672]&	1\\
\\
&2,000& 0.453&	[0,0.615]&	0.999\\
\\
\\
\\
\\
\\
\\
\\
maxscale=1&1500& 0.527&	[0,0.672]&	0.983\\

\label{table:Torus_pv}
\end{tabular}
\end{center}
\footnotesize{Maximum statistic $T_1$ for the real $H_0$ persistence diagram and the simulated $H_0$ persistence diagrams of the 3-torus. The CI is a one-side confidence interval with $95\%$ confidence level. The p-value is also a one-side. Both the CI and the $p$-value are based on 1,000 simulated persistence diagrams. }\\
\normalsize

\section{Examples of Noisy Data}
In this section we examine the influence of adding noise to some fraction $M$ of the $n$ data points. We check the settings of  $M$ = 30\%, 70\%, 80\% and 100\%. The additive noise is generated from the bivariate normal with zero mean, and the identity matrix divided by 9 as the covariance. The examples that we consider here are the one circle (with $r=1$ and $r=3$), the two concentric circles, and the 3-torus of Section 3, but now with adding a noise $M$.
\subsection{One Circle}
\begin{landscape}
\begin{figure}[h!]
\bc
\includegraphics[width=1.45in, height=1.45in]{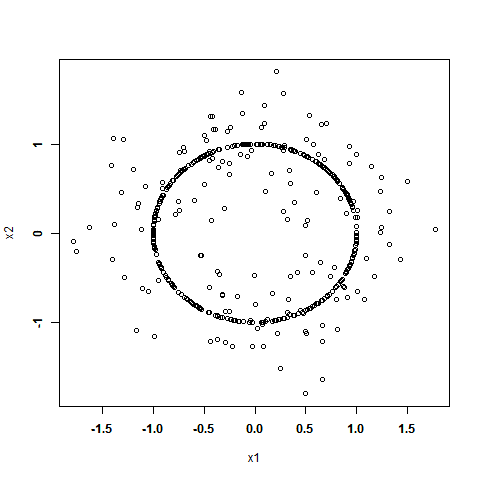} \hskip0.01truein
\includegraphics[width=1.45in, height=1.45in]{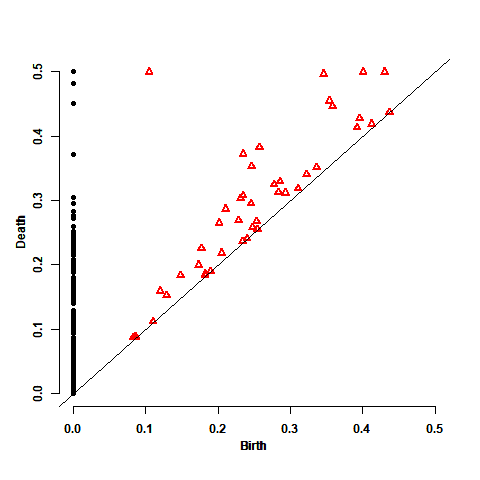} \hskip0.01truein
\includegraphics[width=1.45in, height=1.45in]{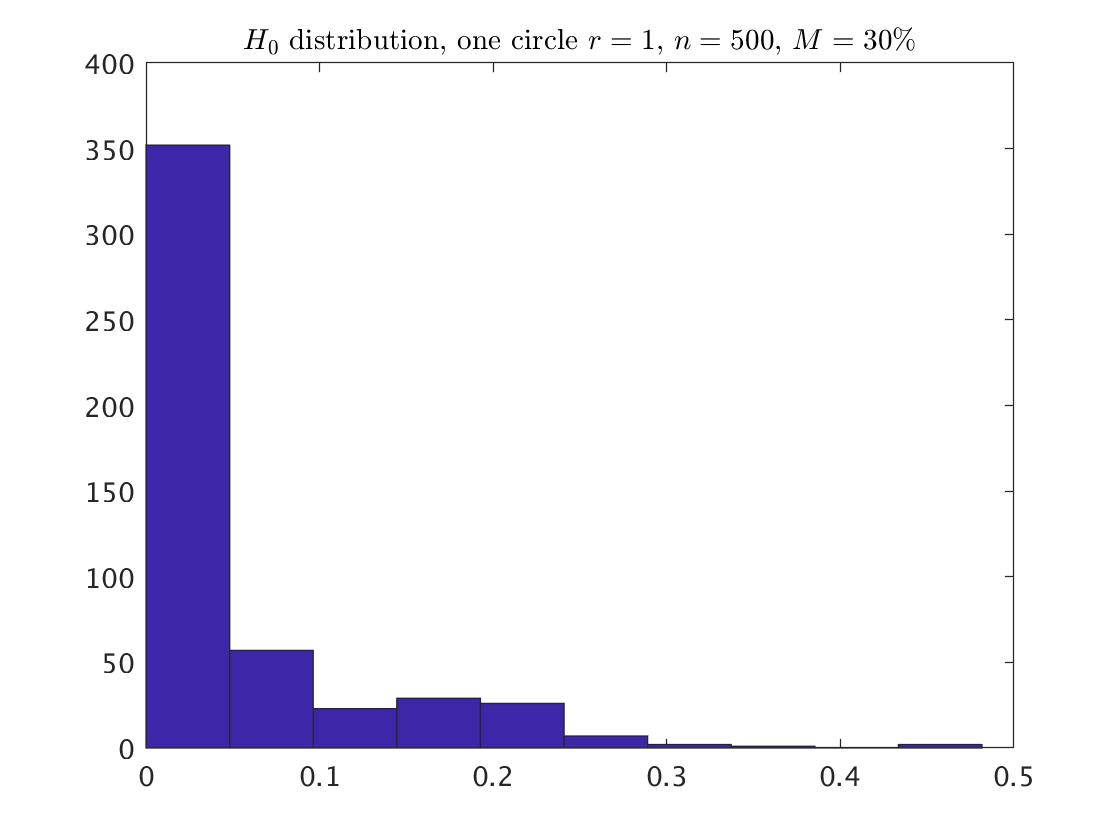} \hskip0.01truein
\includegraphics[width=1.45in, height=1.45in]{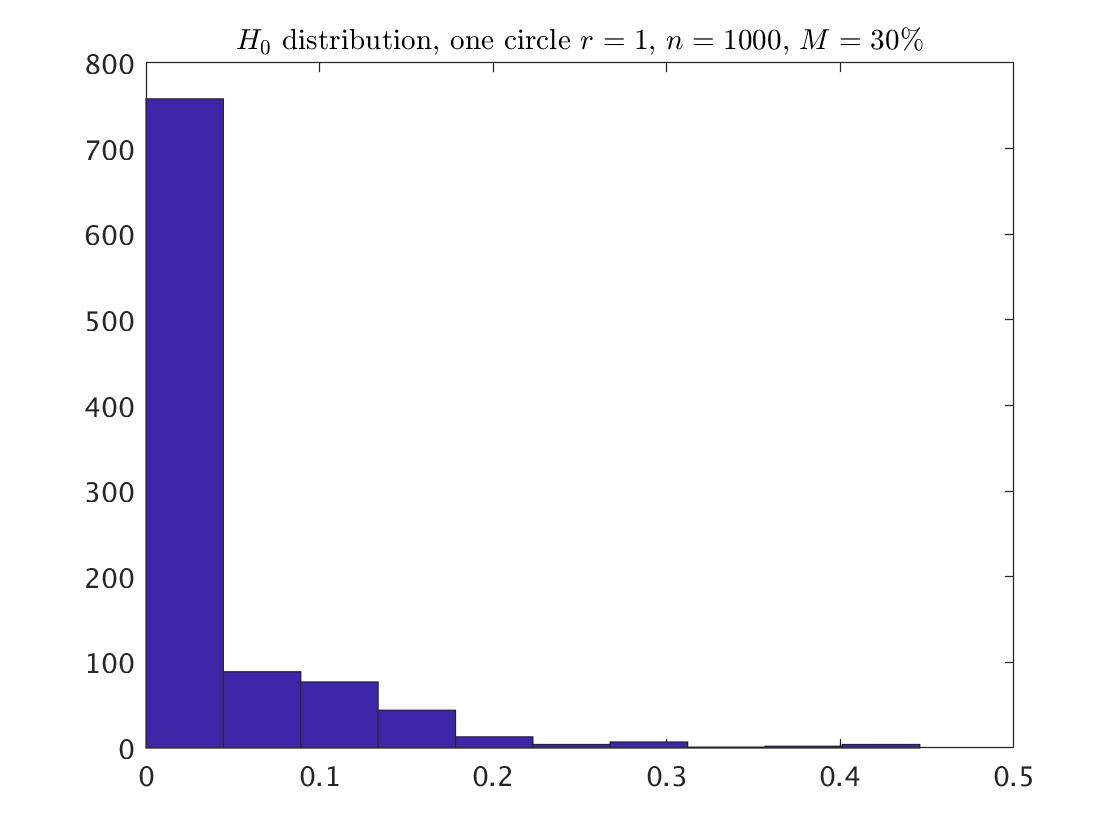} \hskip0.01truein
\includegraphics[width=1.45in, height=1.45in]{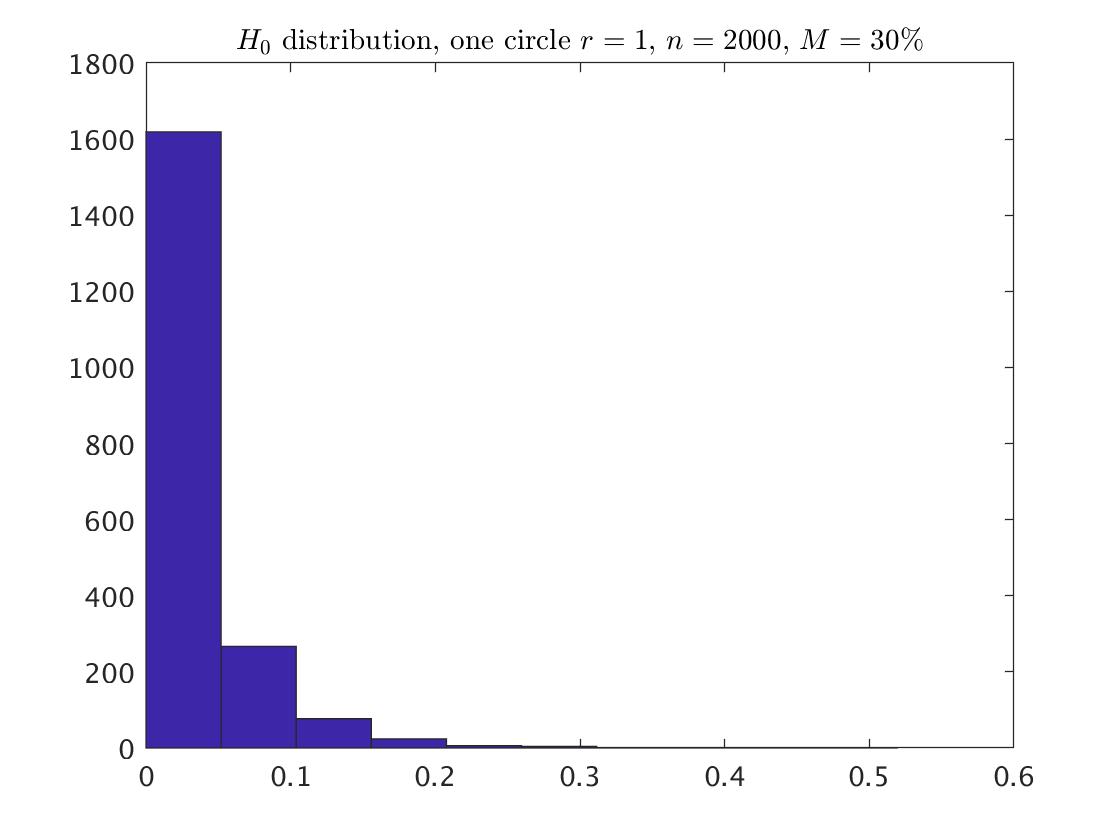} \hskip0.01truein

\includegraphics[width=1.45in, height=1.45in]{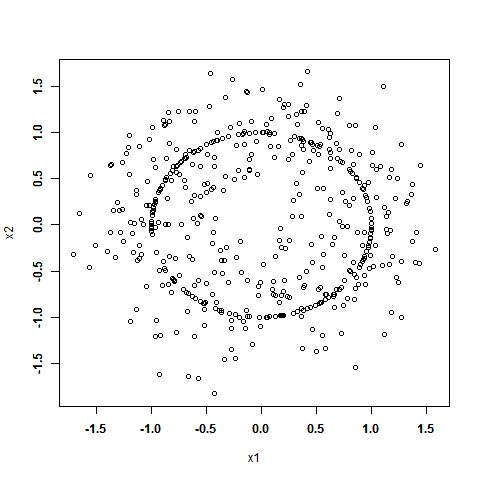} \hskip0.01truein
\includegraphics[width=1.45in, height=1.45in]{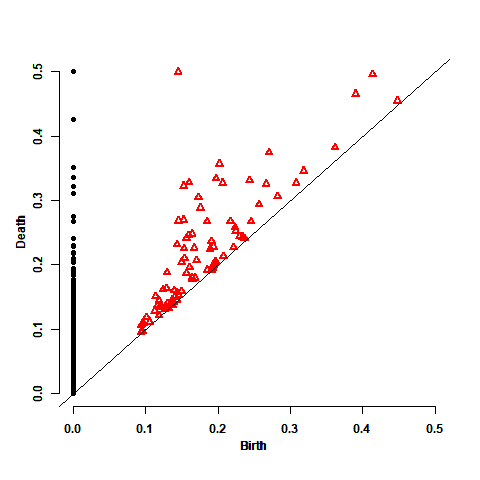} \hskip0.01truein
\includegraphics[width=1.45in, height=1.45in]{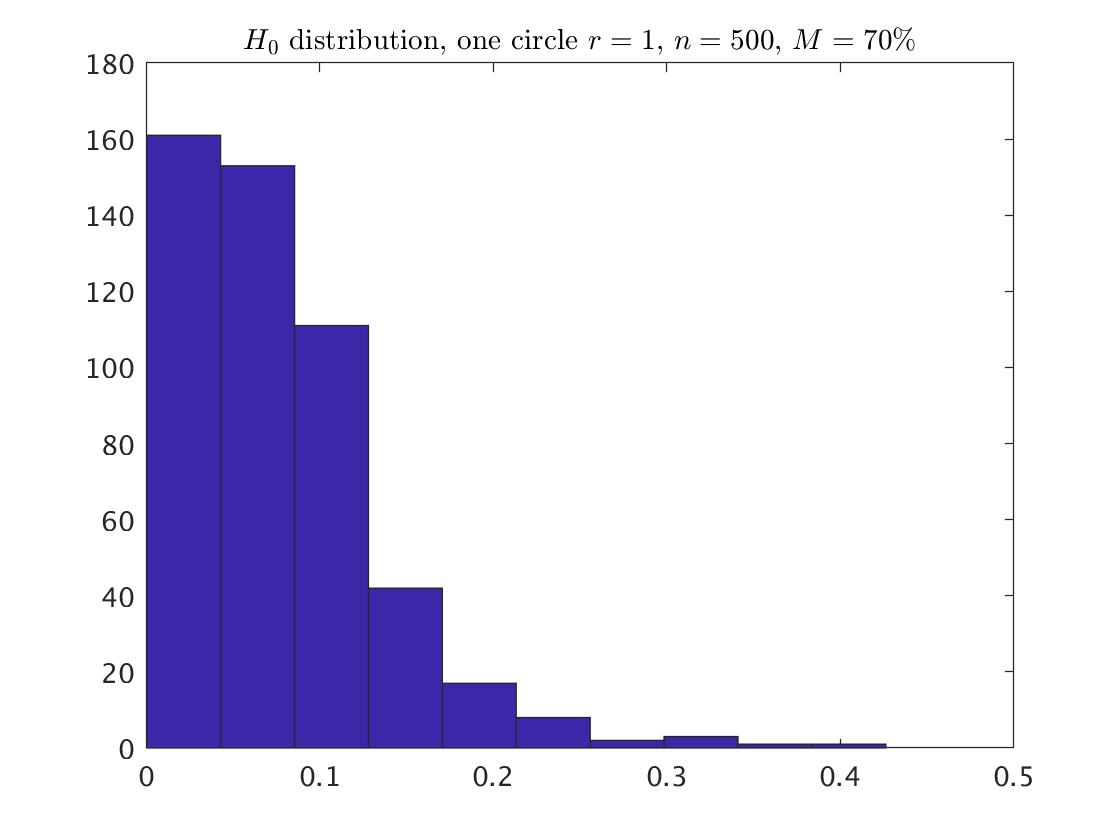} \hskip0.01truein
\includegraphics[width=1.45in, height=1.45in]{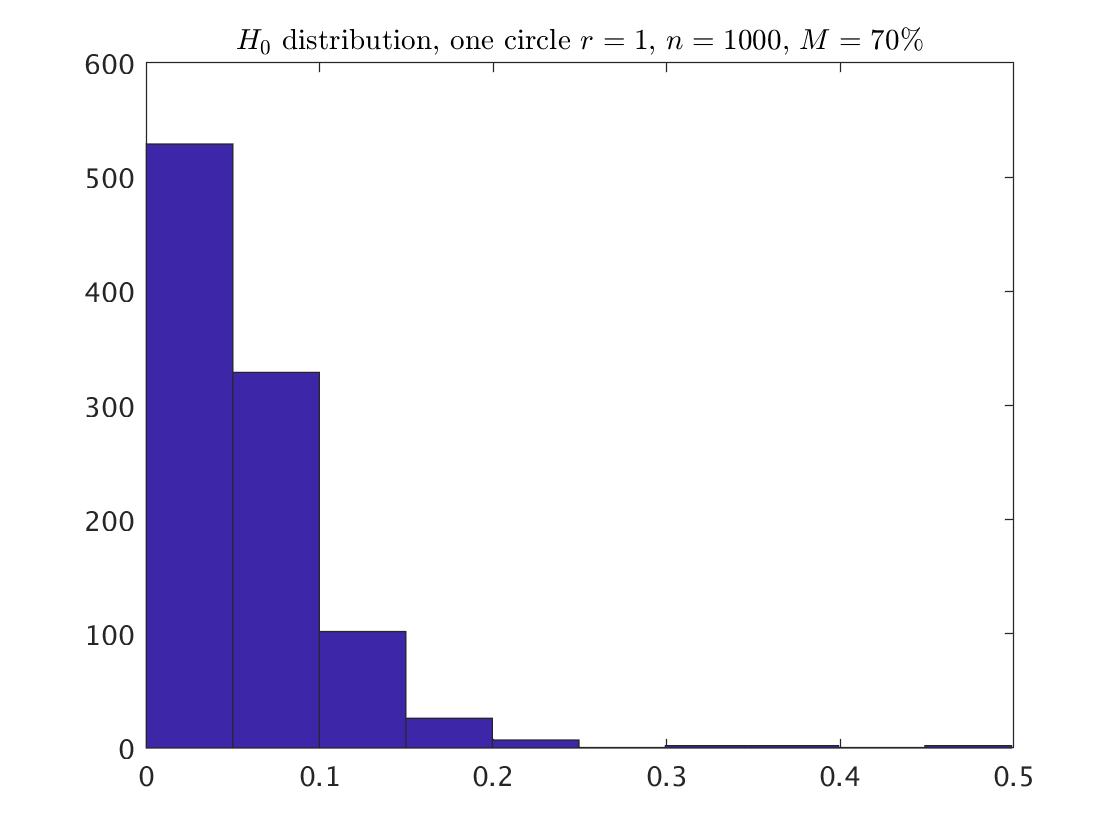} \hskip0.01truein
\includegraphics[width=1.45in, height=1.45in]{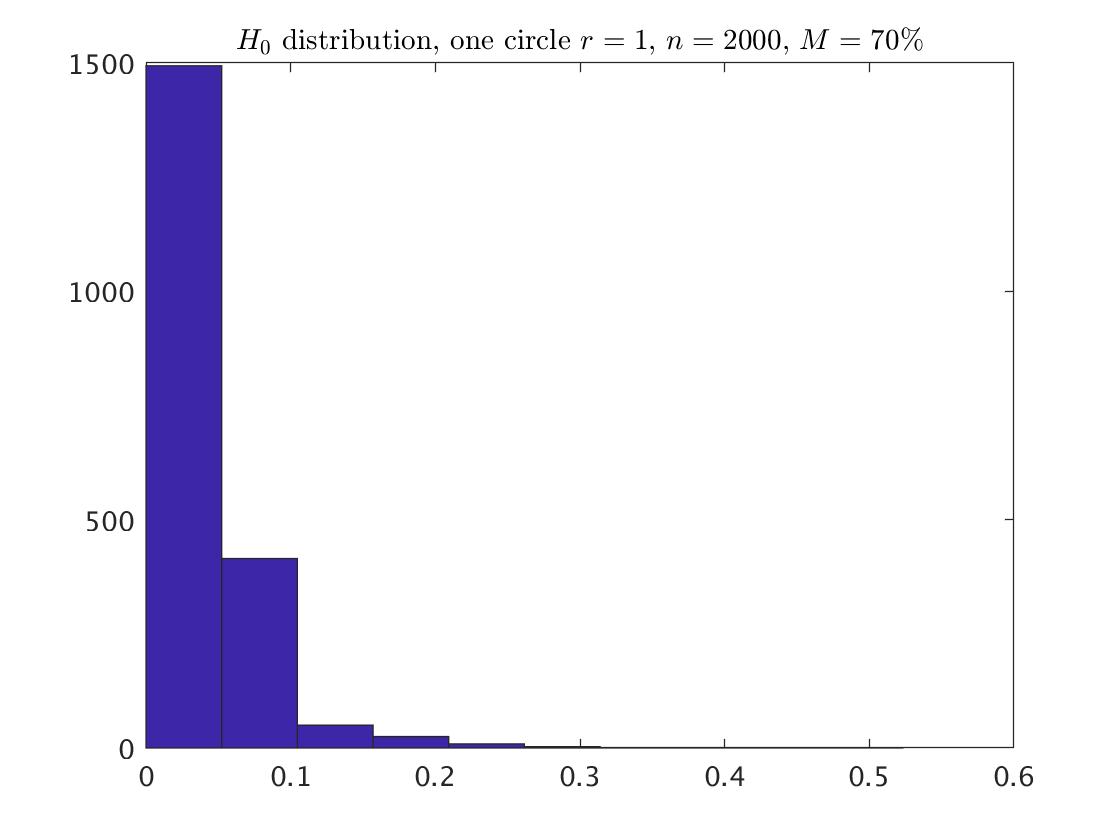} \hskip0.01truein

\includegraphics[width=1.45in, height=1.45in]{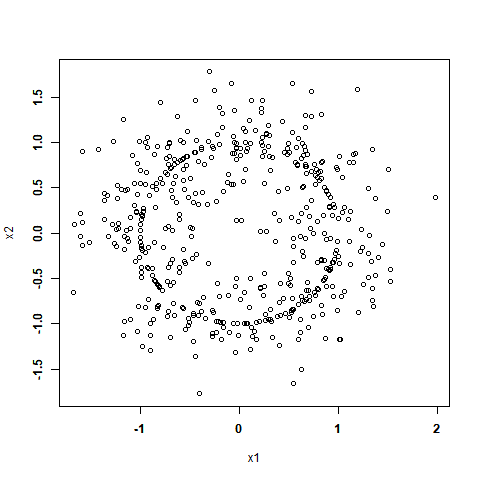} \hskip0.01truein
\includegraphics[width=1.45in, height=1.45in]{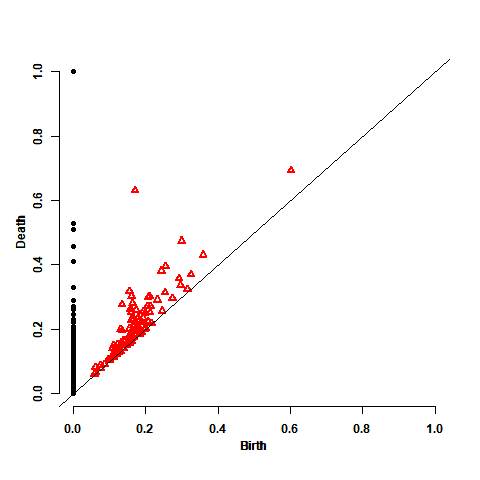} \hskip0.01truein
\includegraphics[width=1.45in, height=1.45in]{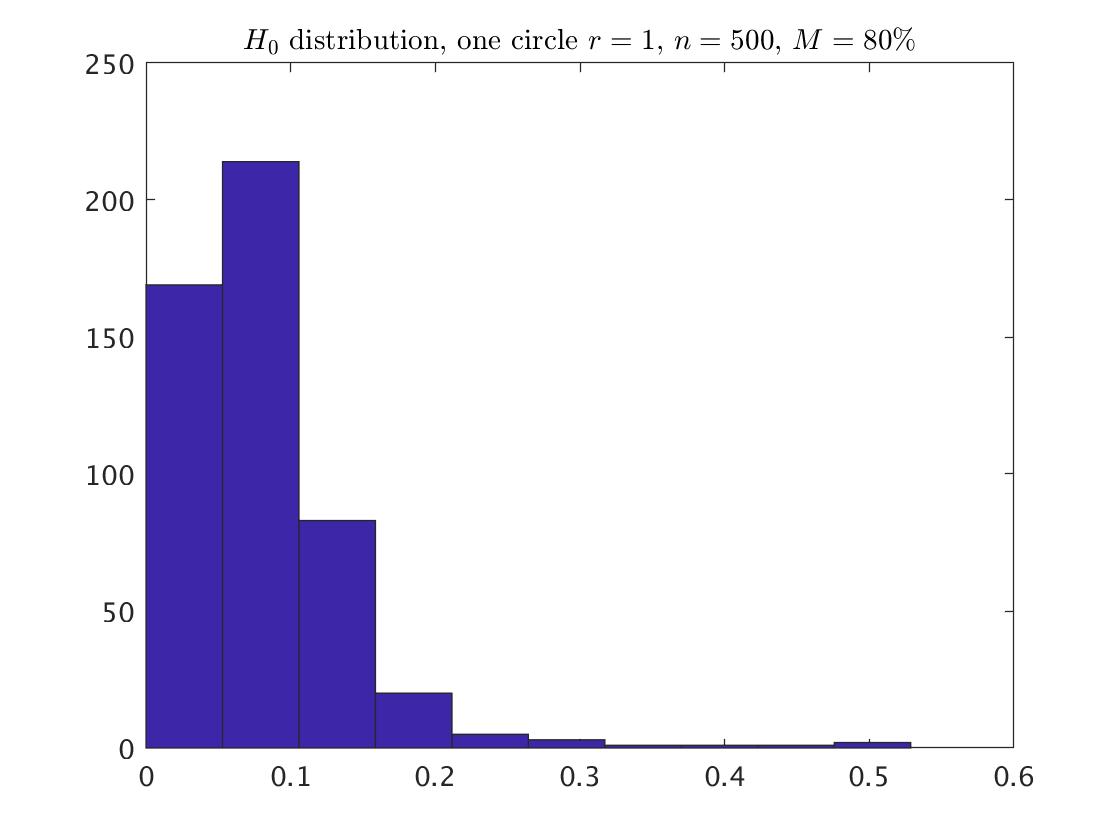} \hskip0.01truein
\includegraphics[width=1.45in, height=1.45in]{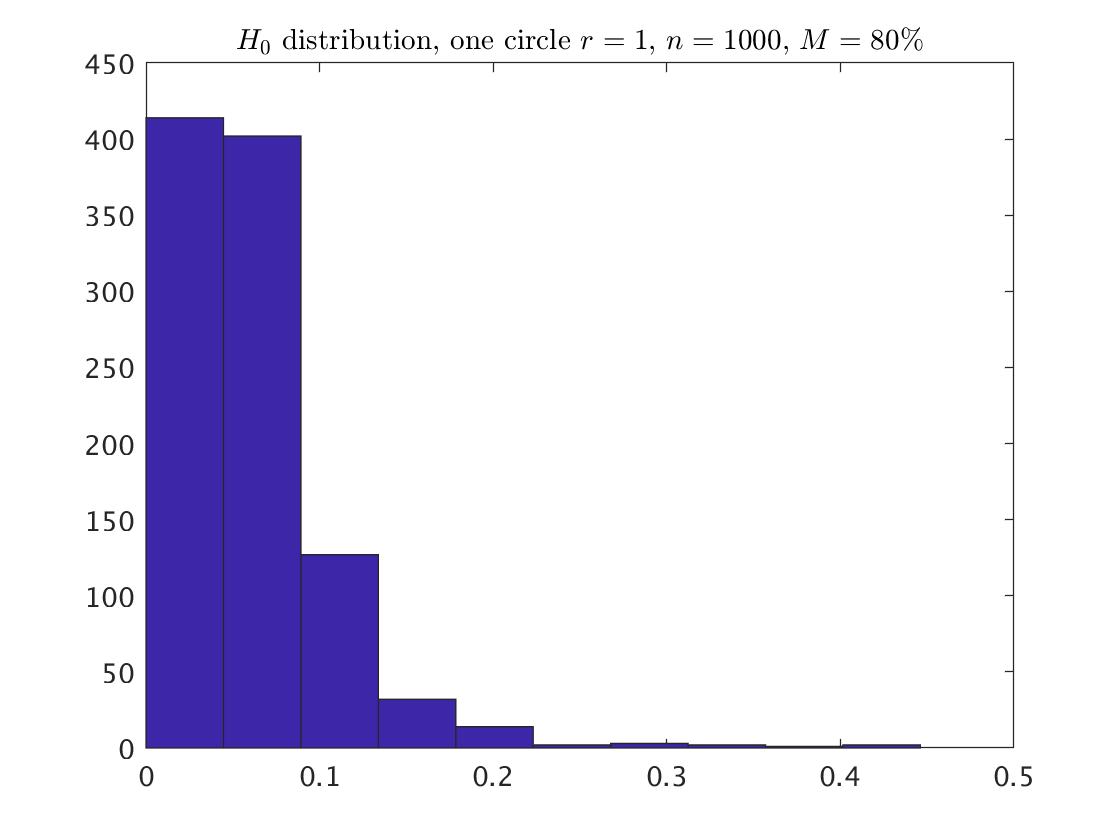} \hskip0.01truein
\includegraphics[width=1.45in, height=1.45in]{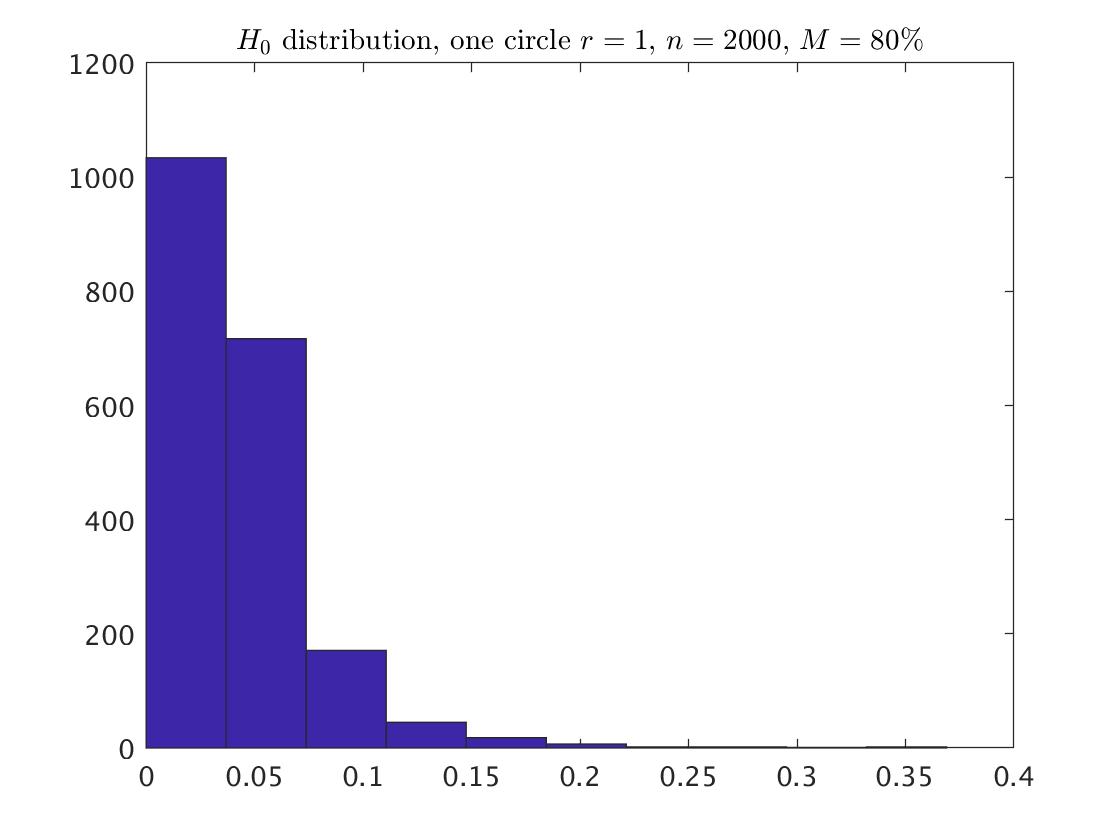} \hskip0.01truein

\includegraphics[width=1.45in, height=1.45in]{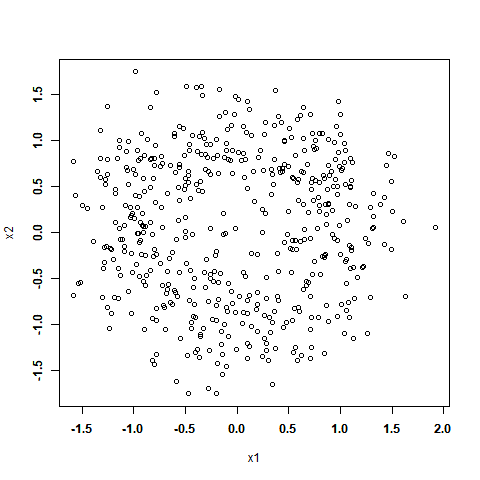} \hskip0.01truein
\includegraphics[width=1.45in, height=1.45in]{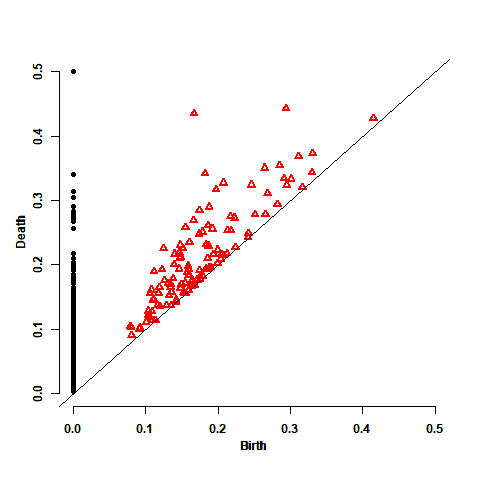} \hskip0.01truein
\includegraphics[width=1.45in, height=1.45in]{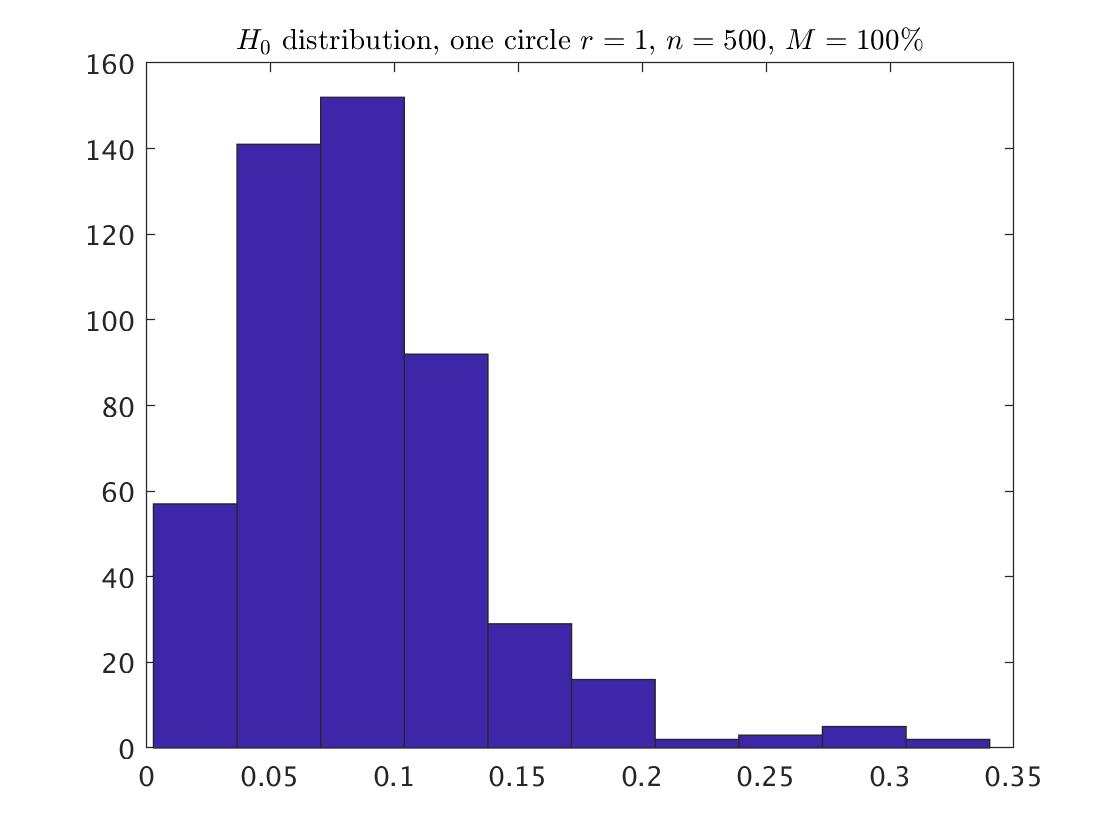} \hskip0.01truein
\includegraphics[width=1.45in, height=1.45in]{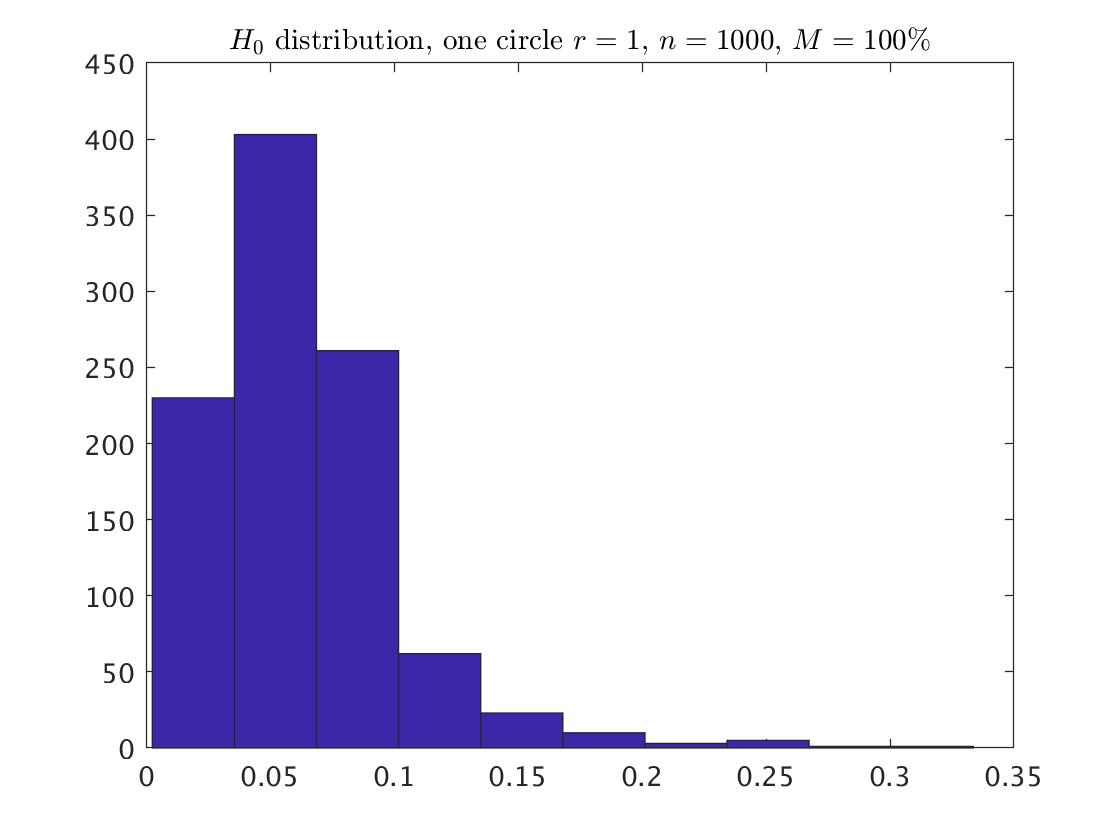} \hskip0.01truein
\includegraphics[width=1.45in, height=1.45in]{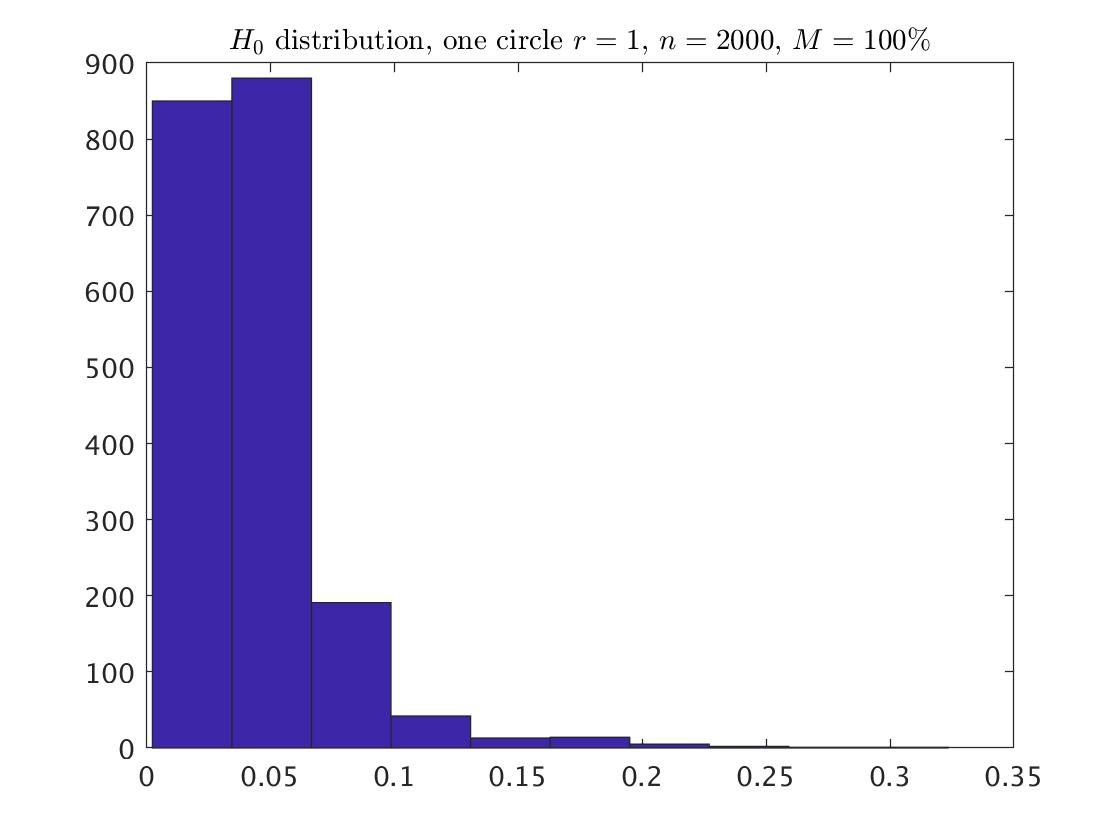} \hskip0.01truein

\ec
\caption{\footnotesize
 One circle with radius 1 and noise of $M\%$ of the sample size $n$. From top to bottom: $M = 30\%, 70\%, 80\%, 100\%$. Each row contains the circle with $n=500$ points, and to its right the histograms of the $H_0$ persistence diagrams for $n=500, 1,000, 2,000$.
}
\label{fig:onecircle_M_noise}
\end{figure}
\end{landscape}

\begin{landscape}
\begin{figure}[h!]
\bc
\includegraphics[width=1.45in, height=1.45in]{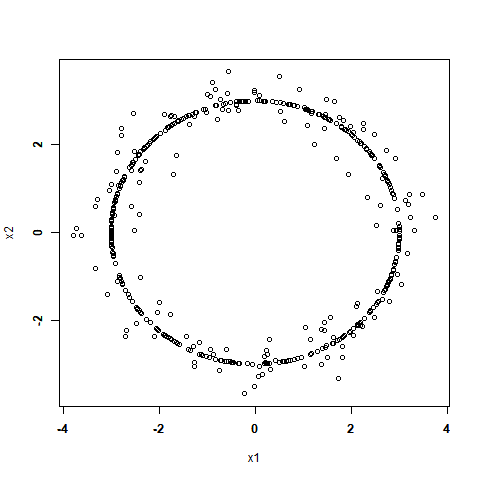} \hskip0.01truein
\includegraphics[width=1.45in, height=1.45in]{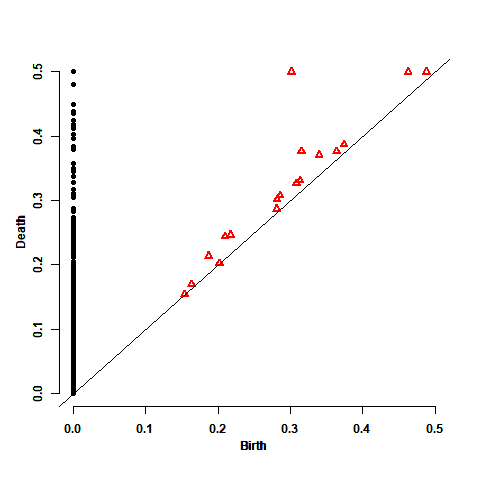} \hskip0.01truein
\includegraphics[width=1.45in, height=1.45in]{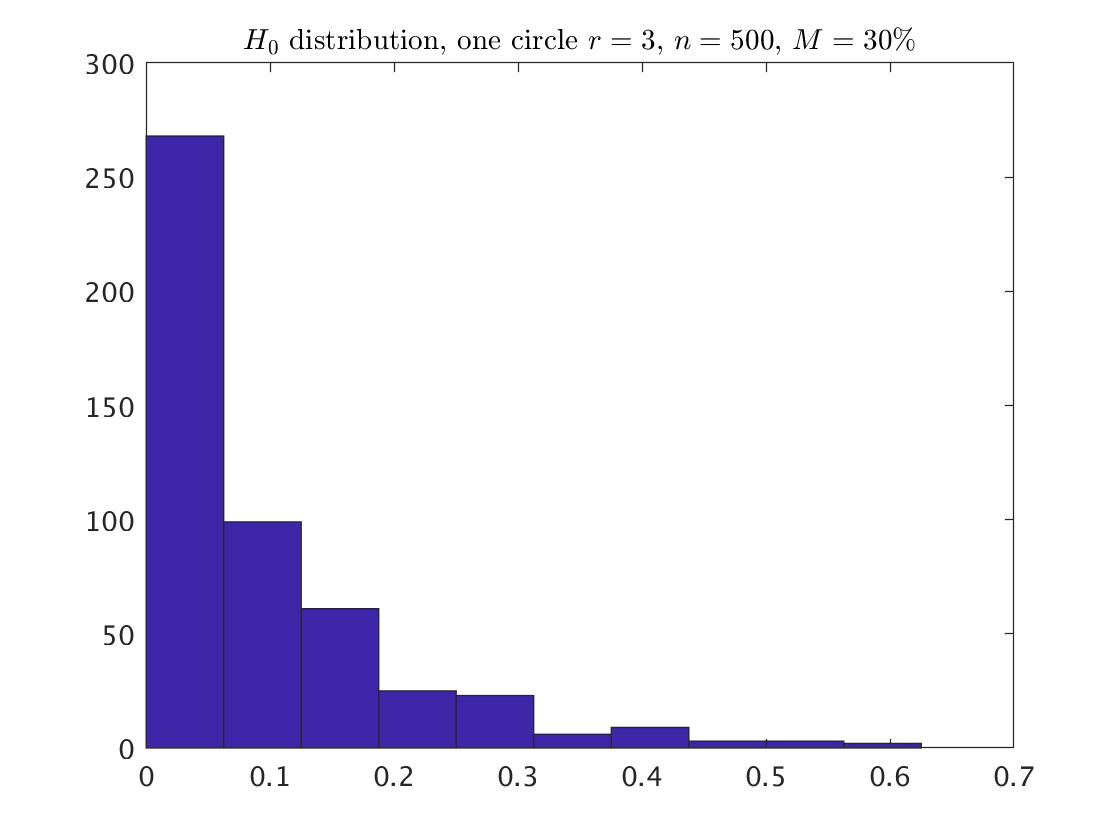} \hskip0.01truein
\includegraphics[width=1.45in, height=1.45in]{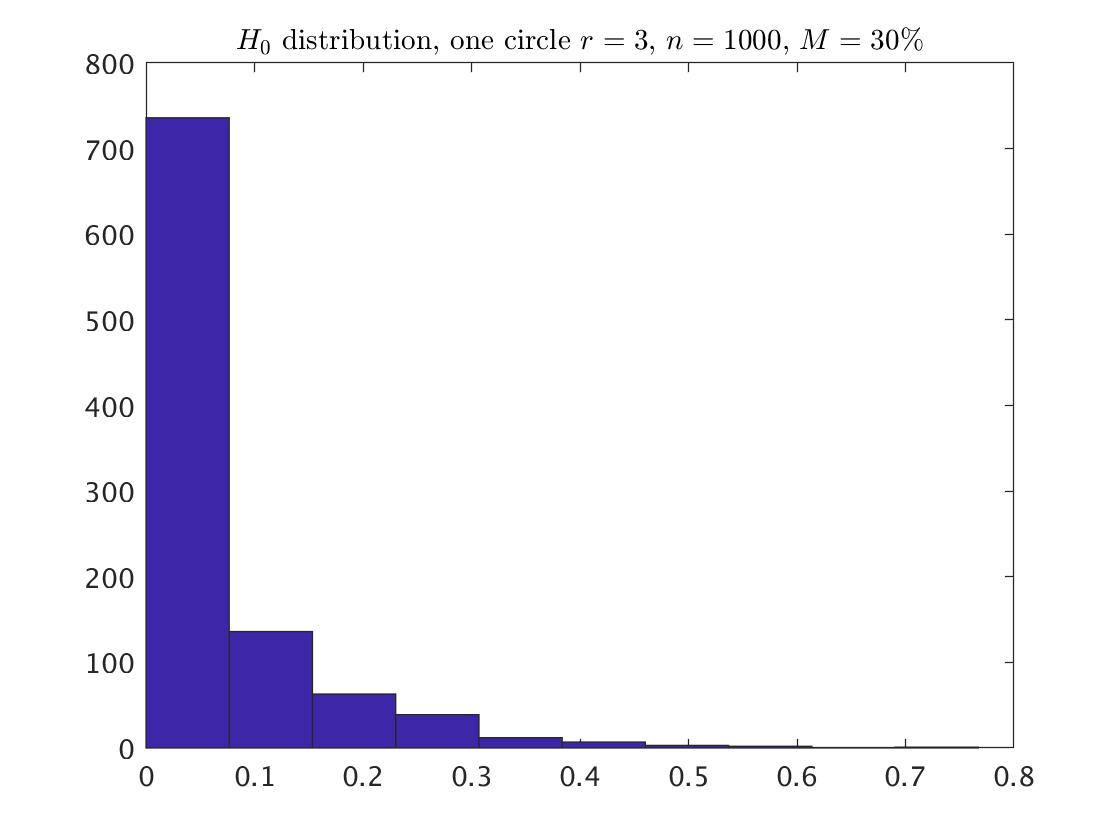} \hskip0.01truein
\includegraphics[width=1.45in, height=1.45in]{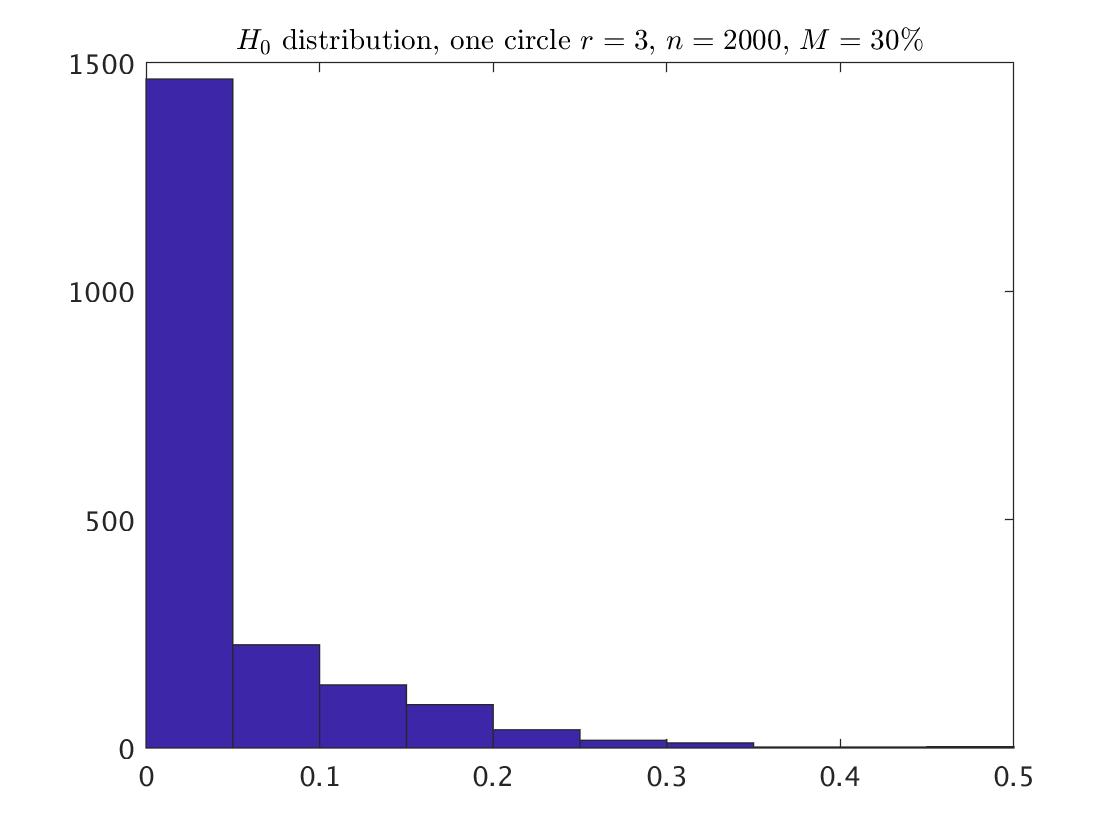} \hskip0.01truein

\includegraphics[width=1.45in, height=1.45in]{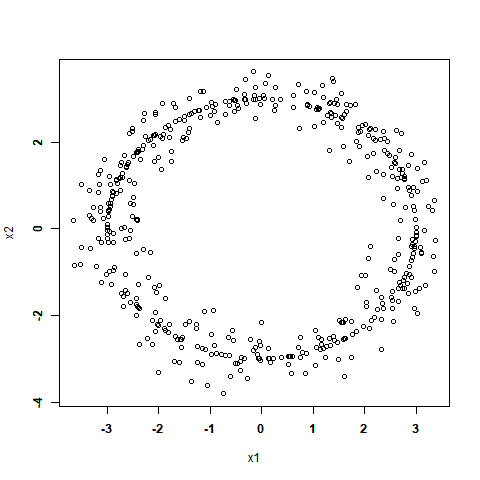} \hskip0.01truein
\includegraphics[width=1.45in, height=1.45in]{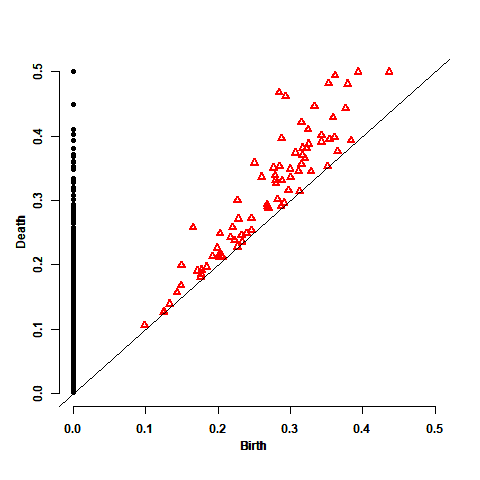} \hskip0.01truein
\includegraphics[width=1.45in, height=1.45in]{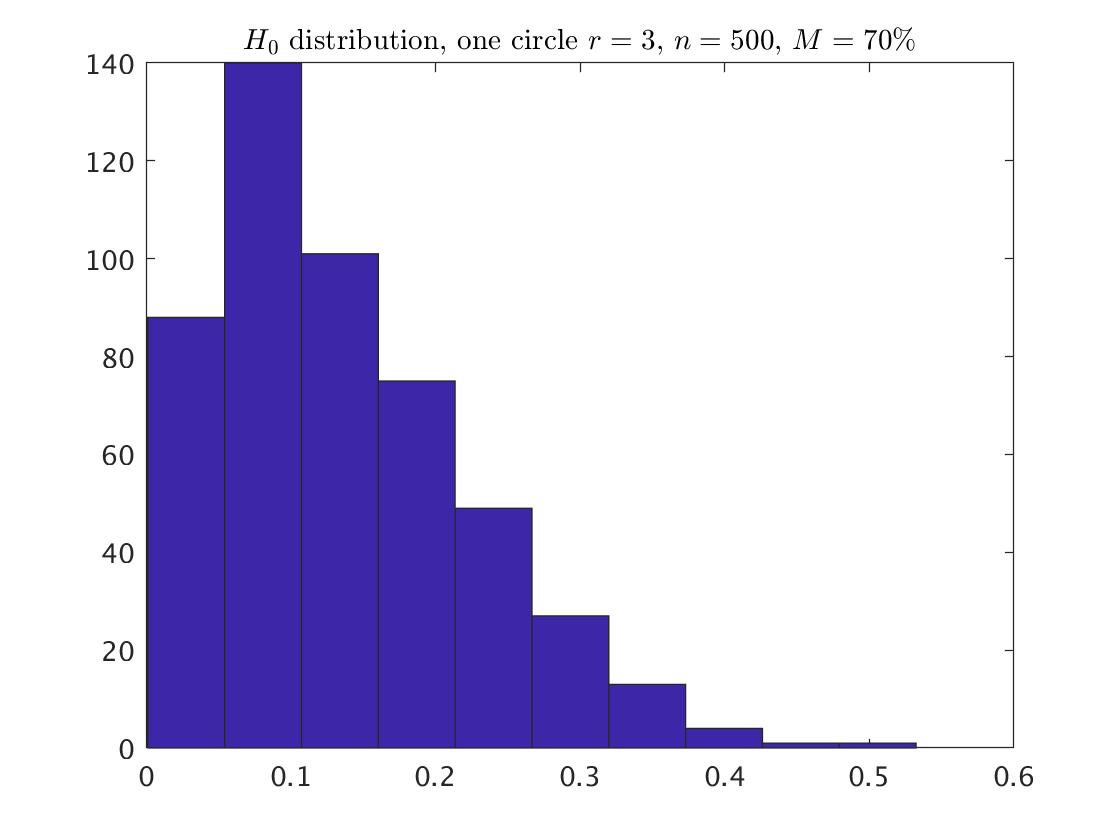} \hskip0.01truein
\includegraphics[width=1.45in, height=1.45in]{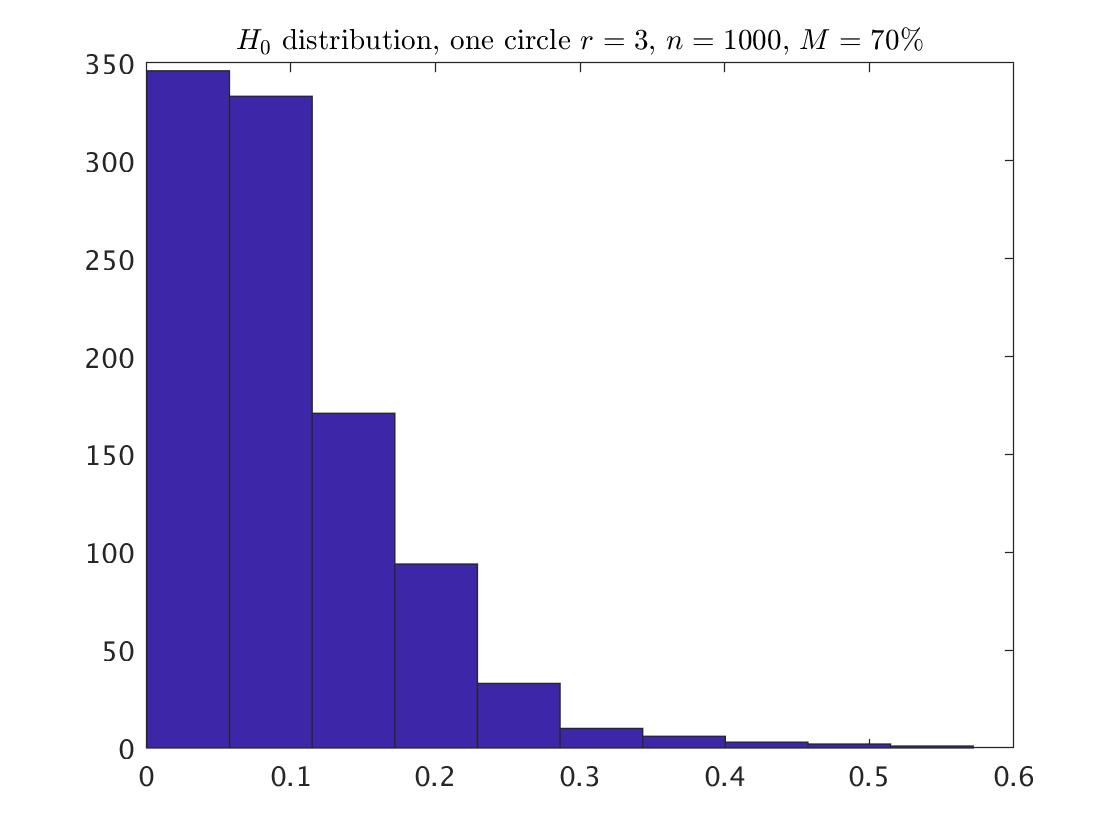} \hskip0.01truein
\includegraphics[width=1.45in, height=1.45in]{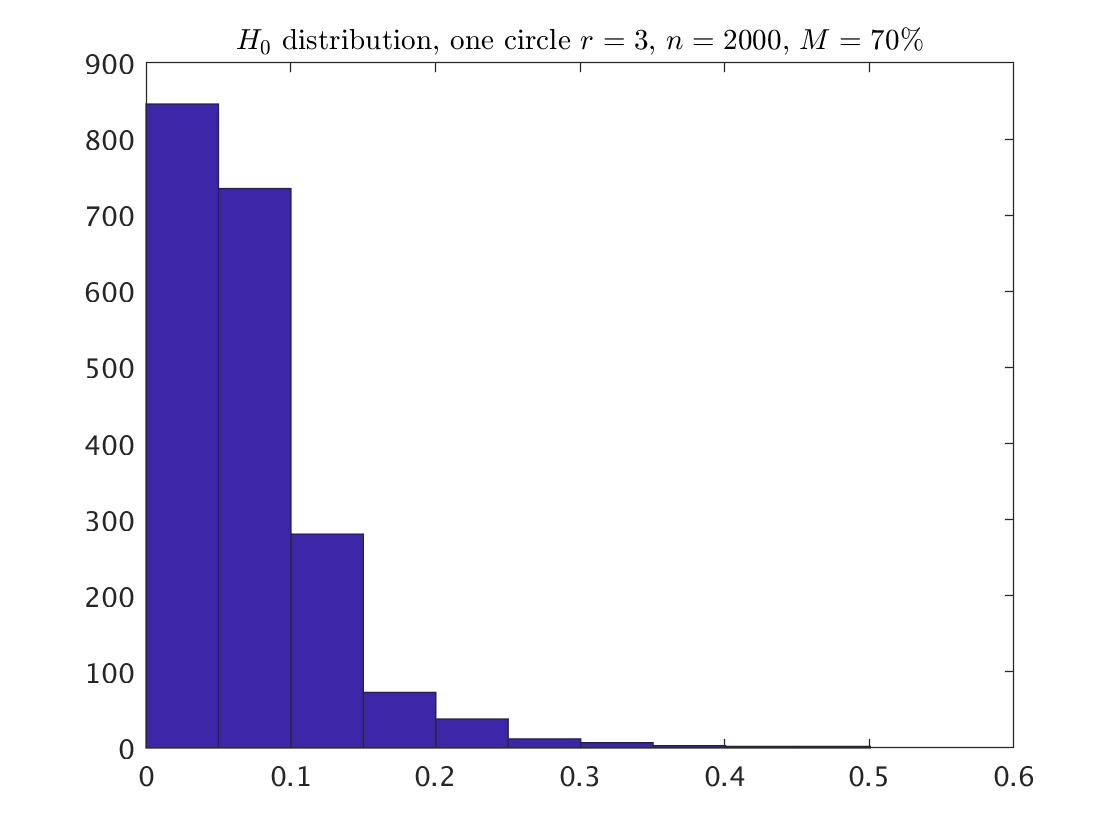} \hskip0.01truein

\includegraphics[width=1.45in, height=1.45in]{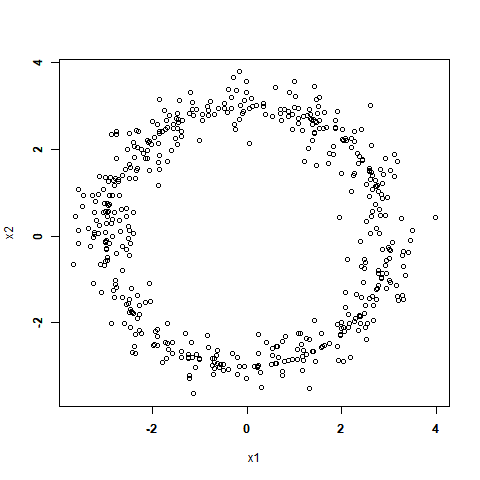} \hskip0.01truein
\includegraphics[width=1.45in, height=1.45in]{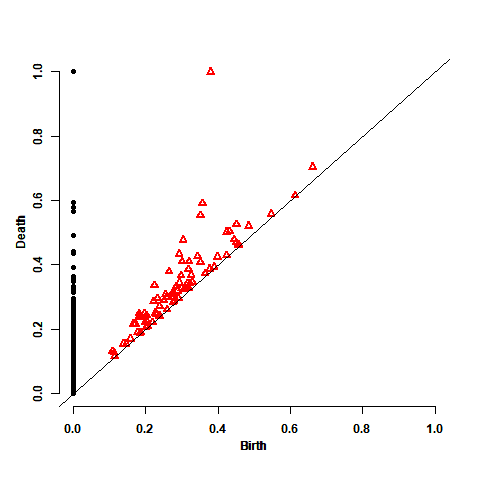} \hskip0.01truein
\includegraphics[width=1.45in, height=1.45in]{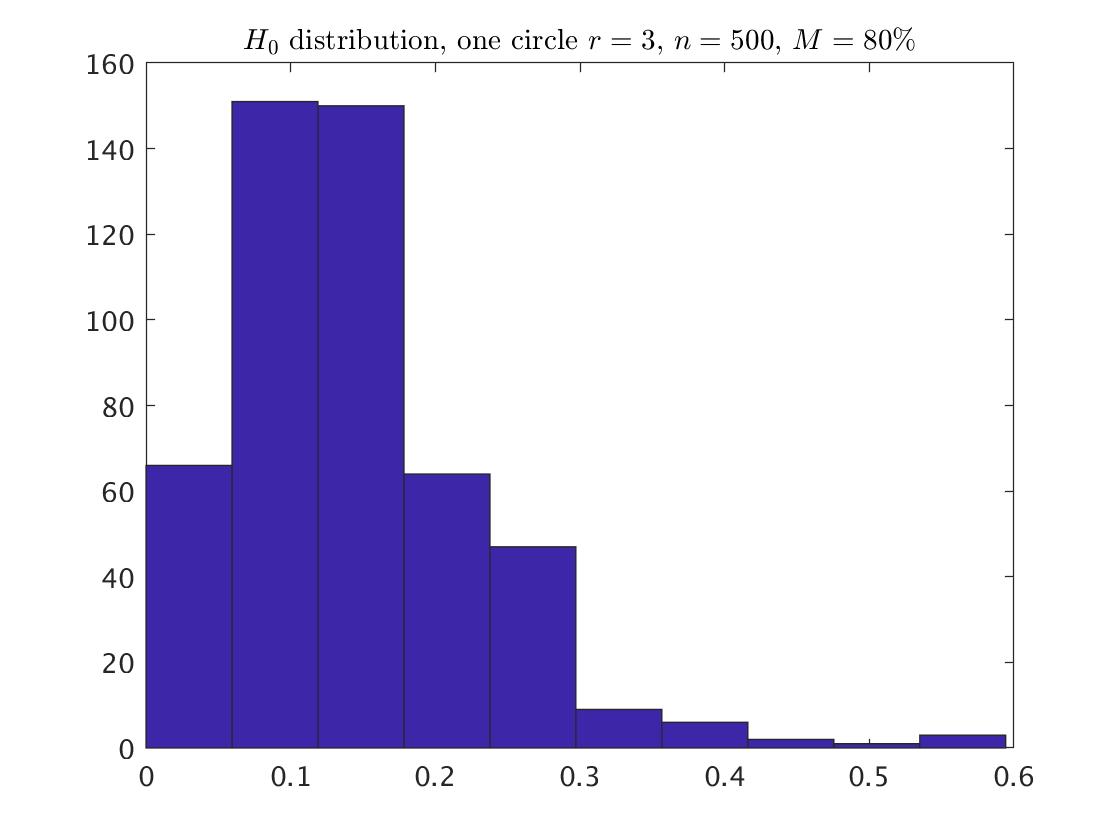} \hskip0.01truein
\includegraphics[width=1.45in, height=1.45in]{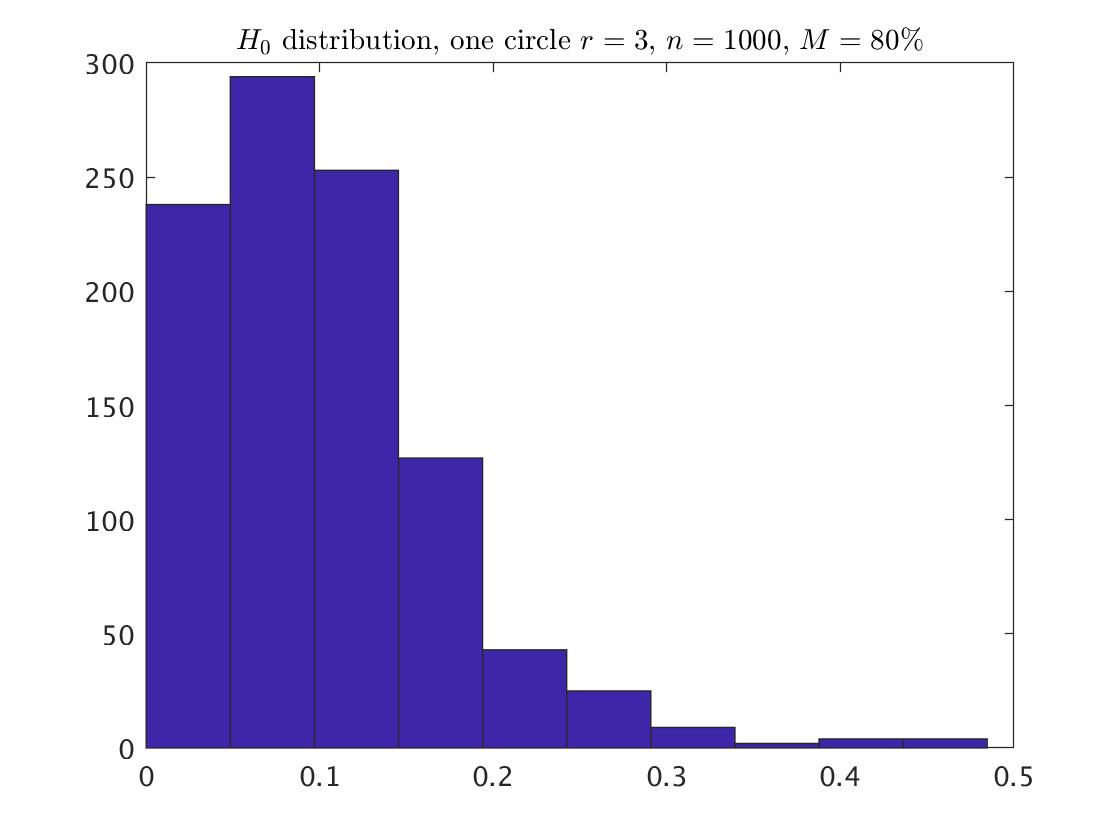} \hskip0.01truein
\includegraphics[width=1.45in, height=1.45in]{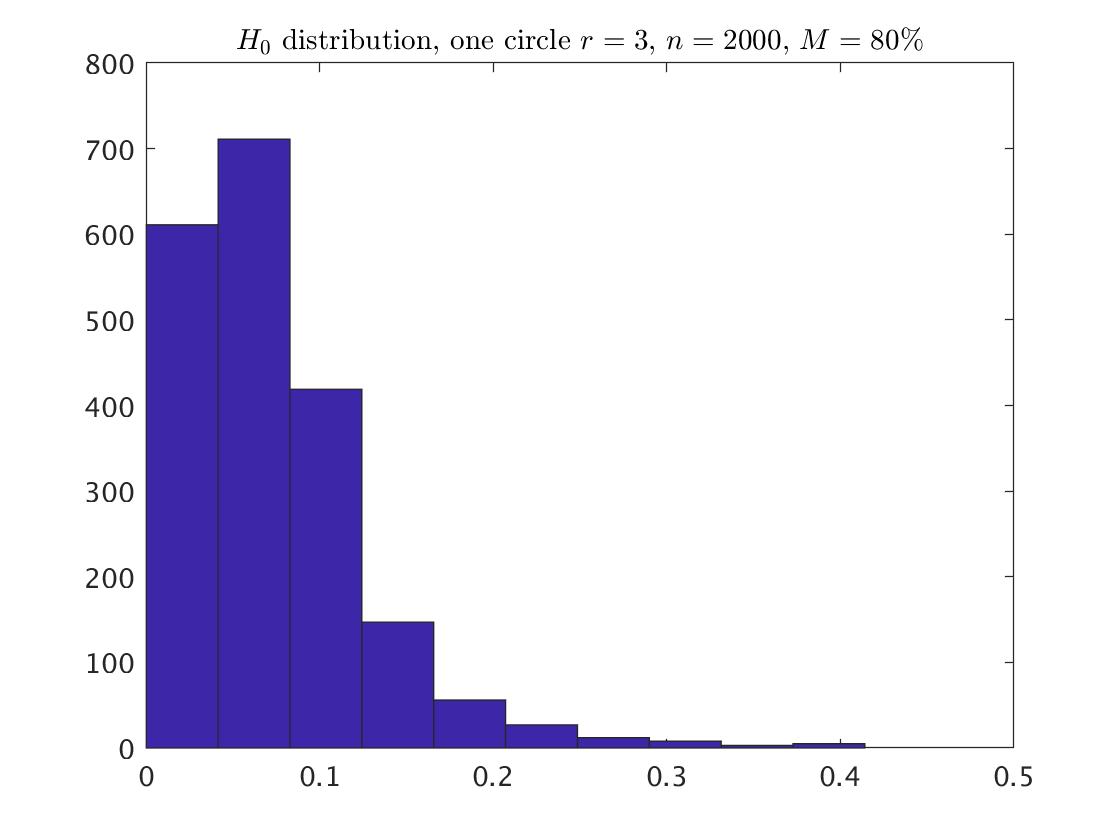} \hskip0.01truein

\includegraphics[width=1.45in, height=1.45in]{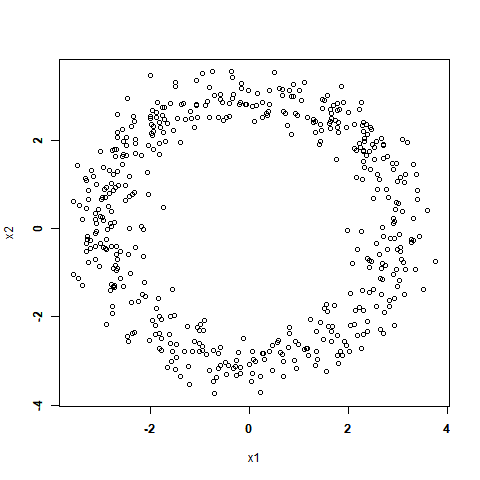} \hskip0.01truein
\includegraphics[width=1.45in, height=1.45in]{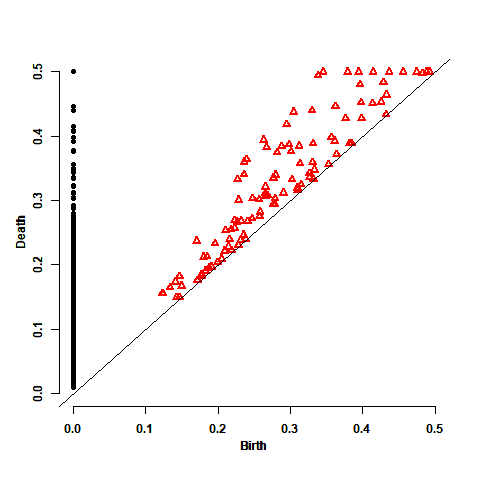} \hskip0.01truein
\includegraphics[width=1.45in, height=1.45in]{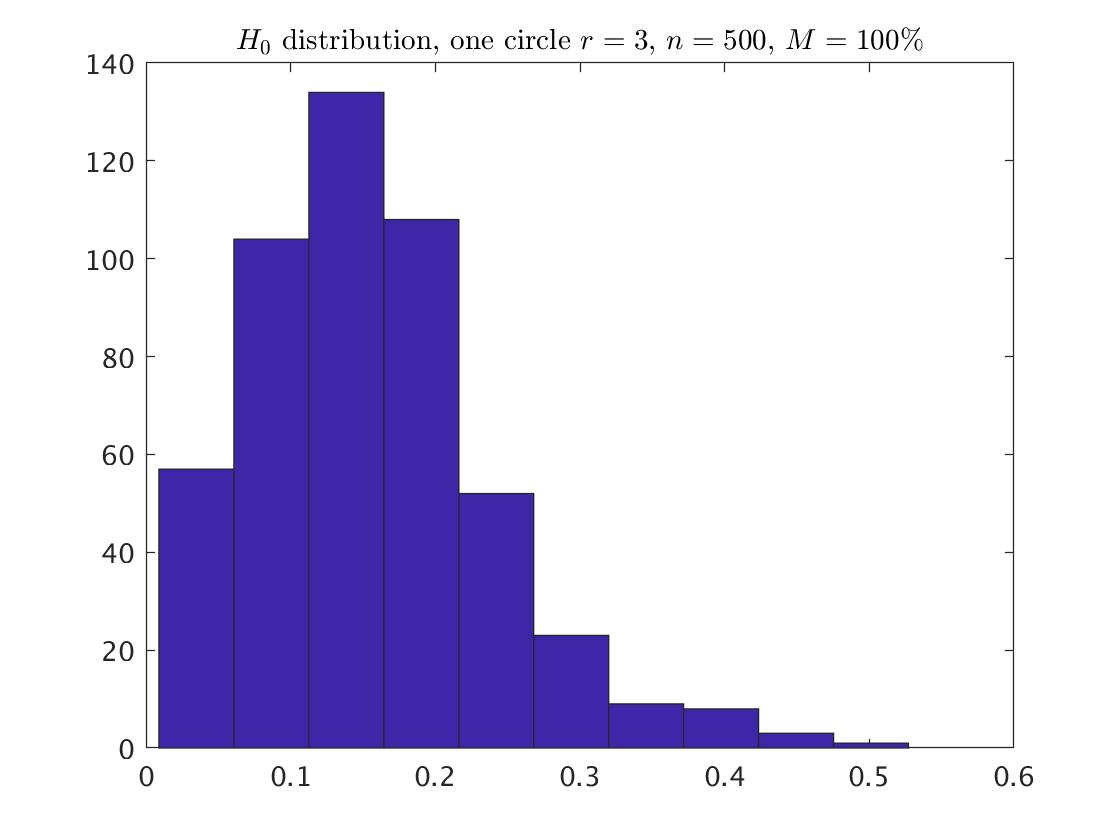} \hskip0.01truein
\includegraphics[width=1.45in, height=1.45in]{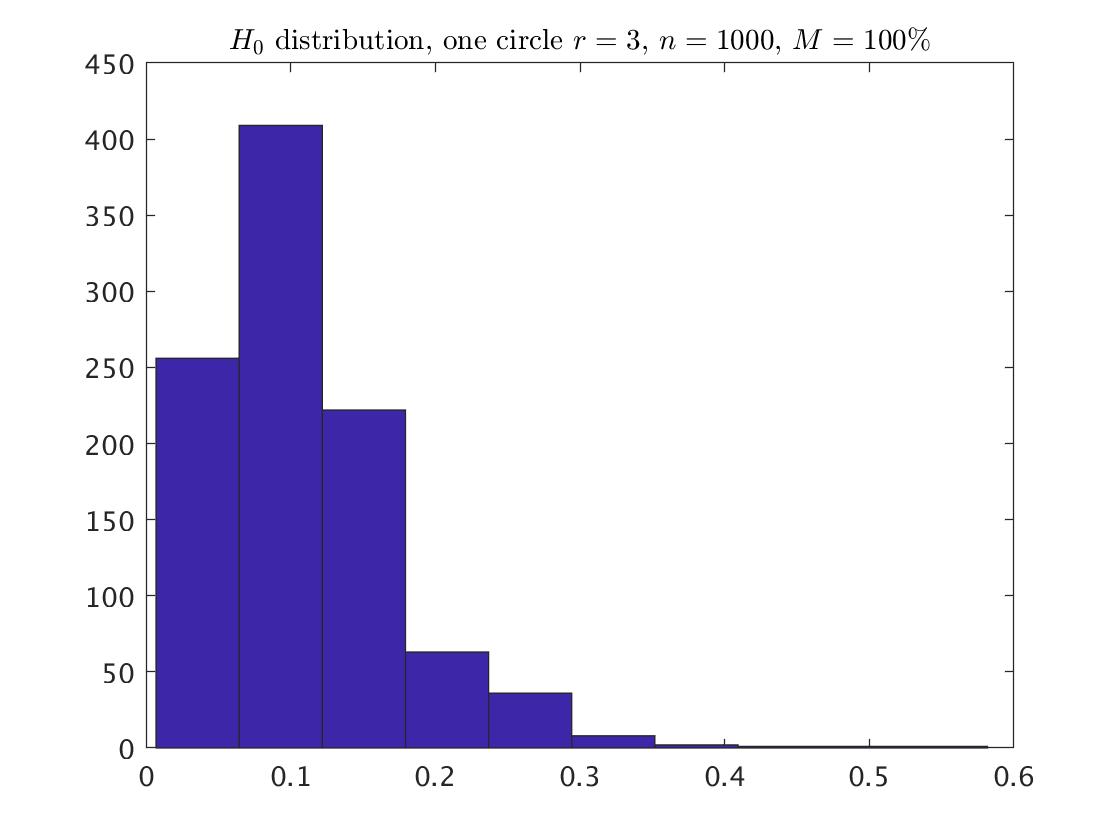} \hskip0.01truein
\includegraphics[width=1.45in, height=1.45in]{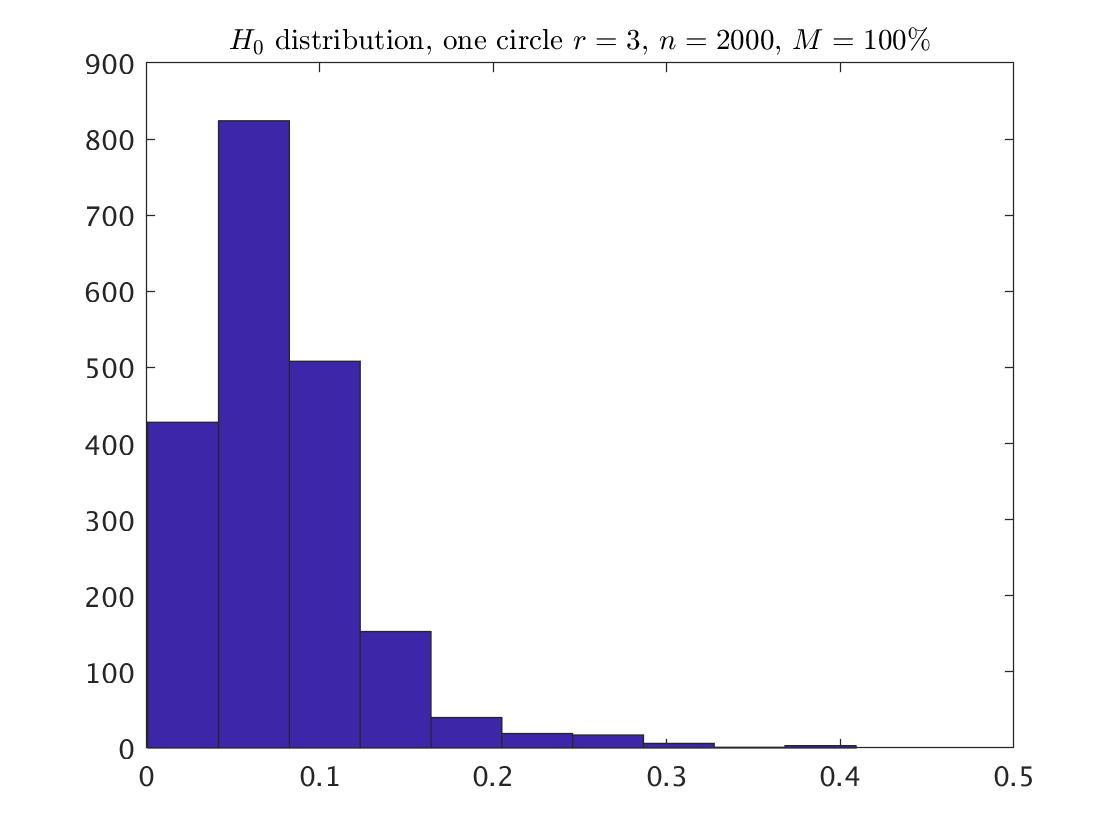} \hskip0.01truein

\ec
\caption{\footnotesize
 One circle with radius 3 and noise of $M\%$ of the sample size $n$. From top to bottom: $M = 30\%, 70\%, 80\%, 100\%$. Each row contains the circle with $n=500$ points, and to its right the histograms of the $H_0$ persistence diagrams for $n=500, 1,000, 2,000$.
}
\label{fig:onecircle3_M_noise}
\end{figure}
\end{landscape}

For noise of $M$=30\% and $M$=70\%, maxscale of 0.5 is large enough to capture the whole data, whereas for noise of $M$=80\%, needs to take maxscale of 1. These values, comparing to the maxscale of 0.3 that was enough for the clean one circle, are larger. The reason for these differences is the points that are located inside the circle. The small maxscale can capture the small connected components that the inside circle points obtain.
Fig. \ref{fig:onecircle_M_noise} and Fig. \ref{fig:onecircle3_M_noise} describe for each value of $M$ the circle with $r=1$ and $r=3$, respectively, and the distribution of the $H_0$ points (without the point at infinity) based on the above relevant value of the maxscale.

\normalsize

\normalsize

For examining the influence of noise on the lengths of the bars, we classified the bars into two groups of small and long bars, depends on the $c_{max}$ of the relevant $n$. The results, depending on noise and maxscale, are presented in Table \ref{table:OneCirc_class}. Clearly, the proportion of the long bars increases as the noise increases. The difference between the various values of maxscale for a given level of noise is in the lengths of the longest bars, whereas the proportion of the long bars is the same. Table \ref{table:OneCirc_percentile} summarizes for each maxscale and a given noise, the $95$-th, $99$-th, and the $100$-th percentiles of the $H_0$ lengths. For a given noise, and a specific maxscale, a percentile that is the same as for a lower maxscale is omitted. The results are: In the circle with $r=1$, the difference between the various values of maxscale is usually only at the $100$-th percentile of the $H_0$ lengths. In the circle with $r=3$, the difference is at the $99$-th percentile, but for $M=30\%$ is at the $95$-th percentile for $n=500$, at the $99$-th percentile for $n=1,000$, and at the $100$-th percentile for $n=2,000$. Checking in more details the lengths under maxscale=0.3 and $M=30\%$ yields that the proportion of lengths that equal to 0.3 is 0.056 for $n=500$, 0.026 for $n=1,000$, and 0.009 for $n=2,000$.

\begin{center}
\fontsize{8.5}{0.9}\selectfont
\captionof{table}{One circle - classification of $H_0$ points}
\begin{tabular}{l|lc|cc|cc}
\\
\\
$r=1$&$n=500$$^a$ &&$n=1,000$$^b$&&$n=2,000$$^c$\\\hline
\\
Noise&short& long&short& long&short& long\\\hline
\\
\\
\\
\\
30\%& 0.747&0.253&0.744& 0.256&0.719& 0.281\\
\\
\\
\\
\\
70\%& 0.445&0.555 &0.439&0.561 &0.368&0.632\\
\\
\\
\\
\\
80\%& 0.399&0.601& 0.363&0.637 &0.267&0.733\\
\\
\\
\\
\\
100\%& 0.283&	0.717& 0.270&0.730 &0.170&0.830\\
\\
\\
\\
\\
\\
\\
\\
\\
\\
\\
\\
\\
& &&&&\\\hline
\\
\\
\\
$r=3$&$n=500$$^d$ &&$n=1,000$$^e$&&$n=2,000$$^f$\\\hline
\\
\\
\\
\\
30\%& 0.840&0.160&0.817	&0.183&0.768&0.232\\
\\
\\
\\
\\
70\%& 0.723&0.277 &0.680&0.320 &0.515&0.485\\
\\
\\
\\
\\
80\%& 0.729&0.271& 0.623&0.377 &0.480&0.520\\
\\
\\
\\
\\
100\%& 0.639&0.361&0.611&0.389 &0.413&0.587\\
\label{table:OneCirc_class}
\end{tabular}
\end{center}
\footnotesize{Classification of the $H_0$ lengths into short and long bars for one circle with radius $r$. 'Short' and 'long' are relative to the relevant $c_{max}$. The classification is based on one persistence diagram, with maxscale$\ge$ 0.3. $^a$ $c_{max}=0.059$, $^b$ $c_{max}=0.038$, $^c$ $c_{max}=0.020$, $^d$ $c_{max}=0.176$, $^e$ $c_{max}=0.115$, $^f$ $c_{max}=0.061$.}\\

\newpage
\normalsize
\begin{center}
\fontsize{8.5}{0.9}\selectfont
\captionof{table}{One circle - percentiles of $H_0$ points}
\begin{tabular}{l|l|c|c|c|c}
\\
\\
$r=1$& &&&\\\hline
\\
Noise&maxscale&percentile $\%$&$n=500$$^a$ &$n=1,000$$^b$&$n=2,000$$^c$\\\hline
\\
\\
\\
\\
30\%&0.3& 95&0.212&	0.156&	0.110\\
\\
\\
\\
\\
&& 99&0.212&	0.156&	0.110\\
\\
\\
\\
\\
&& 100&0.300&	0.300&	0.300\\
\\
\\
\\
\\
& 0.5& 100&0.482&0.446&0.500\\
\\
\\
\\
\\
& 1& 100&0.482&0.446&	0.519\\
& &&&\\\hline
\\
\\
\\
70\%&0.3& 95&0.176&	0.136&	0.099\\
\\
\\
\\
\\
&& 99&0.288&	0.211&	0.182\\
\\
\\
\\
\\
&& 100&0.300&	0.300&	0.300\\
\\
\\
\\
\\
& 0.5& 99&0.294&	0.211&0.182\\
\\
\\
\\
\\
&& 100&0.426&	0.499&	0.500\\
\\
\\
\\
\\
& 1& 100&0.426&	0.499&	0.523\\
& &&&\\\hline
\\
\\
\\
80\%& 0.3& 95&0.174&	0.138&	0.102\\
\\
\\
\\
\\
&& 99&0.294&	0.224&	0.163\\
\\
\\
\\
\\
&& 100&0.300&	0.300&	0.300\\
\\
\\
\\
\\
& 0.5& 99&0.310&	0.224&	0.163\\
\\
\\
\\
\\

&& 100&0.500&	0.446&	0.369\\
\\
\\
\\
\\
&1& 100&0.529&	0.446&0.369\\
& &&&\\\hline
\\
\\
\\
100\%& 0.3& 95&0.177&	0.125&	0.089\\
\\
\\
\\
\\
&& 99&0.281&	0.204&	0.170\\
\\
\\
\\
\\
&& 100&0.300&	0.300&	0.300\\
\\
\\
\\
\\
& 0.5& 100&0.340&	0.334&	0.323\\

& &&&\\\hline
& &&&\\\hline
\\
\\
\\
$r=3$& &&&\\\hline
\\
\\
\\
\\
\\
30\%&0.3& 95&0.300&	0.251&	0.186\\
\\
\\
\\
\\
&& 99&0.300&	0.300&	0.295\\
\\
\\
\\
\\
&& 100&0.300&	0.300&	0.300\\
\\
\\
\\
\\
& 0.5& 95&0.310&0.251&0.186\\
\\
\\
\\
\\
& & 99&0.491&	0.432&	0.295\\
\\
\\
\\
\\
& & 100&0.500&	0.500&	0.500\\
\\
\\
\\
\\
&1 & 99&0.500&	0.432&	0.295\\
\\
\\
\\
\\
& 1& 100&0.625&	0.767&	0.565\\
& &&&\\\hline
\\
\\
\\
70\%&0.3& 95&0.292&	0.239&	0.169\\
\\
\\
\\
\\
&& 99&0.300&	0.300&	0.265\\
\\
\\
\\
\\
&& 100&0.300&	0.300&	0.300\\
\\
\\
\\
\\
& 0.5& 99&0.387&	0.380&	0.265\\
\\
\\
\\
\\
&& 100&0.500&0.500&	0.500\\
\\
\\
\\
\\
& 1& 100&0.532&	0.572&	0.501\\
& &&&\\\hline
\\
\\
\\
80\%& 0.3& 95&0.285&	0.235&	0.171\\
\\
\\
\\
\\
&& 99&0.300&	0.300&	0.267\\
\\
\\
\\
\\
&& 100&0.300&	0.300&	0.300\\
\\
\\
\\
\\
& 0.5& 99&0.440&	0.338&	0.267\\
\\
\\
\\
\\

&& 100&0.500&	0.485&	0.414\\
\\
\\
\\
\\
&1& 100&0.594&	0.485&	0.414\\
& &&&\\\hline
\\
\\
\\
100\%& 0.3& 95&0.300&	0.234&	0.160\\
\\
\\
\\
\\
&& 99&0.300&	0.300&	0.256\\
\\
\\
\\
\\
&& 100&0.300&	0.300&	0.300\\
\\
\\
\\
\\
& 0.5& 95&0.304&	0.234&	0.160\\
\\
\\
\\
\\
&& 99&0.412&	0.303&	0.256\\
\\
\\
\\
\\
&&100&0.500&	0.500&	0.409\\
\\
\\
\\
\\
& 1& 100&0.527&	0.582&	0.409\\
\\
\\
\\
\\
\label{table:OneCirc_percentile}
\end{tabular}
\end{center}
\footnotesize{Percentiles of the $H_0$ lengths. The $H_0$ lengths are based on one persistence diagram, with maxscale$\ge$ 0.3. See text for more details.}\\
\normalsize

\begin{figure}[h!]
\bc
\includegraphics[width=1.45in, height=1.45in]{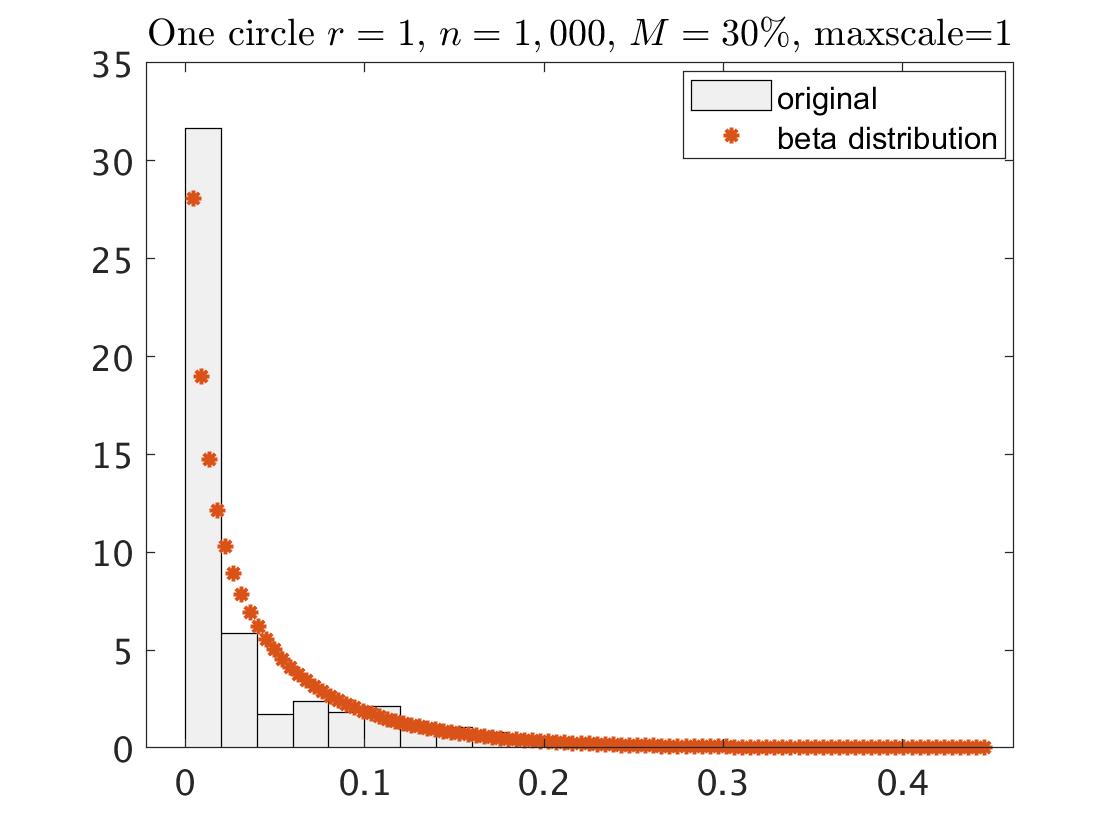} \hskip0.01truein
\includegraphics[width=1.45in, height=1.45in]{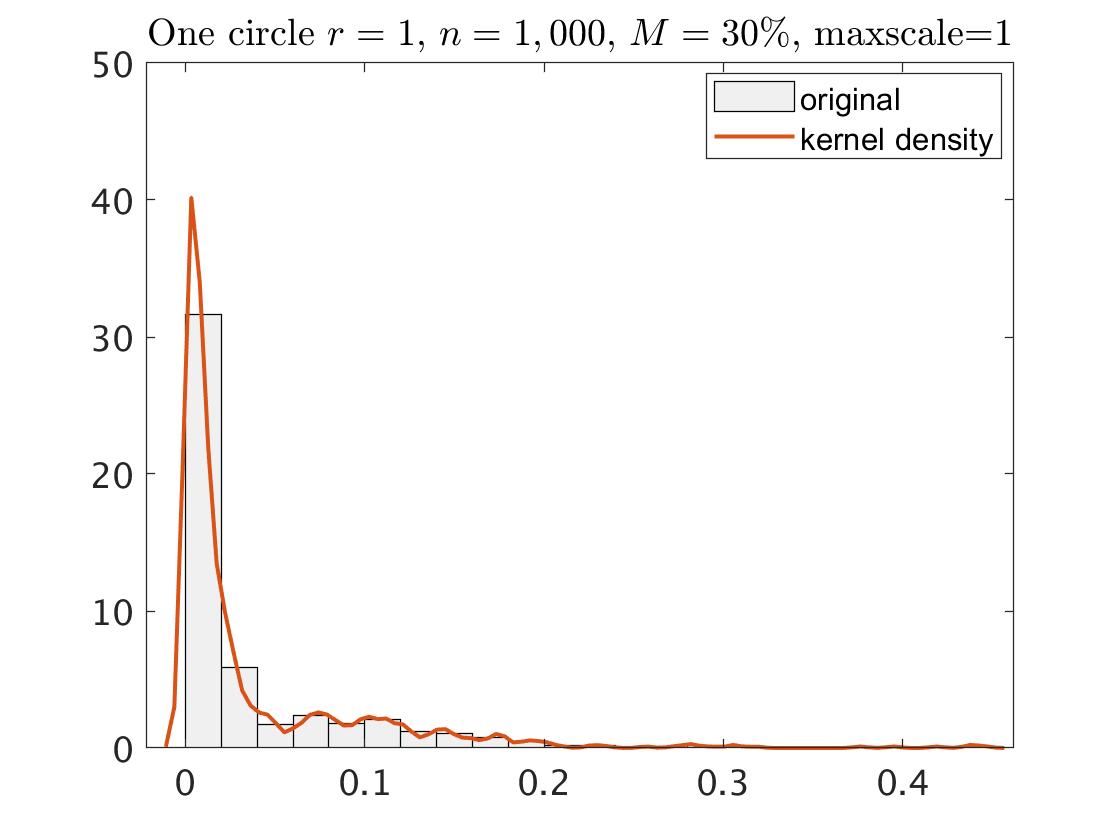} \hskip0.01truein
\includegraphics[width=1.45in, height=1.45in]{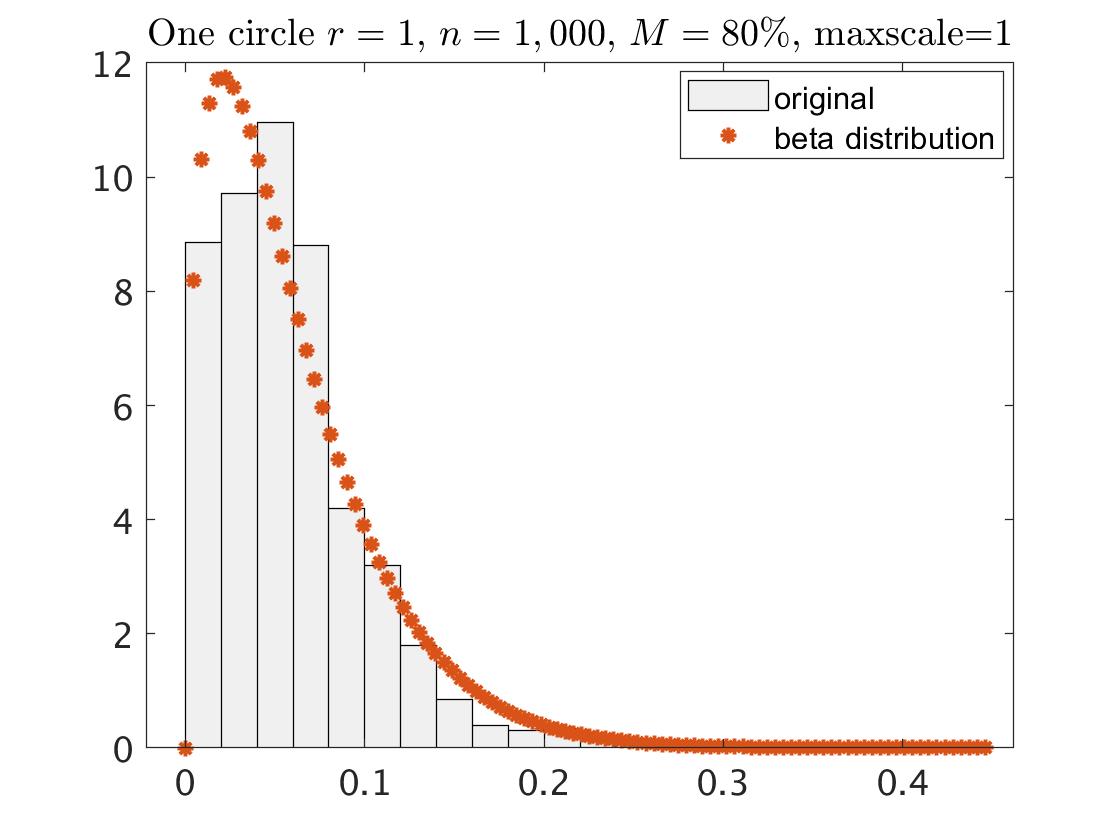} \hskip0.01truein
\includegraphics[width=1.45in, height=1.45in]{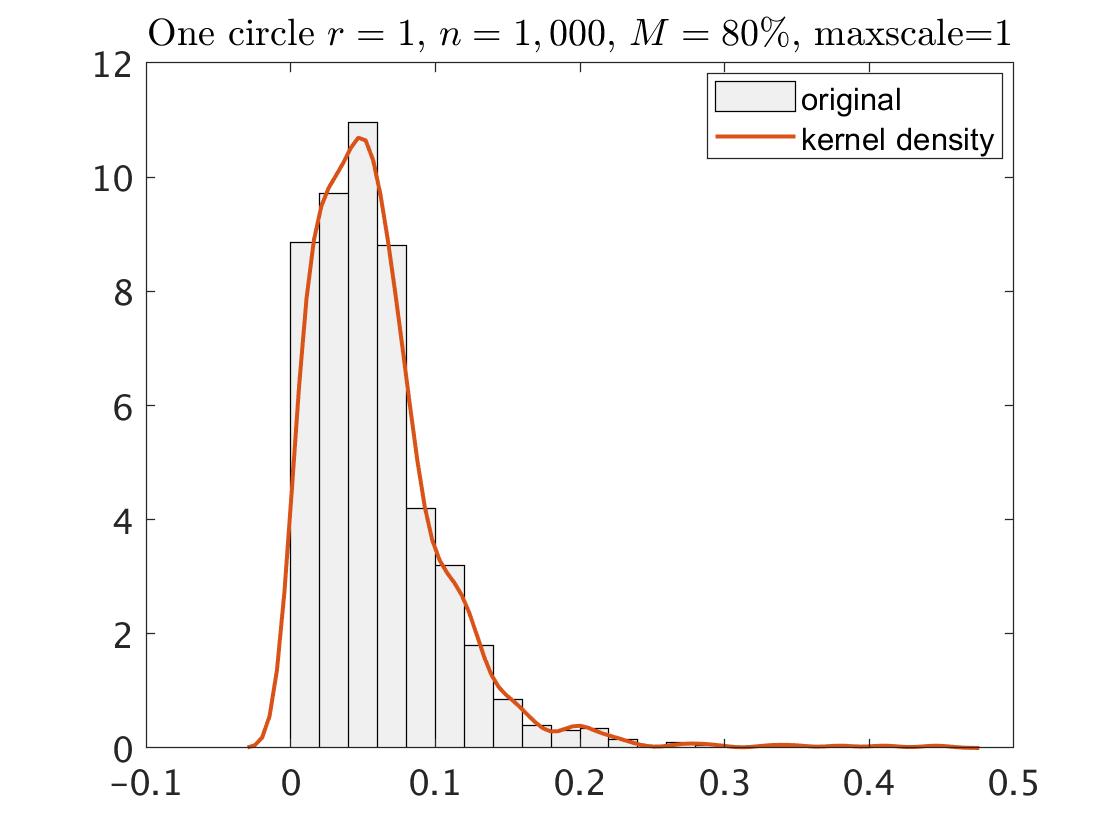} \hskip0.01truein
\ec
\caption{\footnotesize
 Parametric and non-parametric distributions of the $H_0$ points of one circle with $r=1$ and $n=1,000$, based on maxscale=1. The first two plots are based on data with $M = 30\%$ noise, whereas the last two plots are based on data with  $M = 80\%$ noise. }
\label{fig:onecircle_kernel}
\end{figure}

\begin{figure}[h!]
\bc
\includegraphics[width=1.55in, height=1.55in]{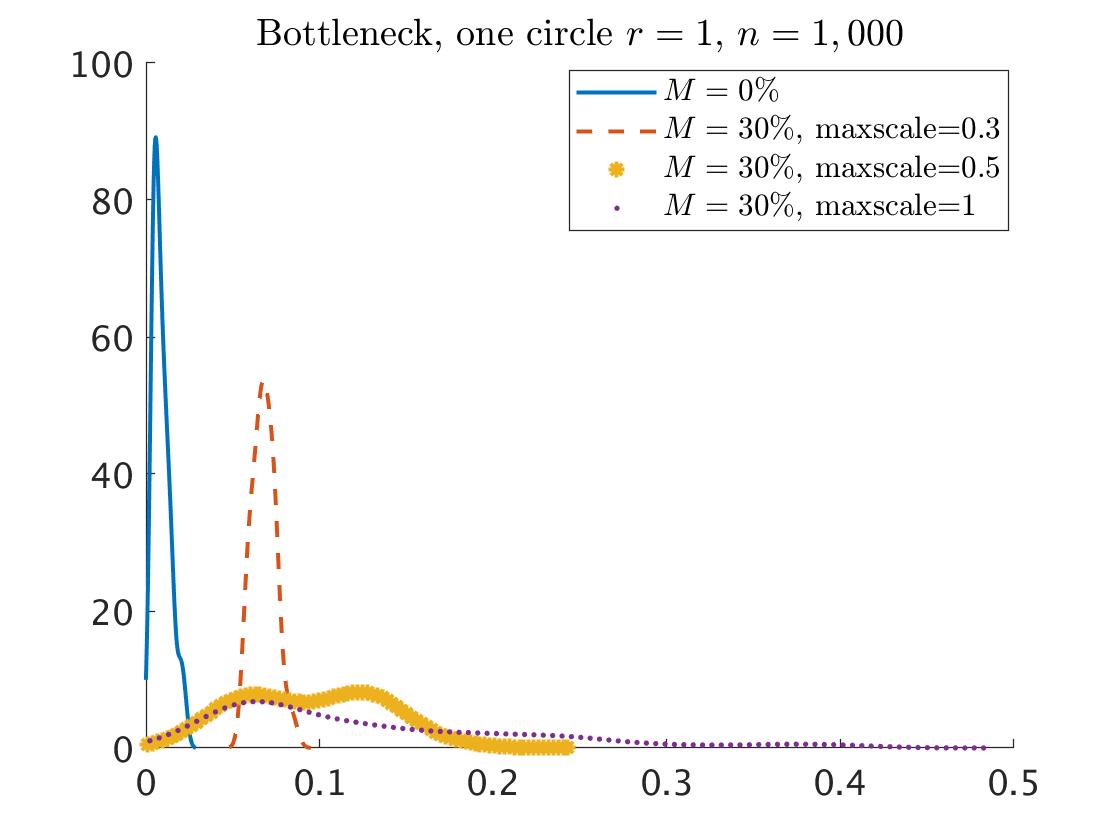} \hskip0.01truein
\ec
\caption{\footnotesize
 Smoothed empirical densities of the bottleneck distance between 100 real persistence diagrams and their corresponded simulated persistence diagrams, for one circle with $r=1$ and $n=1,000$.}
\label{fig:onecircle_bottle}
\end{figure}

Fitting a parametric distribution for the $H_0$ points in the same way as we did in the case of zero noise, yields that the best fitting is the beta distribution. The different parameters $a$ and $b$ of the beta distribution for each case are summarized in Appendix A.2. But, the goodness of the fit is not good, that is, the beta distribution  is no longer appropriate for describing the behaviour of the $H_0$ points in the setting of noisy data. For example, for the circle with $r=1$, $n=1,000$, and $M=30\%$, the distribution of the $H_0$ points together with the fitted beta distribution are described in the left plot of Fig. \ref{fig:onecircle_kernel}. The third plot of Fig. \ref{fig:onecircle_kernel} describes the same thing but for $M=80\%$. We can see that in both cases, the beta distribution does not fit well the distribution of $H_0$ points. In addition, we generated 100 collections of samples from one circle with $r=1$. Each sample contained $n=1,000$ points with fraction of $M=0\%$ noisy points. This is according to the same procedure that generated the original data with $r=1$, $n=1,000$, and $M=0\%$. For each sample we calculated the $H_0$ persistence diagram based on maxscale=1, and fitted the beta distribution to that points. In the next step we calculated for each $H_0$ persistence diagram its simulated diagram (using 999 random numbers from the fitted distribution). Finally, we computed the bottleneck distance (see, for example, \cite{Edelsbrunner2010}) between each pair of the real $H_0$ persistence diagram and its corresponded simulated $H_0$ persistence diagram. The distribution of the bottleneck distances over the 100 collections is described in the blue curve in Fig. \ref{fig:onecircle_bottle}. That is, for clean data, the distance between the real and the simulated $H_0$ persistence diagrams is close to zero, which indicate on a good matching between the real and the simulated persistence diagrams. We repeated this procedure for 100 samples that include $M=30\%$ noise, when we based each of the persistence diagrams on maxscale=0.3, 0.5, 1. The distributions of the bottleneck distances for each maxscale are described in the red, orange and purple curves in Fig. \ref{fig:onecircle_bottle}. We can see that the bottleneck distance is larger for $M=30\%$ relative to that of $M=0\%$, and it is getting larger as the maxscale increases. That is, the simulated persistence diagram becomes more different from the real persistence diagram as the maxscale increases.
Therefore needs to do something else, for example to examine non-parametric distribution such as the kernel density estimator, as describes in the second and the fourth plots of Fig. \ref{fig:onecircle_kernel}. We will
investigate it in further research.
\newpage

But still, using the parametric beta distribution fitting for the $H_0$ points, and following the statistical procedure for topological signals identifications that we used in the clean data, obtains that in most of the considered scenarios of the one circle we could identify correctly the one connected component. Table \ref{table:OneCirc_noisyStat} summarizes the results.
\newpage
\begin{center}
\fontsize{8.5}{0.9}\selectfont
\captionof{table}{One circle - confidence interval and p-value}
\begin{tabular}{l|lcc|ccc|cccccc}
\\
\textbf{$r=1$}&& $n=500$ &&& $n=1,000$&&& $n=2,000$\\\hline
\\
\\
\\
\\
Noise&$T_1$ real PD & CI & $p$-value& $T_1$ real PD & CI & $p$-value&  $T_1$ real PD & CI & $p$-value\\\hline
\\
\\
\\
\\
30\%, maxscale=0.3& 0.300 &[0,0.507]&0.976&0.300&[0,0.456]&0.870&0.300&[0,0.409]&0.550\\
\\
30\%, maxscale=0.5& 0.482 &[0,0.523]&0.142 &0.446&	[0,0.465]&	0.090&0.500&	[0,0.399]&	0.001\\
\\
30\%, maxscale=1& & &&	&	&&0.519&	[0,0.395]&	0.001\\
\\
\\
\\
\\
\\
\\
70\%, maxscale=0.3& 0.300 &[0,0.436]&0.917  &0.300&	[0,0.394]&	0.644&0.300&	[0,0.345]&	0.222 \\
\\
70\%, maxscale=0.5& 0.426 &[0,0.447]&0.104  &0.499&	[0,0.406]&	0.001&0.500&[0,0.341]&0 \\
\\
70\%, maxscale=1&  &&  &&	&	&0.523&[0,0.346]&0  \\
\\
\\
\\
\\
\\
\\
80\%, maxscale=0.3& 0.300 &[0,0.405]&0.776  &0.300&[0, 0.371]&0.474&0.300&[0,0.304]&0.056\\
\\
80\%, maxscale=0.5& 0.500 &[0,0.432]&0.004  &0.446&[0, 0.378]&0.001&0.369&[0,0.303]&0.004\\
\\
80\%, maxscale=1& 0.529 &[0,0.434]&0.001  \\
\\
\\
\\
\\
\\
\\
100\%, maxscale=0.3& 0.300 &[0,0.356]&0.370  &0.300&[0,0.287]&0.030&0.300&0.223&0\\
\\
100\%, maxscale=0.5& 0.340 &[0,0.357]&0.096  &0.334&[0,0.288]&0.007&0.323&0.223&0\\

\\
\\
\\
\\
\\
\\
\\
\\
\\
\textbf{$r=3$}&& $n=500$ &&& $n=1,000$&&& $n=2,000$\\\hline
\\
\\
\\
\\
Noise&$T_1$ real PD & CI & $p$-value& $T_1$ real PD & CI & $p$-value&  $T_1$ real PD & CI & $p$-value\\\hline
\\
\\
\\
\\
30\%, maxscale=0.3& 0.300 &[0,0.629]&1&0.300&[0,0.590]&1&0.300&[0,0.526]&1 \\
\\
30\%, maxscale=0.5& 0.500 &[0,0.682]&0.816 &0.500&[0,0.616] &0.508&0.500&[0,0.544]&0.164 \\
\\
30\%, maxscale=1& 0.625 &[0,0.687]&0.193 &0.767&[0,0.634] &0 \\
\\
\\
\\
\\
\\
\\
\\
70\%, maxscale=0.3& 0.300 &[0,0.304]&0.942  &0.300&	[0,0.540]&	1&0.300&[0,0.467]&0.997  \\
\\
70\%, maxscale=0.5&0.500&[0,0.625]&0.680& 0.500 &[0,0.567]&0.238  &0.500&	[0,0.477]&	0.028  \\
\\
70\%, maxscale=1&0.532&[0,0.635]&0.433& 0.572& [0,0.569]&0.047  &0.501&	[0,0.477]&	0.028 \\
\\
\\
\\
\\
\\
\\
80\%, maxscale=0.3& 0.300 &[0,0.302]&0.863  &0.300&[0,0.508]&1&0.300&[0,0.433]&0.989\\
\\
80\%, maxscale=0.5& 0.500 &[0,0.588]&0.470  &0.485&[0,0.524]&0.154\\
\\
80\%, maxscale=1& 0.594 &[0,0.597]&0.056  \\
\\
\\
\\
\\
\\
\\
100\%, maxscale=0.3& 0.300 &[0,0.300]&0.105  &0.300&[0,0.451]&0.999&0.300&[0,0.368]&0.623\\
\\
100\%, maxscale=0.5& 0.500 &[0,0.563]&0.283 &0.485&[0,0.471]&0.017&0.409&[0,0.377]&0.010\\
\\
100\%, maxscale=1& 0.527 &[0,0.567]&0.146 \\
\\
\label{table:OneCirc_noisyStat}
\end{tabular}
\end{center}
\footnotesize{Maximum statistic $T_1$ for the real $H_0$ persistence diagram and the simulated $H_0$ persistence diagrams of a sample $n$ of the one circle with radius $r$. The CI is a one-side confidence interval with $95\%$ confidence level. The $p$-value is also a one-side. Both the CI and the $p$-value are based on 1,000 simulated persistence diagrams.
The noise is an additive noise for some fraction $M$ of $n$, see text for more details.}\\
\normalsize

\subsection{Two Concentric Circles}
\begin{landscape}
\begin{figure}[h!]
\bc
\includegraphics[width=1.45in, height=1.45in]{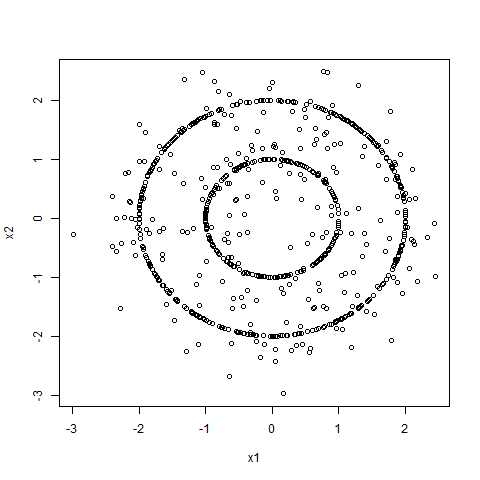} \hskip0.01truein
\includegraphics[width=1.45in, height=1.45in]{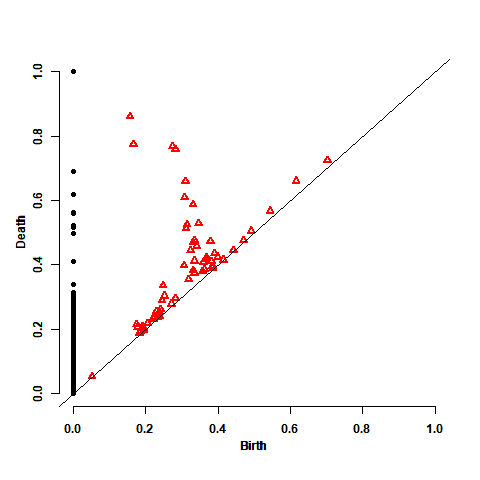} \hskip0.01truein
\includegraphics[width=1.45in, height=1.45in]{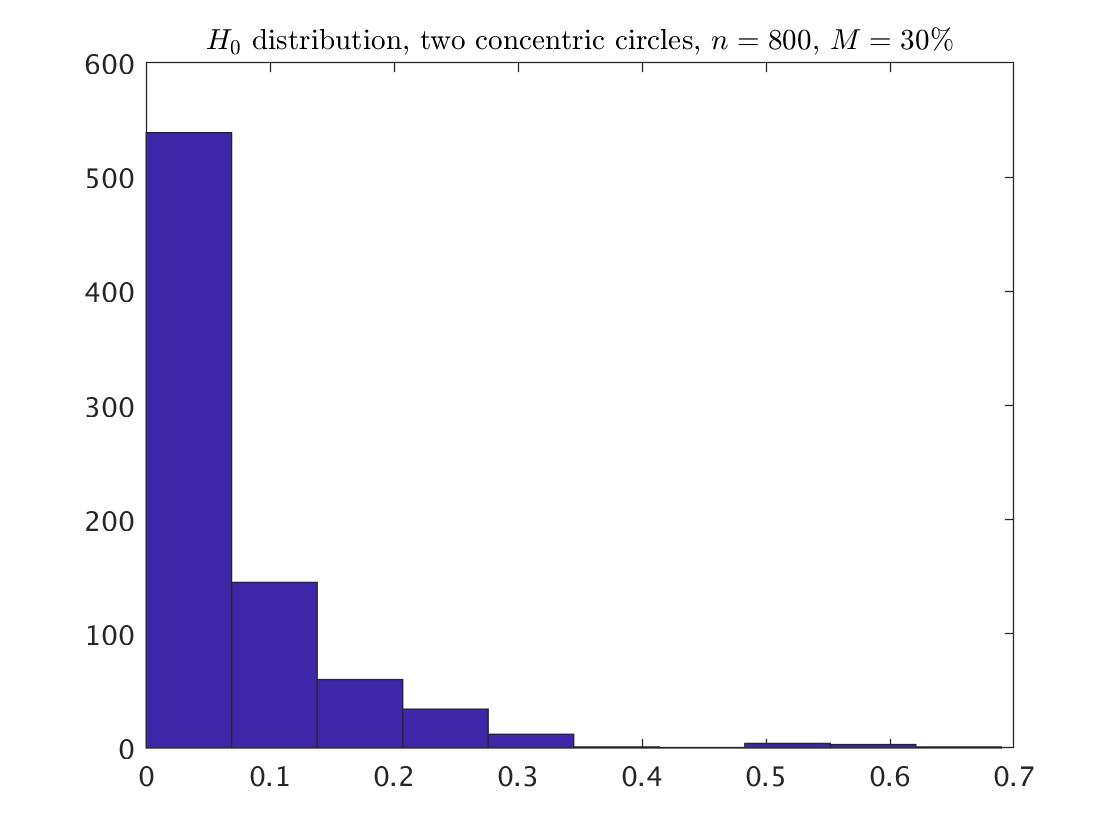} \hskip0.01truein
\includegraphics[width=1.45in, height=1.45in]{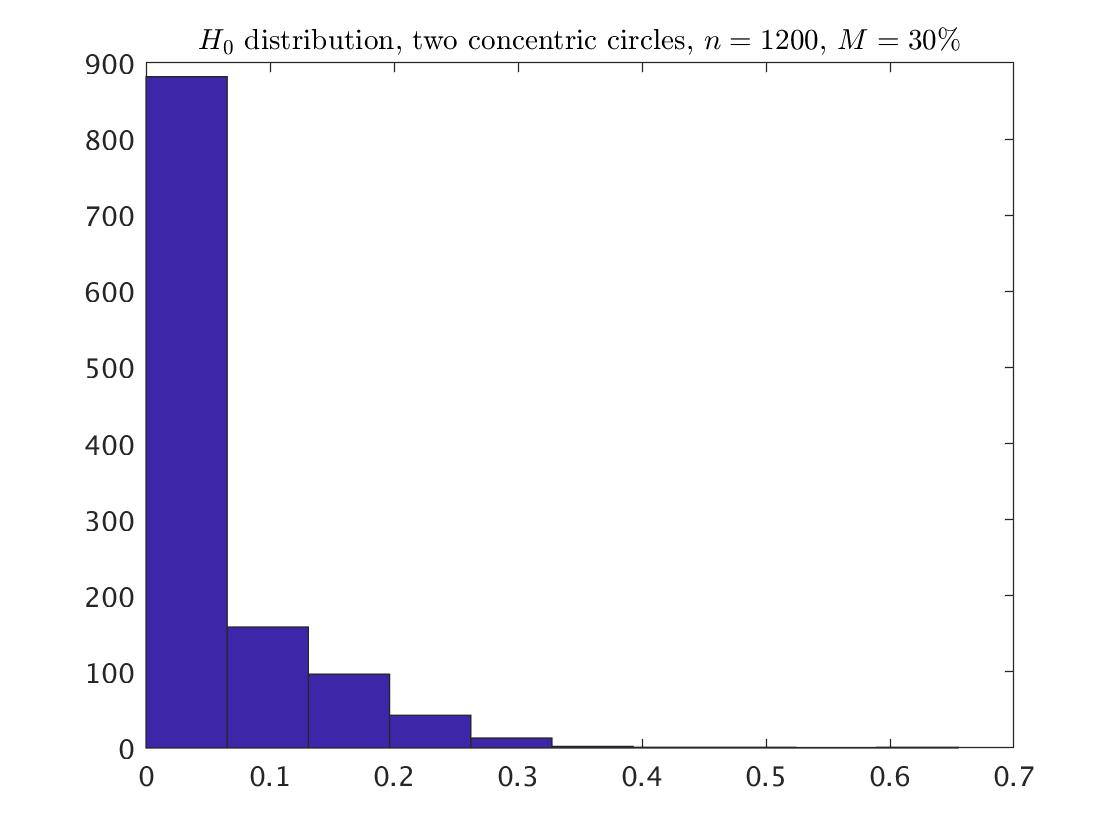} \hskip0.01truein
\includegraphics[width=1.45in, height=1.45in]{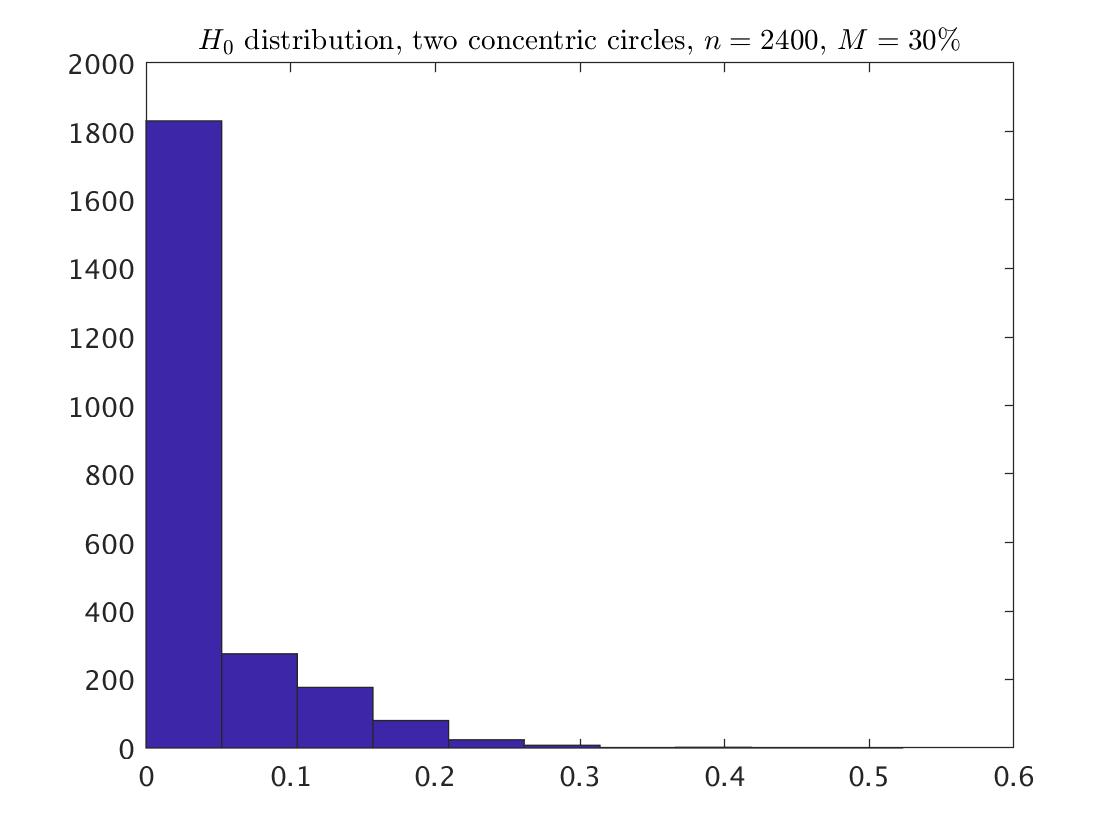} \hskip0.01truein

\includegraphics[width=1.45in, height=1.45in]{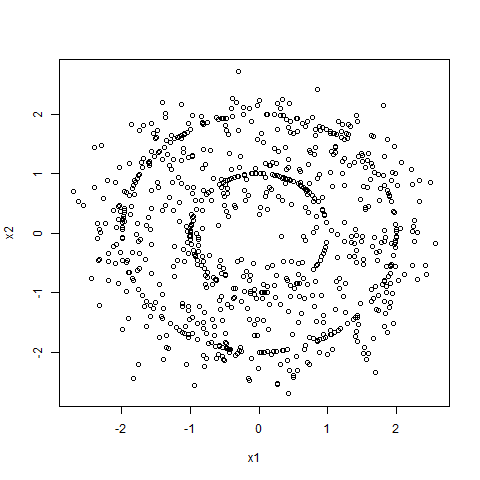} \hskip0.01truein
\includegraphics[width=1.45in, height=1.45in]{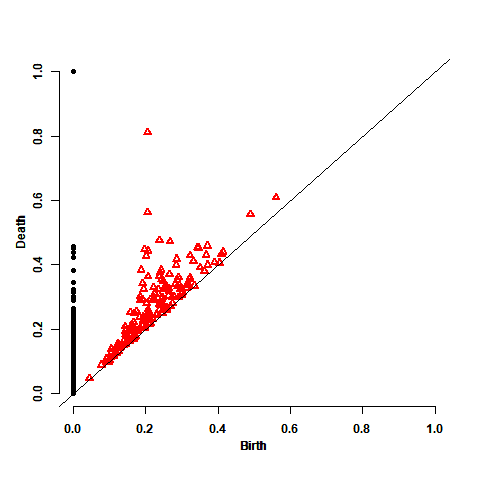} \hskip0.01truein
\includegraphics[width=1.45in, height=1.45in]{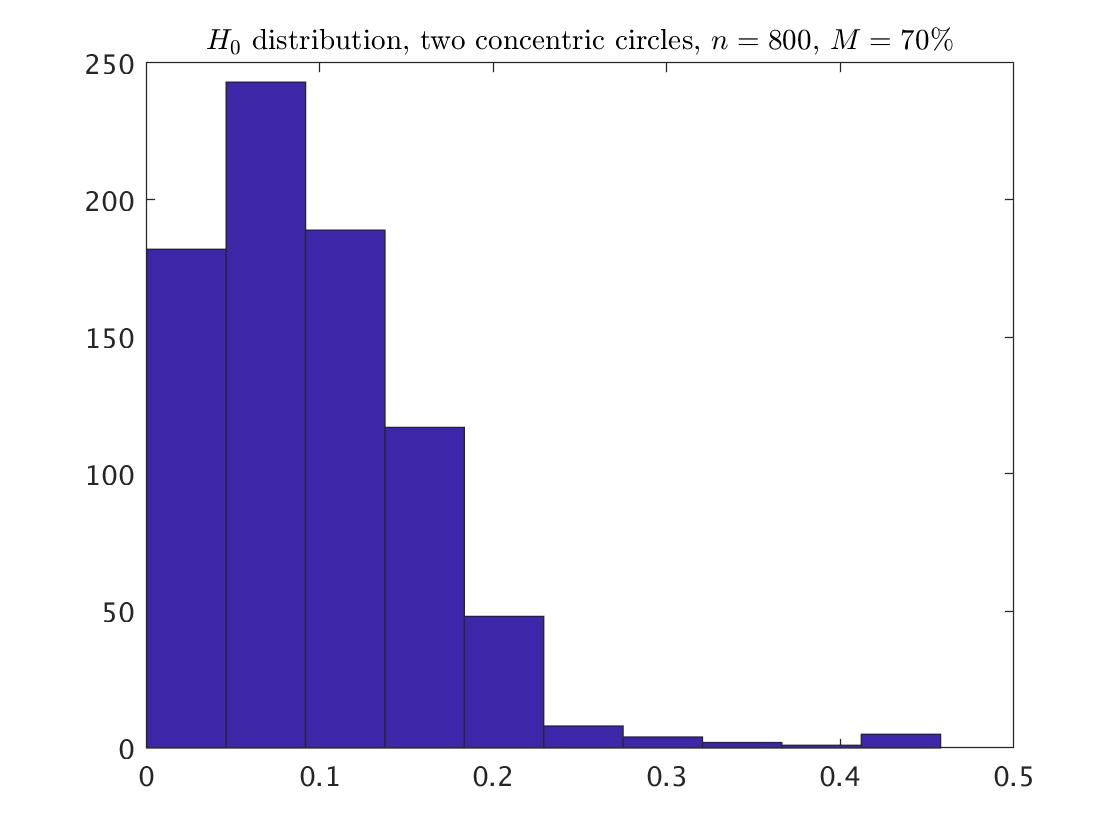} \hskip0.01truein
\includegraphics[width=1.45in, height=1.45in]{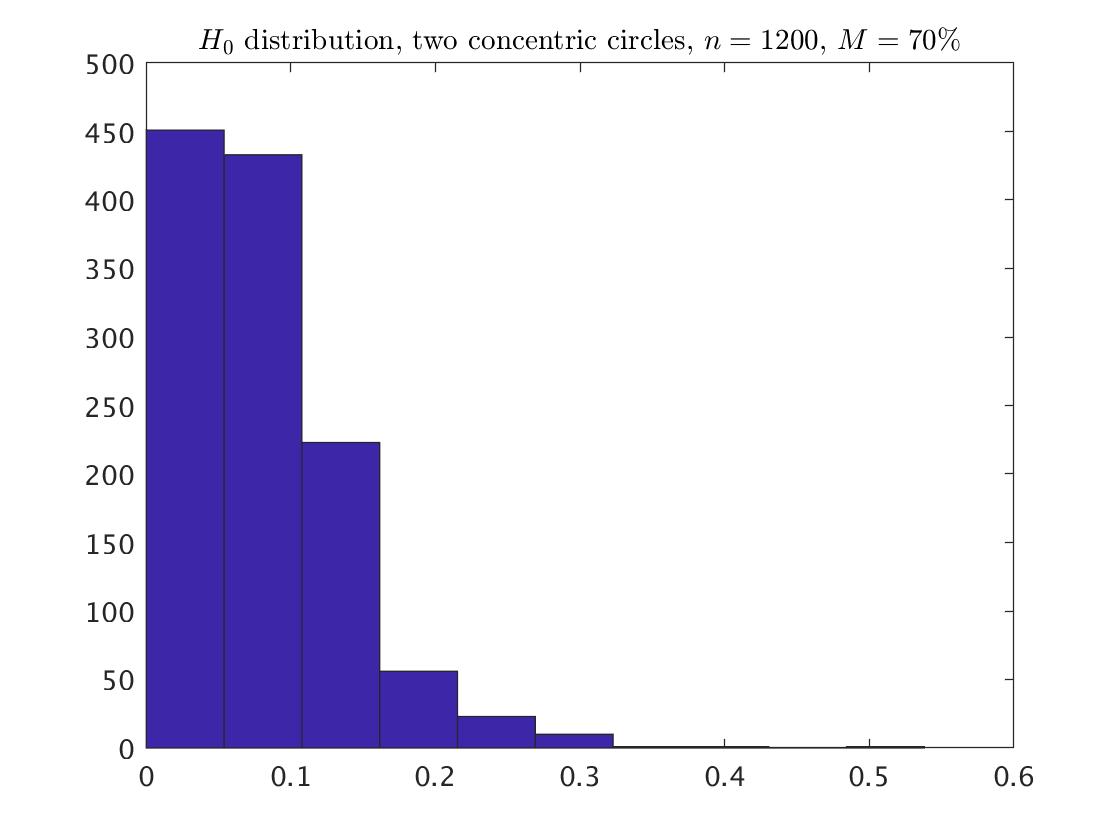} \hskip0.01truein
\includegraphics[width=1.45in, height=1.45in]{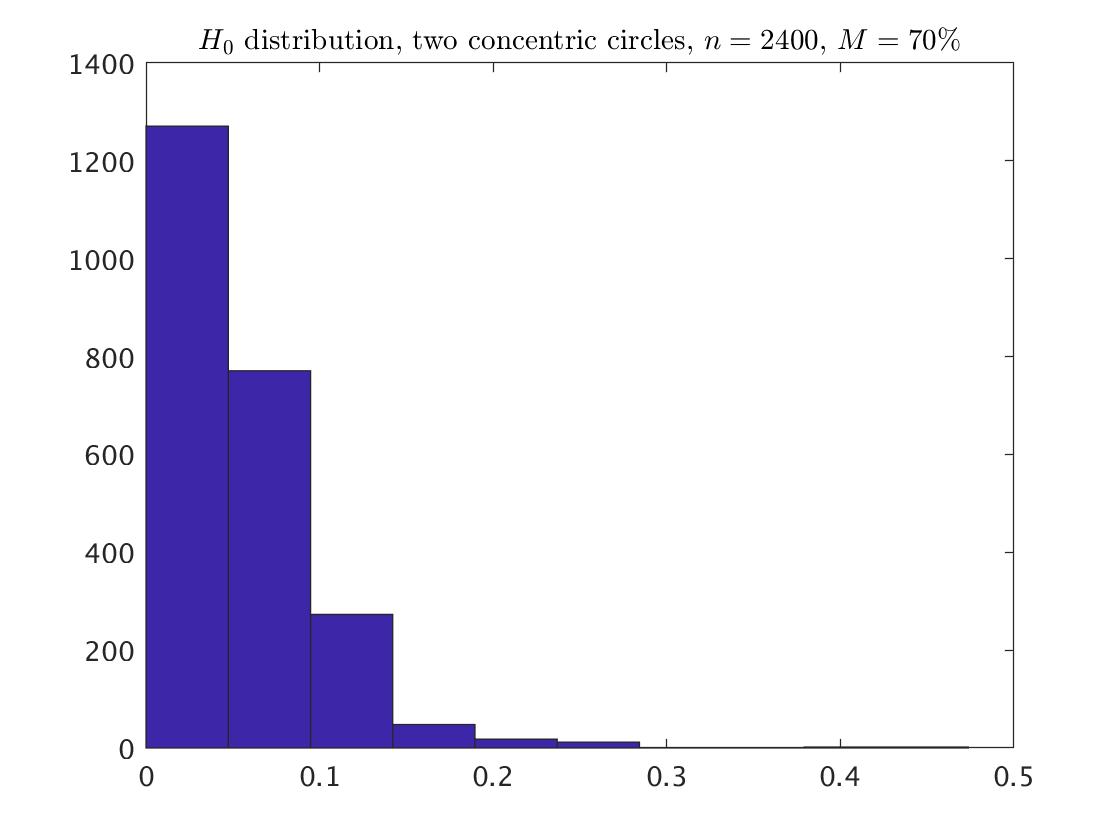} \hskip0.01truein

\includegraphics[width=1.45in, height=1.45in]{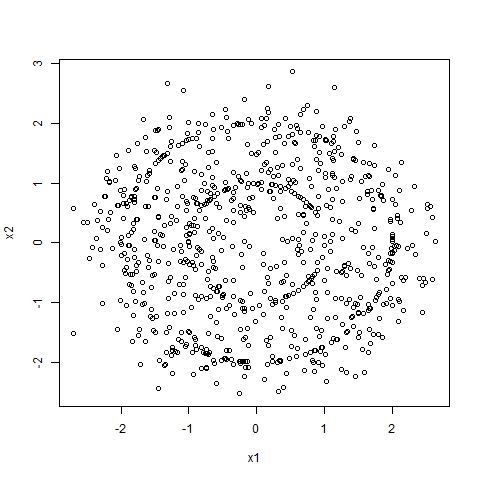} \hskip0.01truein
\includegraphics[width=1.45in, height=1.45in]{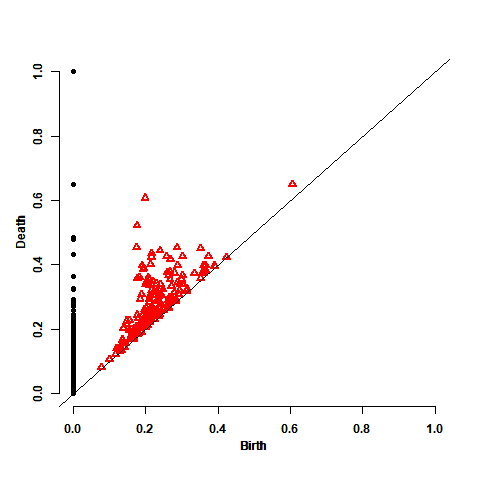} \hskip0.01truein
\includegraphics[width=1.45in, height=1.45in]{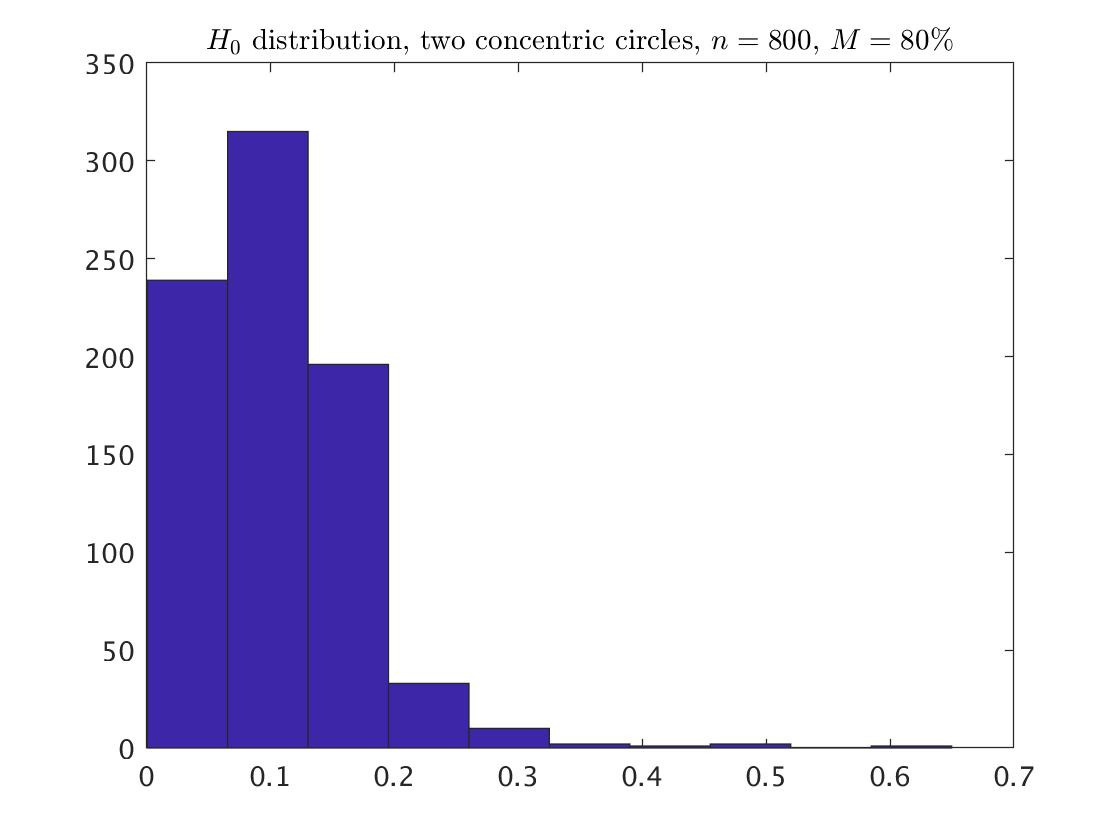} \hskip0.01truein
\includegraphics[width=1.45in, height=1.45in]{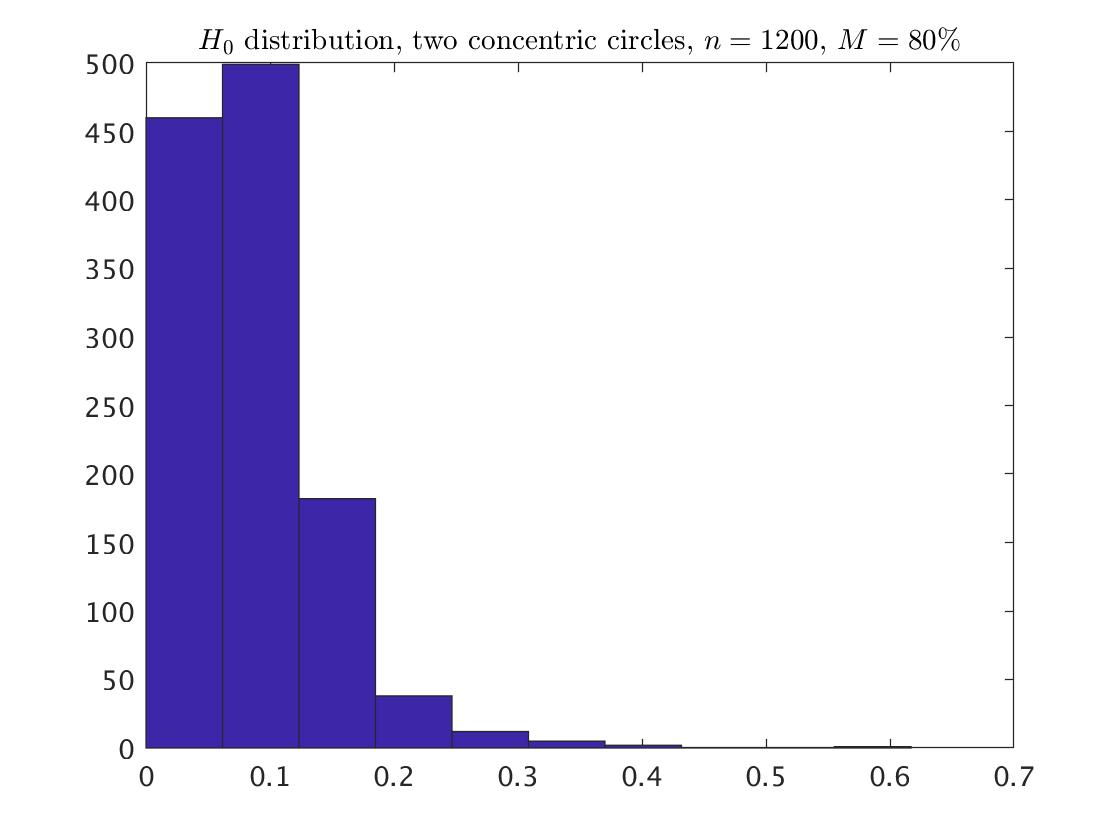} \hskip0.01truein
\includegraphics[width=1.45in, height=1.45in]{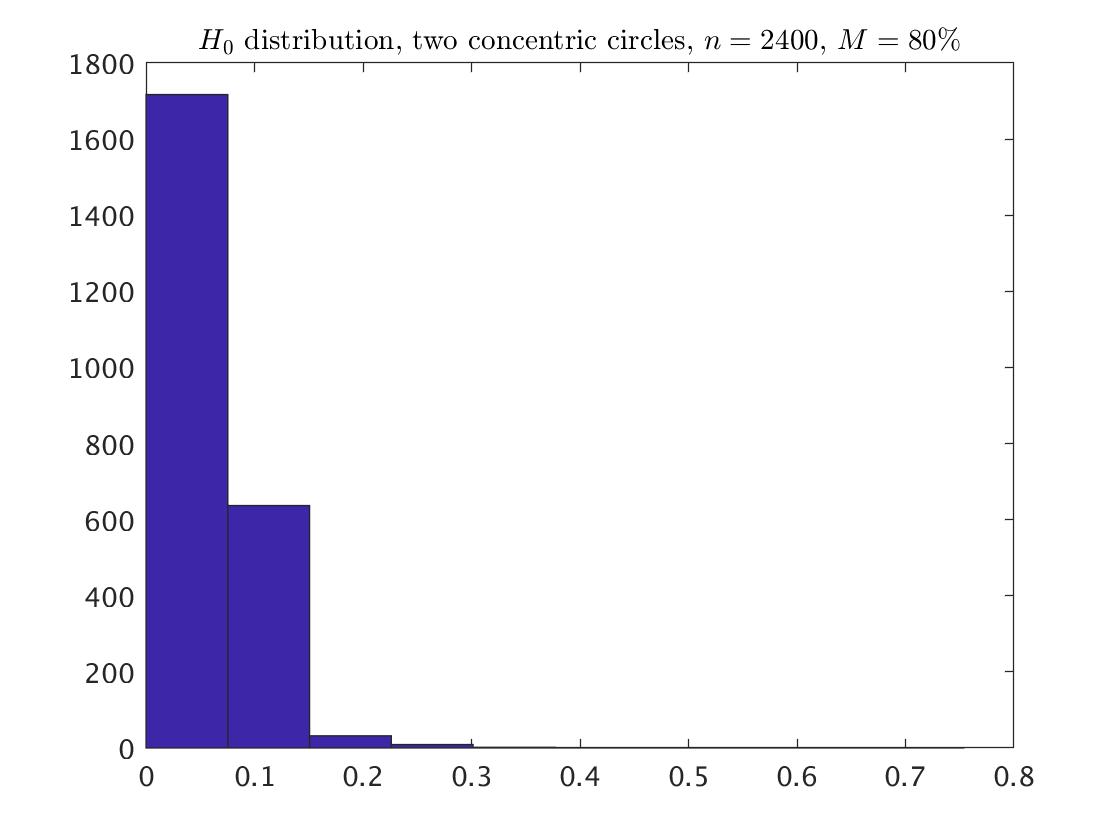} \hskip0.01truein

\includegraphics[width=1.45in, height=1.45in]{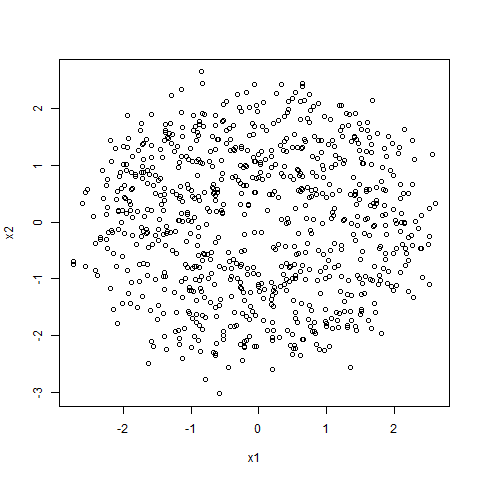} \hskip0.01truein
\includegraphics[width=1.45in, height=1.45in]{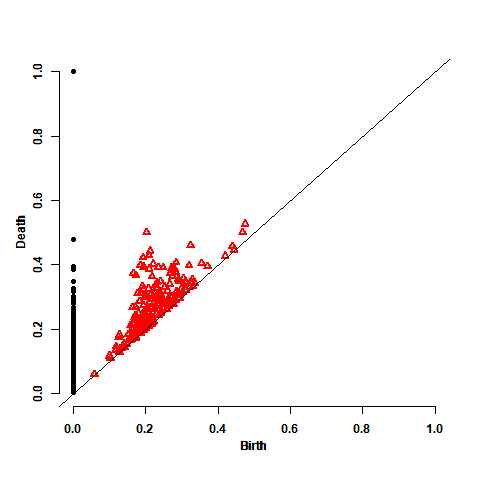} \hskip0.01truein
\includegraphics[width=1.45in, height=1.45in]{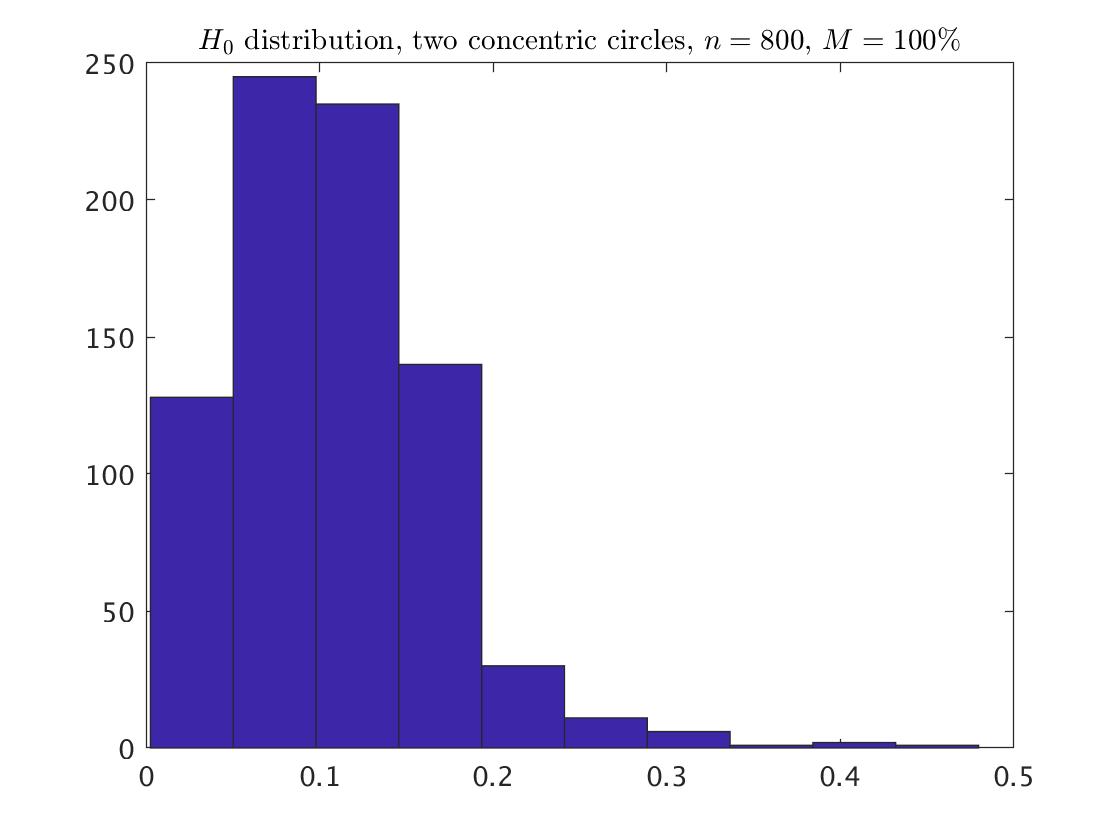} \hskip0.01truein
\includegraphics[width=1.45in, height=1.45in]{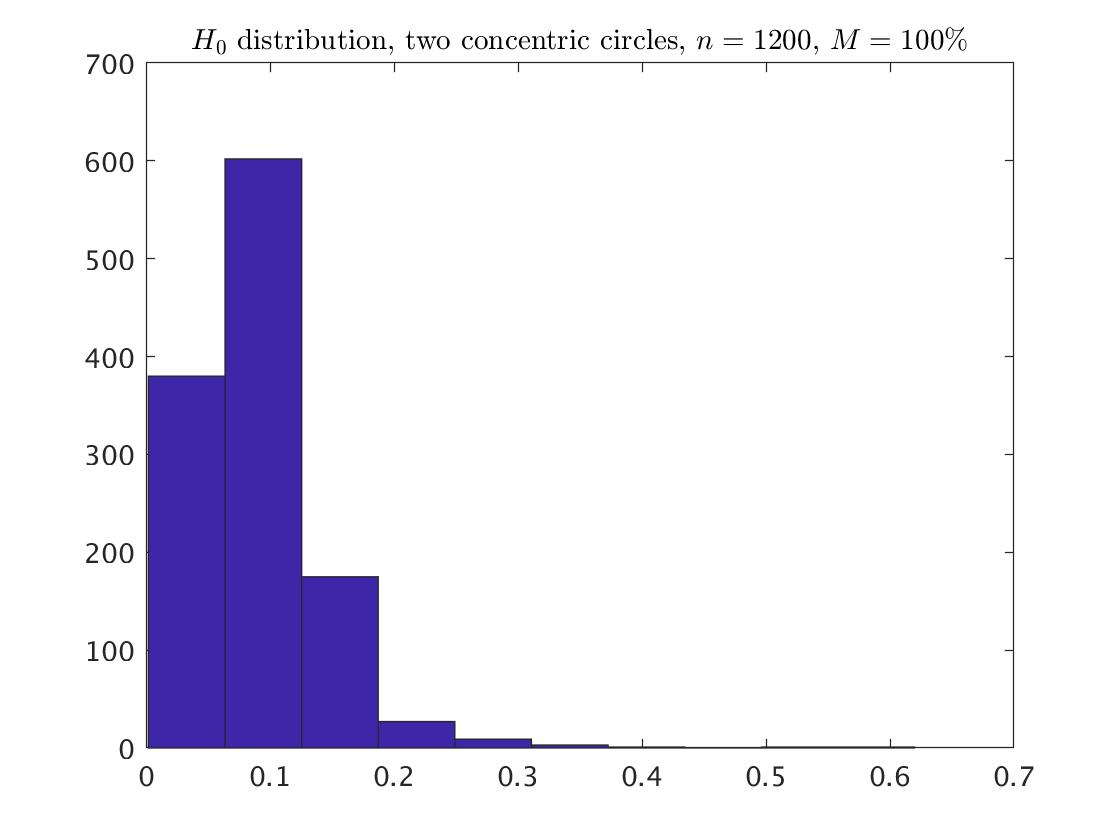} \hskip0.01truein
\includegraphics[width=1.45in, height=1.45in]{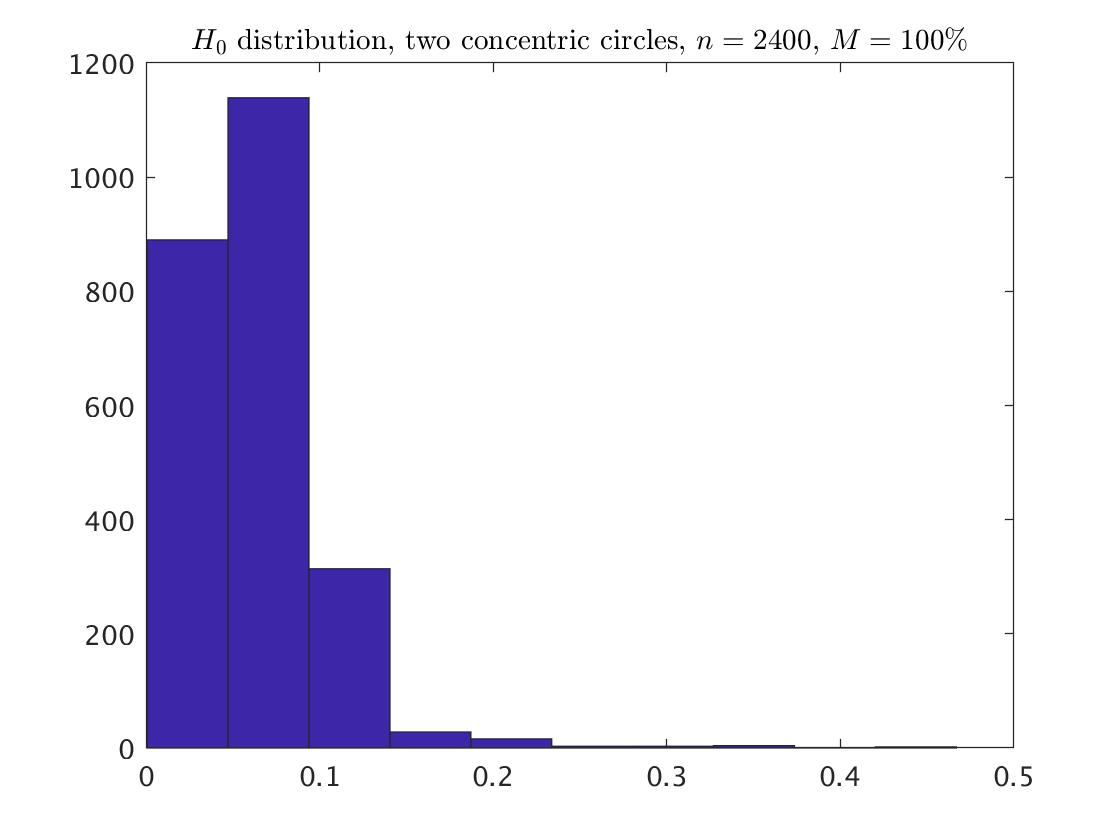} \hskip0.01truein

\ec
\caption{\footnotesize
 Two concentric circles with noise of $M\%$ of the sample size $n$. From top to bottom: $M = 30\%, 70\%, 80\%, 100\%$. Each row contains the circle with $n=800$ points, and to its right the histograms of the $H_0$ persistence diagrams for $n=800, 1,200, 2,400$.
}
\label{fig:Twocircles_M_noise}
\end{figure}
\end{landscape}

We examine now the influence of noise in the example of the two concentric circles, for $n=800, 1,200, 2,400$, and maxscale= 0.3, 0.5, 1.
Fig. \ref{fig:Twocircles_M_noise} describes for each value of $M$ the data, the persistence diagram  based on maxscale=1, and the distribution of the $H_0$ points (without the point at infinity) based on maxscale=1.
For examining the influence of noise on the length of the bars, we again, as in the previous example, classified the bars into two groups of small and long bars, depends on the $c_{max}$ of the relevant $n$. The results, depending on noise and maxscale, are presented in Table \ref{table:circlesLong}. As in the previous example, the proportion of the long bars increases as the noise increases, and the difference between the various values of maxscale for a given level of noise is in the lengths of the longest bars, where the proportion of the long bars is the same. Table \ref{table:circlesPercntiles} summarizes for each maxscale and a given noise, the $95$-th, $99$-th, and the $100$-th percentiles of the $H_0$ lengths. We can see that the difference between the various values of maxscale is usually at the $100$-th or the $99$-th percentile of the $H_0$ lengths.

\normalsize

\begin{center}
\fontsize{8.5}{0.9}\selectfont
\captionof{table}{Two concentric circles - classification of $H_0$ points}
\begin{tabular}{l|lc|cc|cc}
\\
\\
&$n=800$$^a$ &&$n=1,200$$^b$&&$n=2,400$$^c$\\\hline
\\
Noise&short& long&short& long\\\hline
\\
\\
\\
\\
30\%& 0.850&0.152&0.801& 0.199&0.800&0.200\\
\\
\\
\\
\\
70\%& 0.757&0.244 &0.678&0.322&0.678&0.322\\
\\
\\
\\
\\
80\%& 0.707&0.294& 0.649&0.351&0.637&0.363 \\
\\
\\
\\
\\
100\%& 0.695&0.306&	0.614& 0.386&0.603&0.397\\
\label{table:circlesLong}
\end{tabular}
\end{center}
\footnotesize{Classification of the $H_0$ death times into short and long bars for two concentric circles. The classification is based on one persistence diagram, with maxscale$\ge$ 1. $^a$ $c_{max}=0.134$, $^b$ $c_{max}=0.097$, $^c$ $c_{max}=0.067$.}\\

\normalsize
\begin{center}
\fontsize{8.5}{0.9}\selectfont
\captionof{table}{Two concentric circles - percentiles of $H_0$ points}
\begin{tabular}{l|l|c|c|c|c}
\\
\\
Noise&maxscale&percentile $\%$&$n=800$$^a$ &$n=1,200$$^b$ &$n=2,400$$^c$\\\hline
\\
\\
\\
\\
30\%&0.3& 95&0.228&	0.198&	0.156\\
\\
\\
\\
\\
&& 99&0.300&	0.282&	0.226\\
\\
\\
\\
\\
&& 100&0.300&	0.300&0.300\\
\\
\\
\\
\\
& 0.5& 100&0.500&0.500&0.500\\
\\
\\
\\
\\
& 1& 100&0.690&	0.655&	0.523\\
\\
\\
\\
\\
& &&&\\\hline
\\
\\
\\
70\%&0.3& 95&0.198&	0.186&	0.128\\
\\
\\
\\
\\
&& 99&0.300&	0.276&	0.221\\
\\
\\
\\
\\
&& 100&0.300&	0.300&	0.300\\
\\
\\
\\
\\
& 0.5& 99&0.319&	0.276&	0.221\\
\\
\\
\\
\\
&& 100&0.458&	0.500&	0.474\\
\\
\\
\\
\\
& 1& 100&0.458&	0.538&	0.474\\
\\
\\
\\
\\
& &&&\\\hline
\\
\\
\\
80\%& 0.3& 95&0.202&	0.183&	0.123\\
\\
\\
\\
\\
&& 99&0.294&	0.300&	0.192\\
\\
\\
\\
\\
&& 100&0.300&	0.300&	0.300\\
\\
\\
\\
\\
&0.5& 100&0.500&	0.500&	0.500\\
\\
\\
\\
\\
&1& 100&0.650&	0.617&	0.754\\
& &&&\\\hline
\\
\\
\\
100\%& 0.3& 95&0.203&	0.168&	0.116\\
\\
\\
\\
\\
&& 99&0.300&	0.276&	0.199\\
\\
\\
\\
\\
&& 100&0.300&	0.300&	0.300\\
\\
\\
\\
\\
& 0.5& 99&0.311&	0.276&	0.199\\
\\
\\
\\
\\
& & 100&0.480&	0.500&	0.467\\
\\
\\
\\
\\
&1& 100&0.480&	0.620&	0.467\\
\\
\\
\\
\\
\label{table:circlesPercntiles}
\end{tabular}
\end{center}
\footnotesize{Percentiles of the $H_0$ lengths. The $H_0$ lengths are based on one persistence diagram, with maxscale$\ge$ 0.3.}\\
\normalsize

\subsection{3-Torus}
\begin{landscape}
\begin{figure}[h!]
\bc
\includegraphics[width=1.45in, height=1.45in]{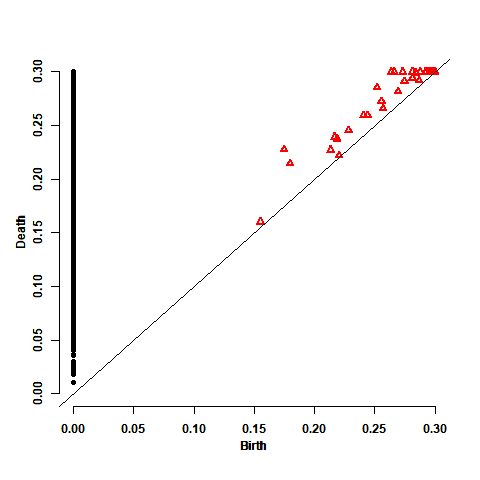} \hskip0.01truein
\includegraphics[width=1.45in, height=1.45in]{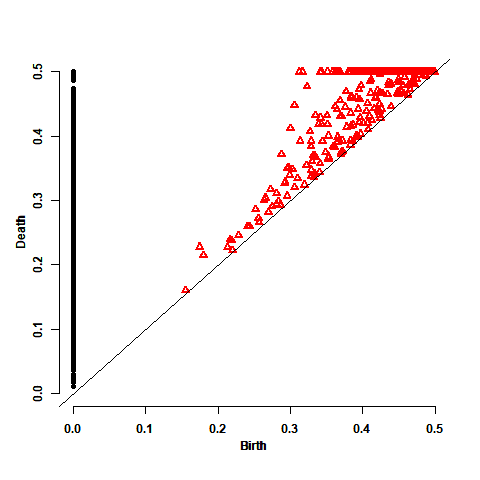} \hskip0.01truein
\includegraphics[width=1.45in, height=1.45in]{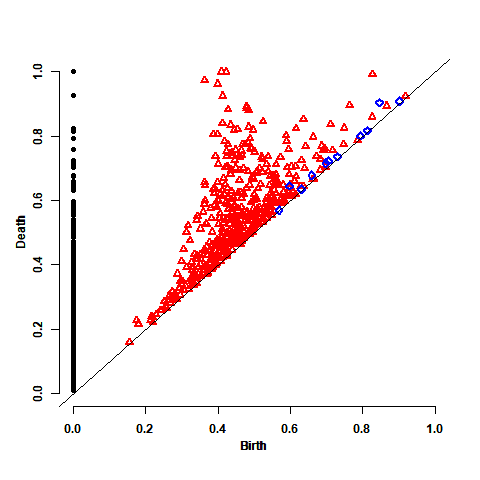} \hskip0.01truein
\includegraphics[width=1.45in, height=1.45in]{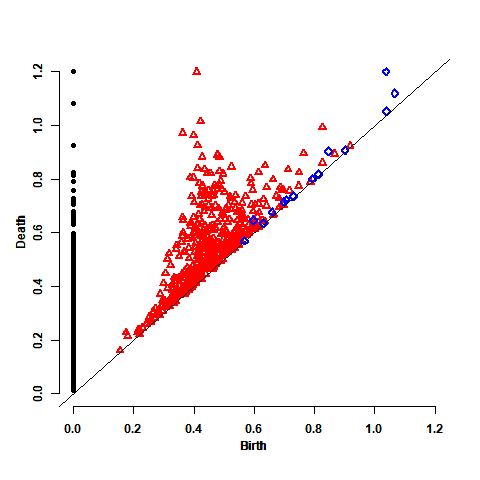} \hskip0.01truein
\includegraphics[width=1.45in, height=1.45in]{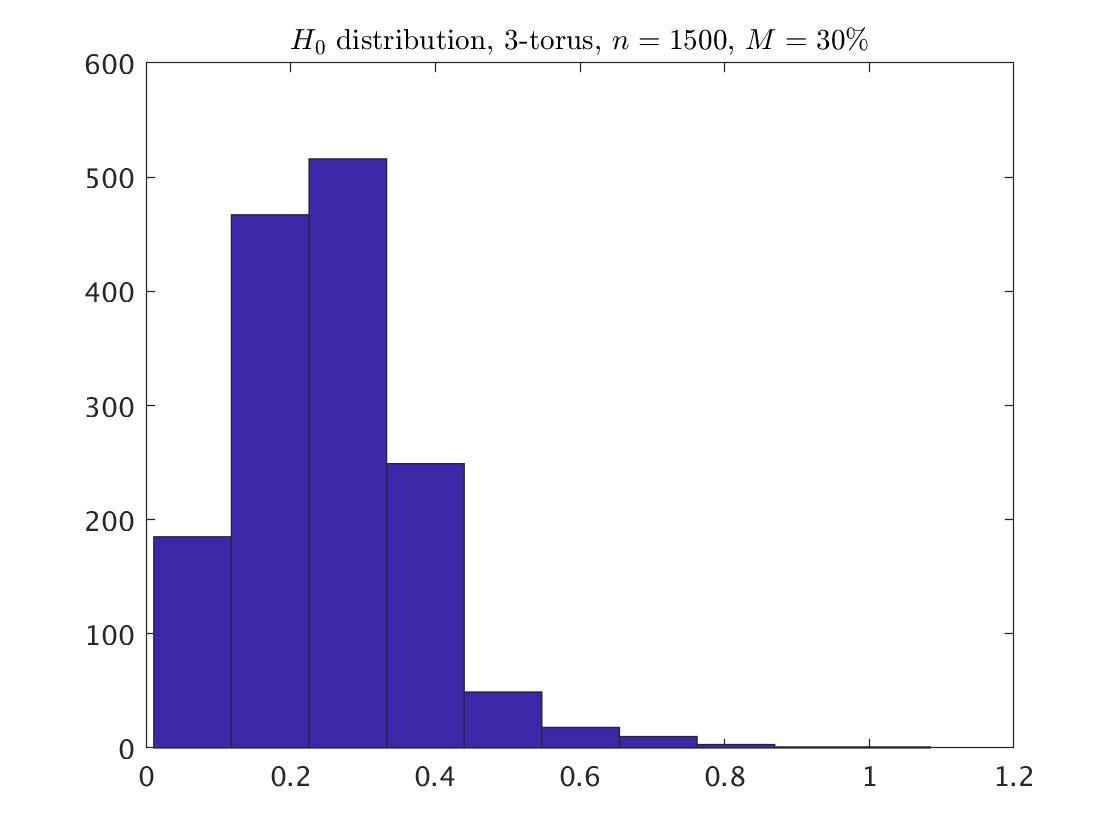} \hskip0.01truein

\includegraphics[width=1.45in, height=1.45in]{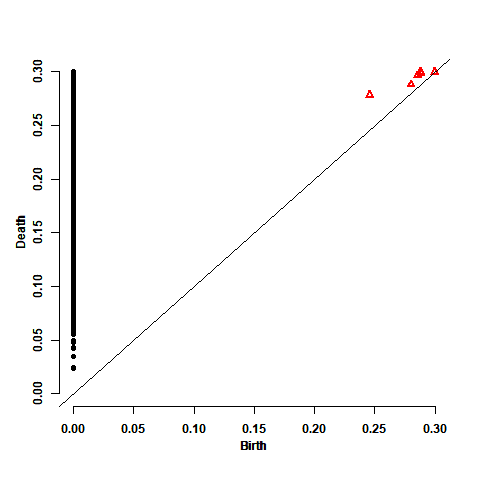} \hskip0.01truein
\includegraphics[width=1.45in, height=1.45in]{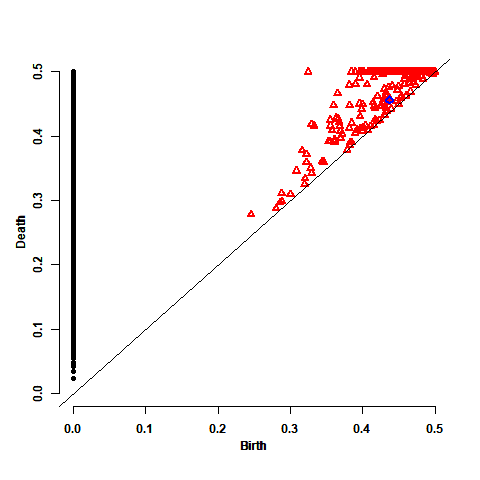} \hskip0.01truein
\includegraphics[width=1.45in, height=1.45in]{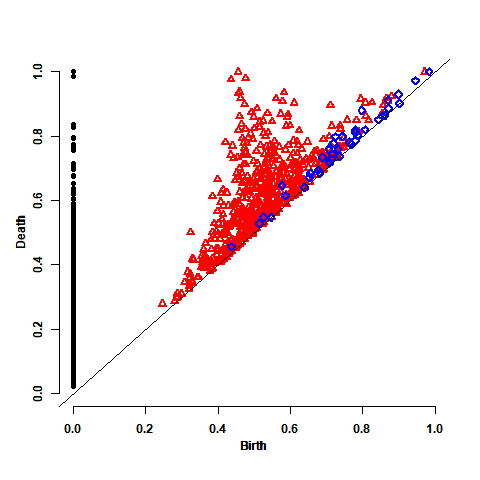} \hskip0.01truein
\includegraphics[width=1.45in, height=1.45in]{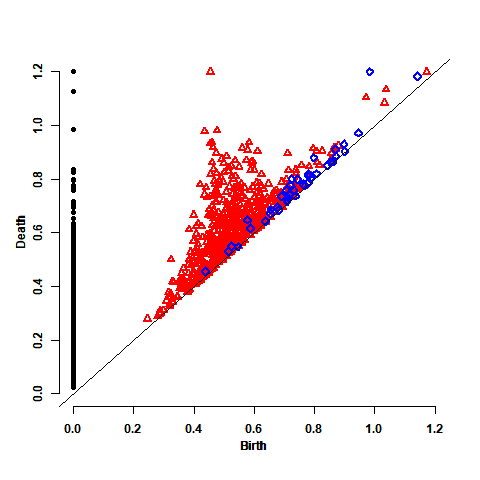} \hskip0.01truein
\includegraphics[width=1.45in, height=1.45in]{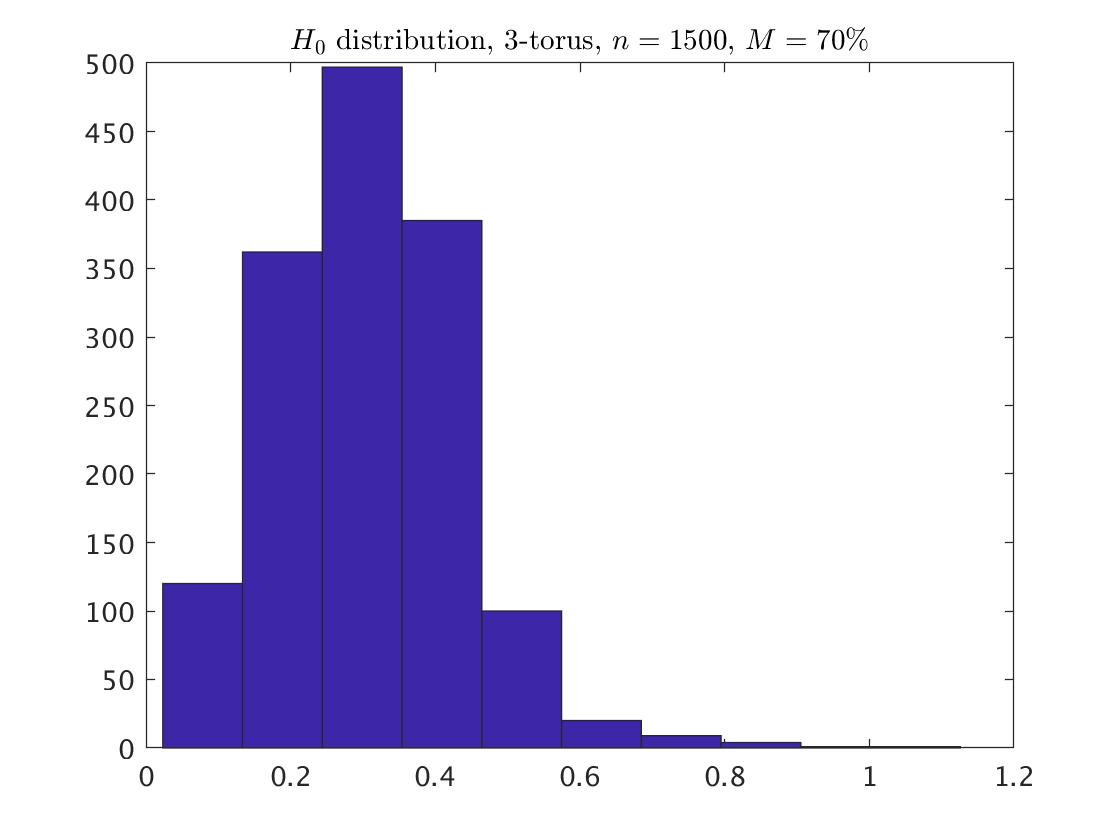} \hskip0.01truein

\includegraphics[width=1.45in, height=1.45in]{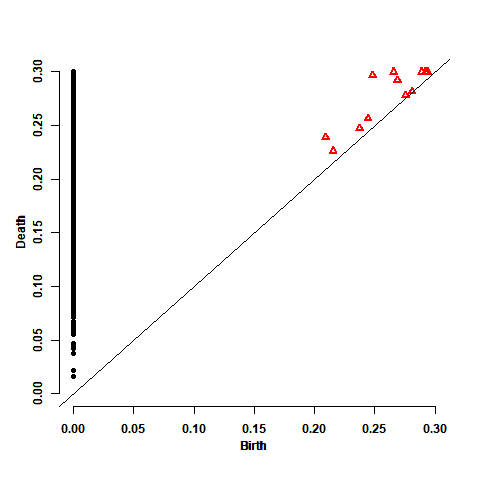} \hskip0.01truein
\includegraphics[width=1.45in, height=1.45in]{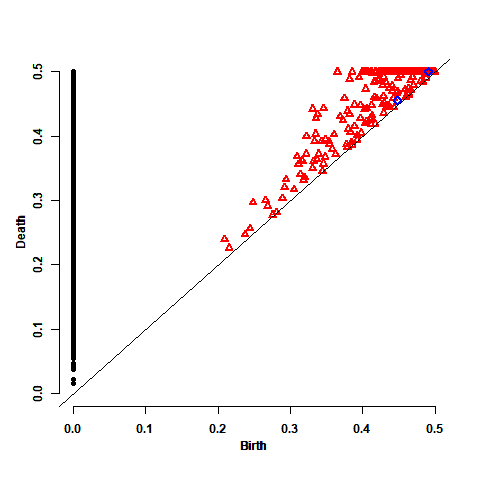} \hskip0.01truein
\includegraphics[width=1.45in, height=1.45in]{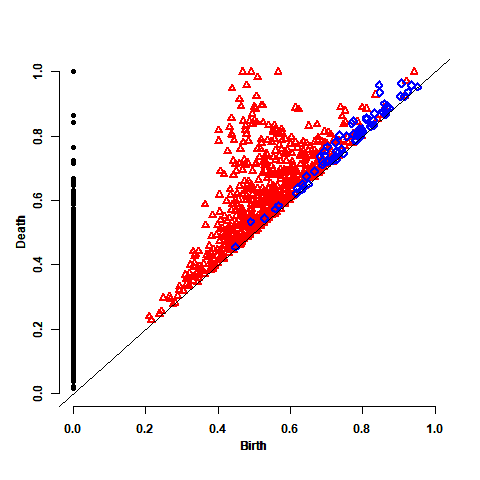} \hskip0.01truein
\includegraphics[width=1.45in, height=1.45in]{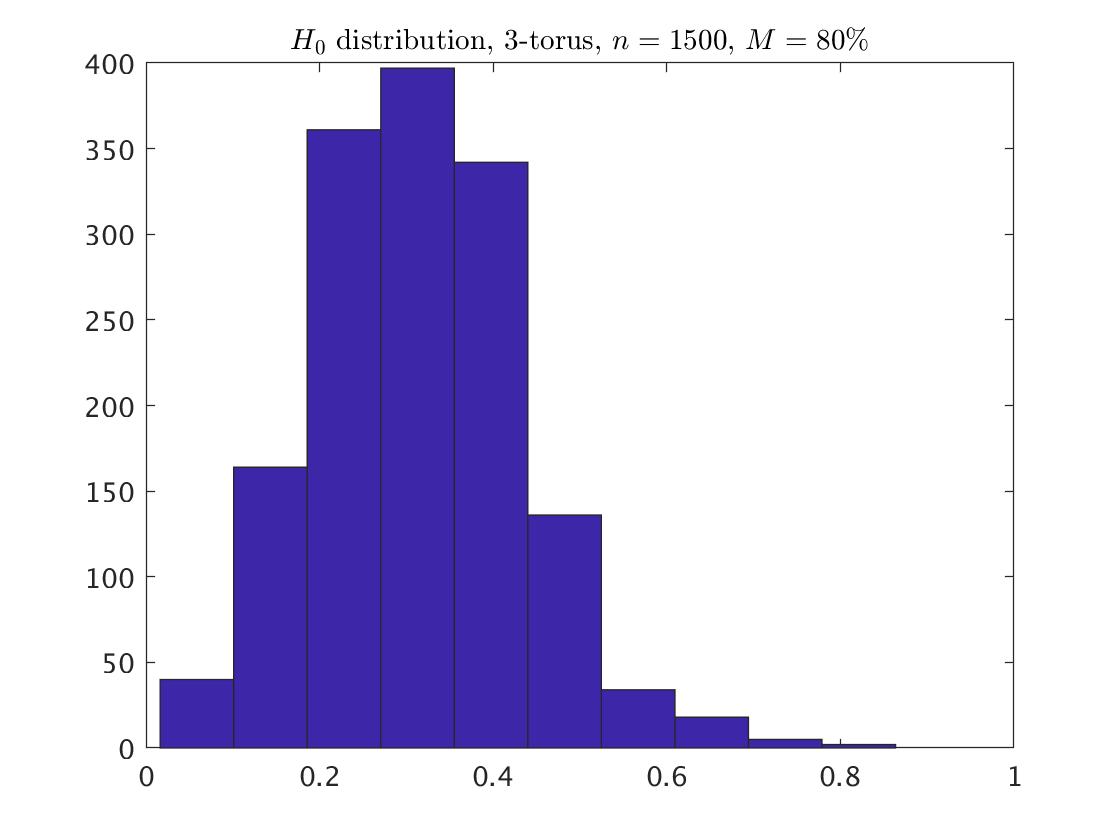} \hskip0.01truein

\includegraphics[width=1.45in, height=1.45in]{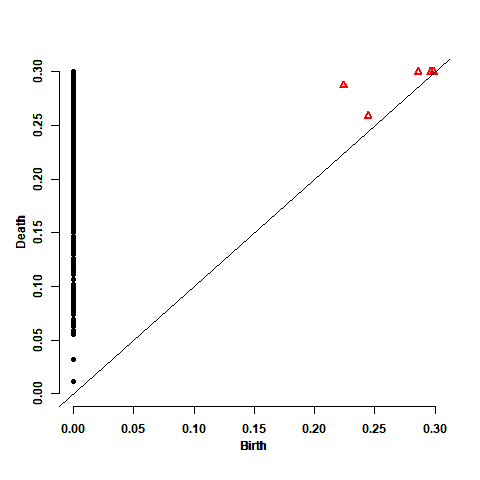} \hskip0.01truein
\includegraphics[width=1.45in, height=1.45in]{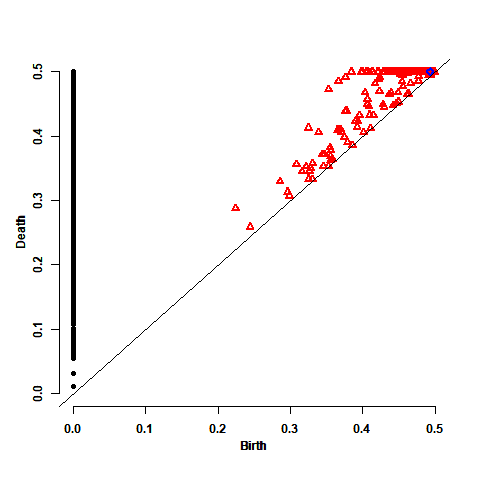} \hskip0.01truein
\includegraphics[width=1.45in, height=1.45in]{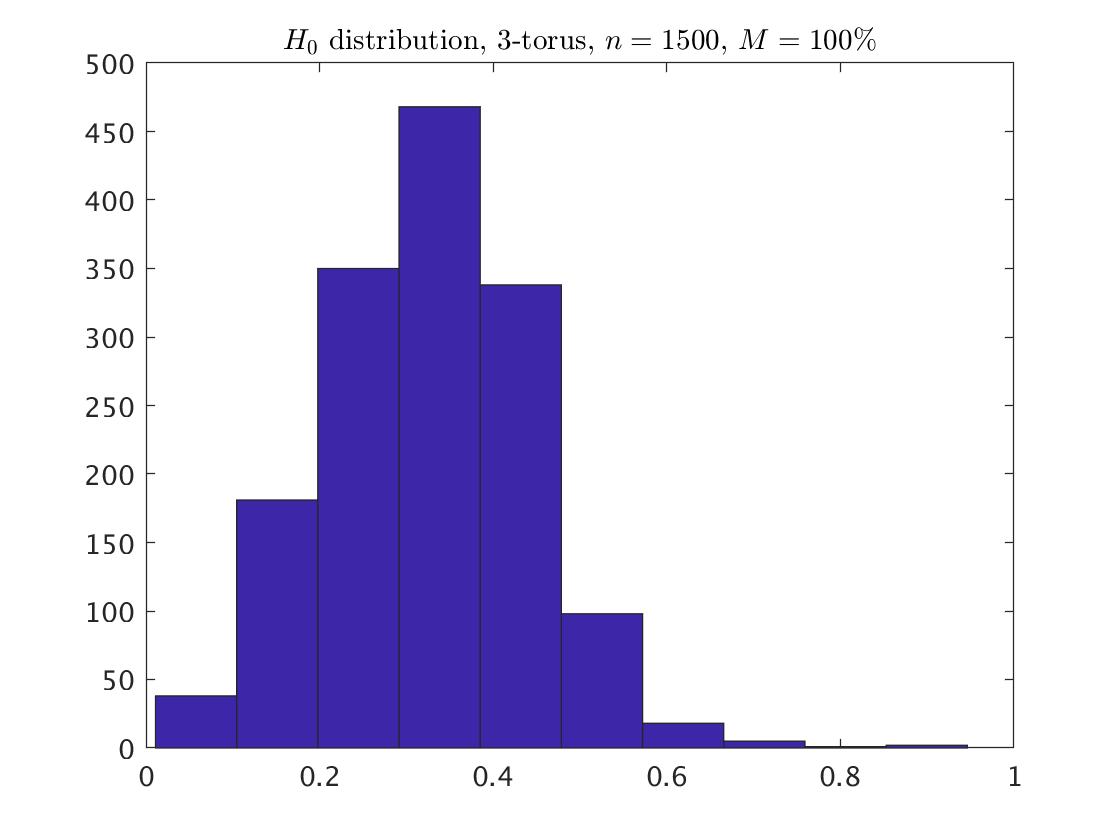} \hskip0.01truein
\includegraphics[width=1.45in, height=1.45in]{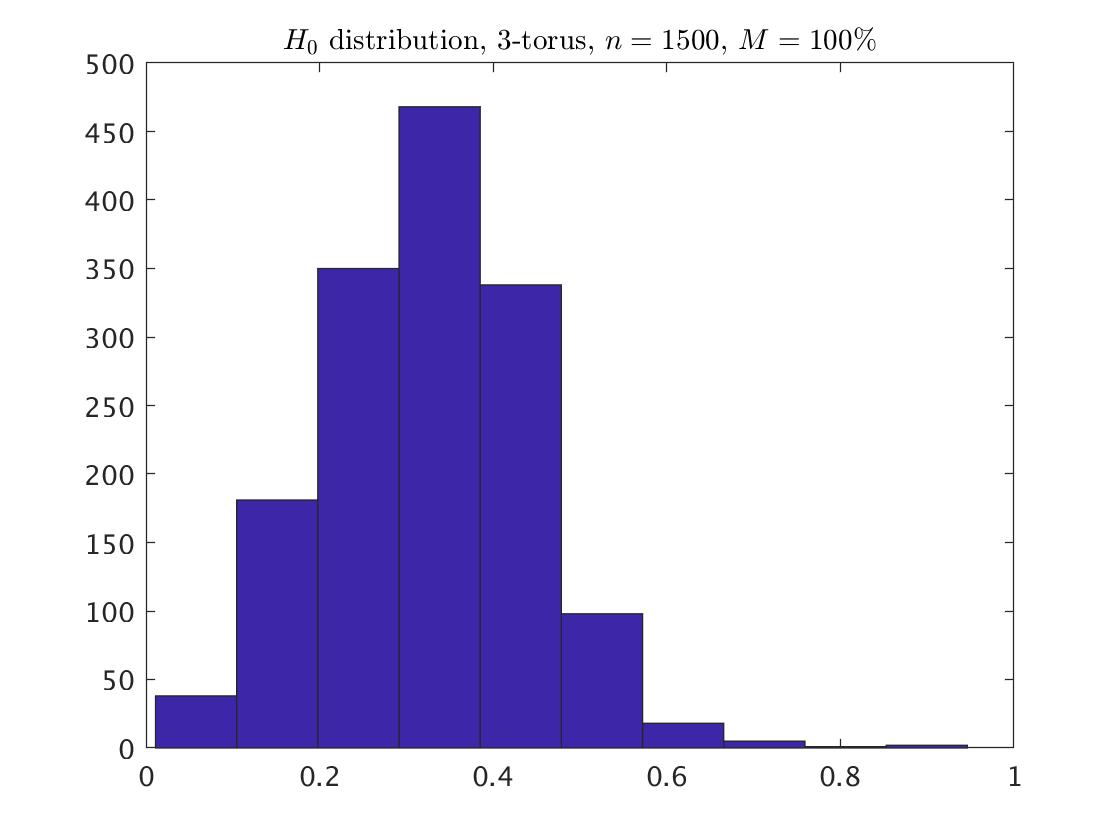} \hskip0.01truein
\ec
\caption{\footnotesize
 3-Torus with noise of $M\%$ of the sample size $n=1,500$. From top to bottom: $M = 30\%, 70\%, 80\%, 100\%$. Each row contains the persistence diagram based on maxscale=0.3, 0.5, and 1, where the additional persistence diagram in $M = 30\%, 70\%$ is based on maxscale=1.2.
}
\label{fig:torus_M_noise}
\end{figure}
\end{landscape}

\begin{landscape}
\begin{figure}[h!]
\bc
\includegraphics[width=1.45in, height=1.45in]{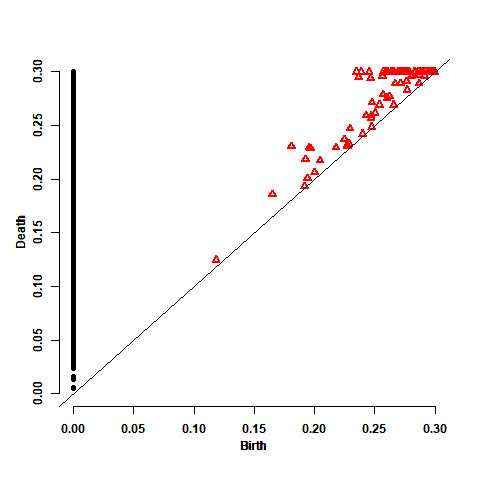} \hskip0.01truein
\includegraphics[width=1.45in, height=1.45in]{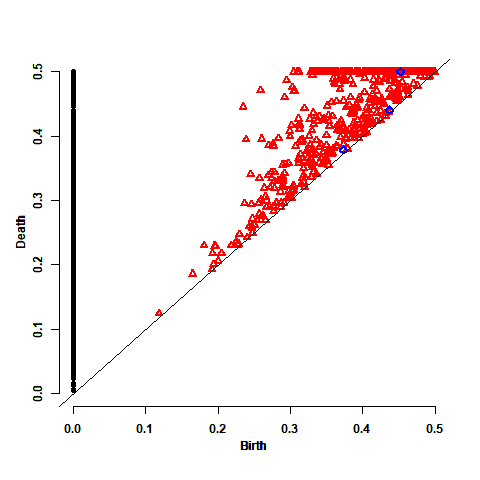} \hskip0.01truein
\includegraphics[width=1.45in, height=1.45in]{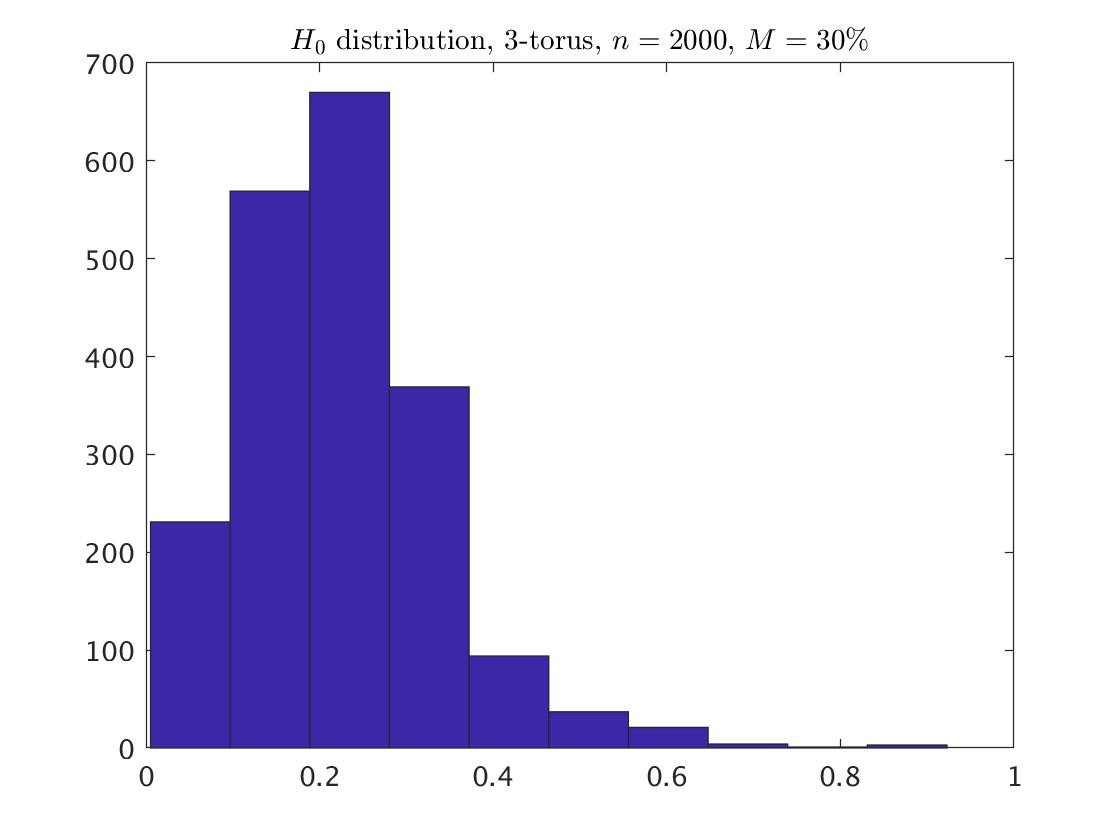} \hskip0.01truein
\includegraphics[width=1.45in, height=1.45in]{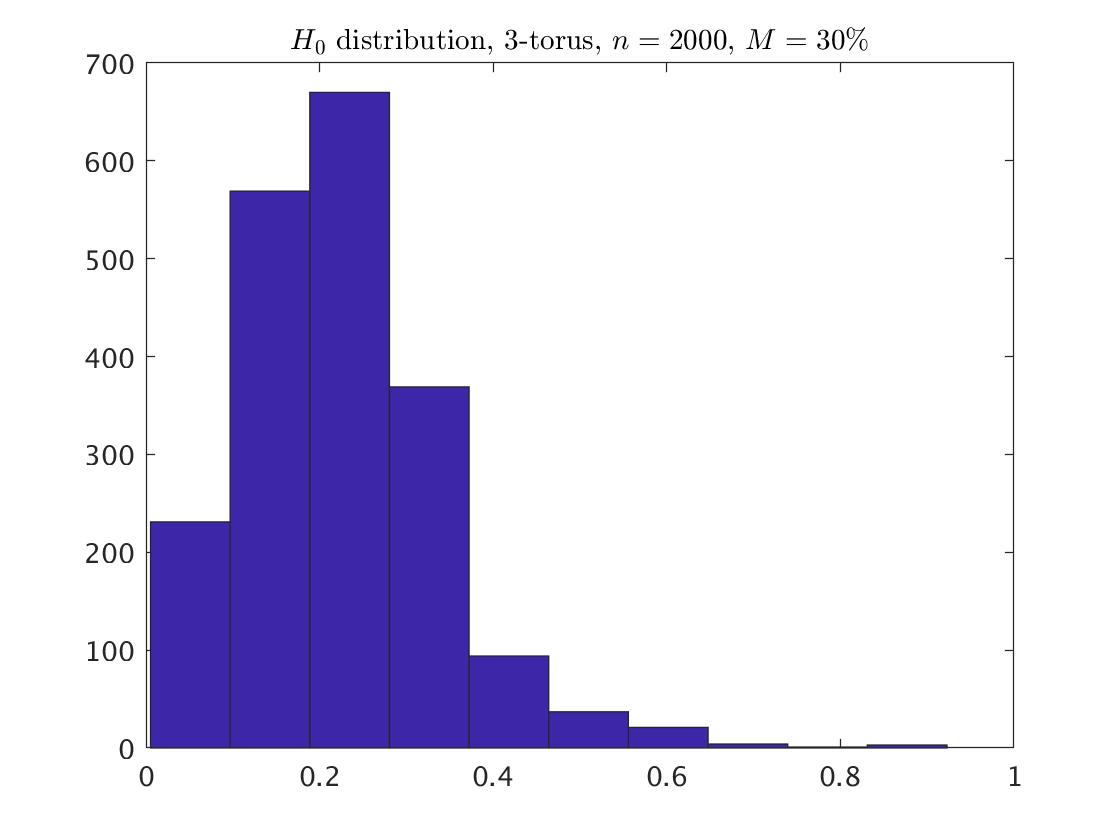} \hskip0.01truein

\includegraphics[width=1.45in, height=1.45in]{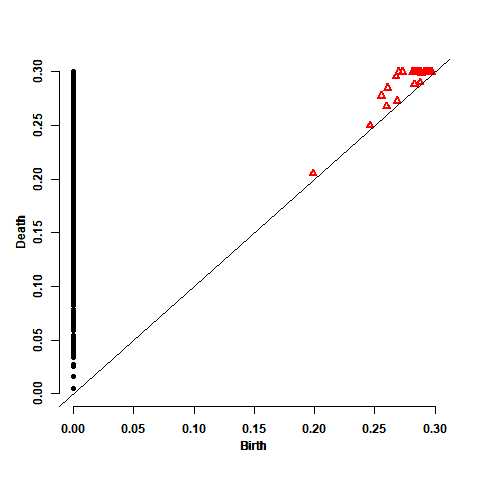} \hskip0.01truein
\includegraphics[width=1.45in, height=1.45in]{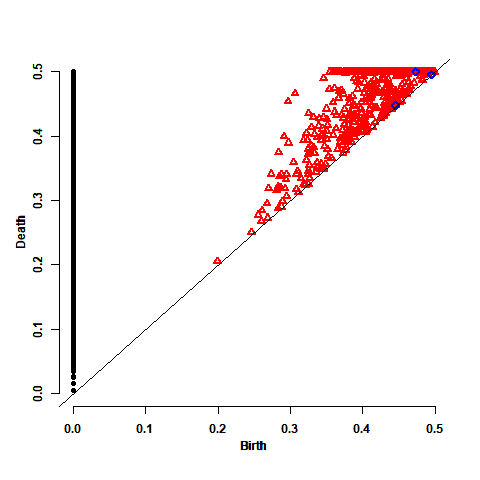} \hskip0.01truein
\includegraphics[width=1.45in, height=1.45in]{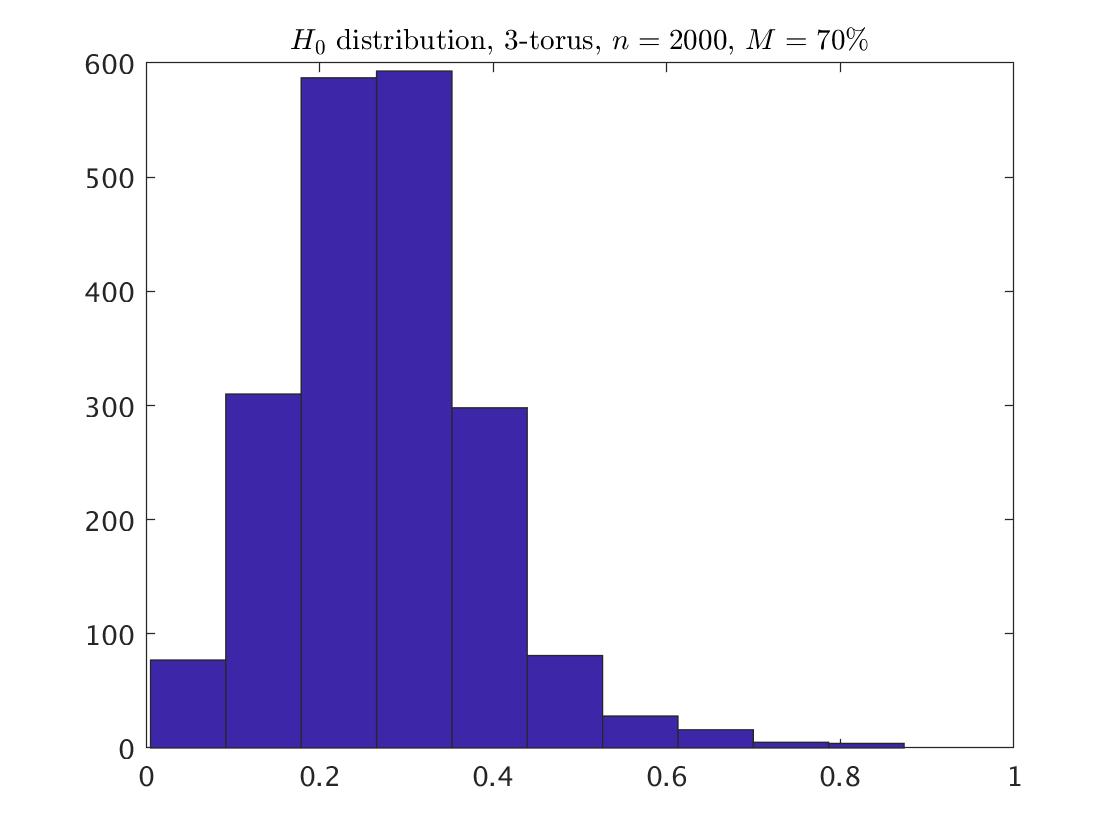} \hskip0.01truein
\includegraphics[width=1.45in, height=1.45in]{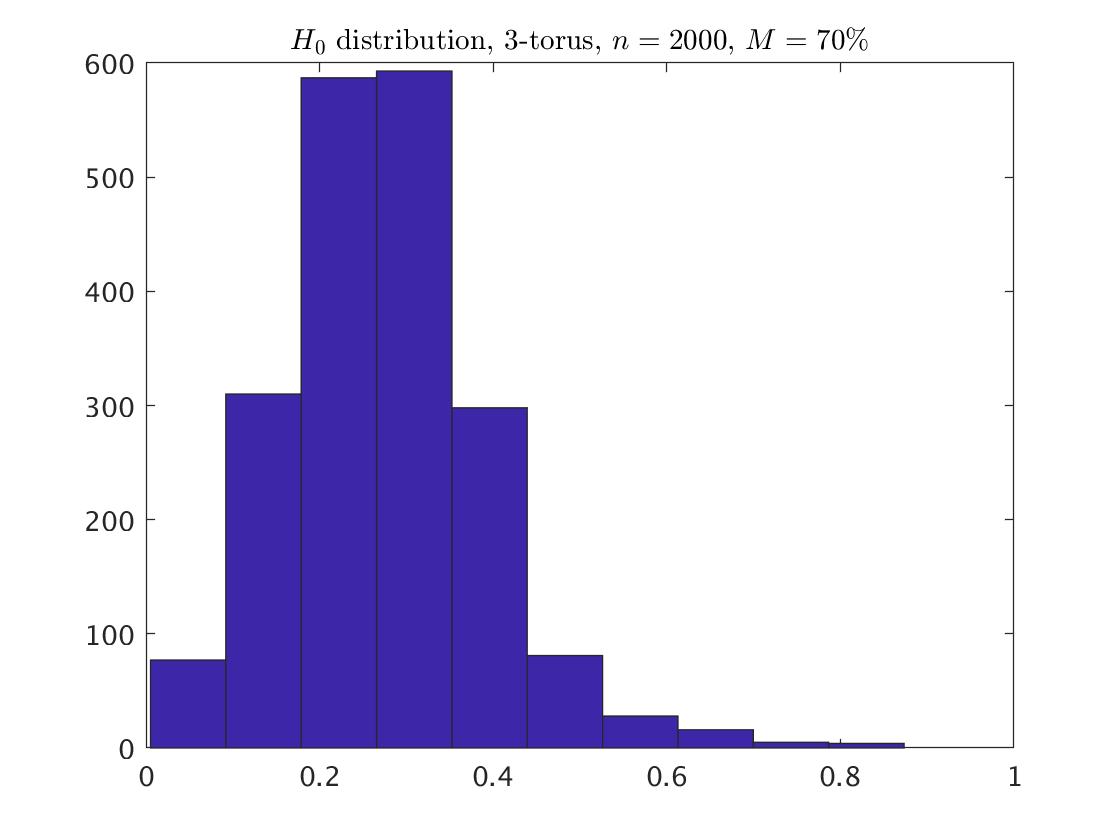} \hskip0.01truein

\includegraphics[width=1.45in, height=1.45in]{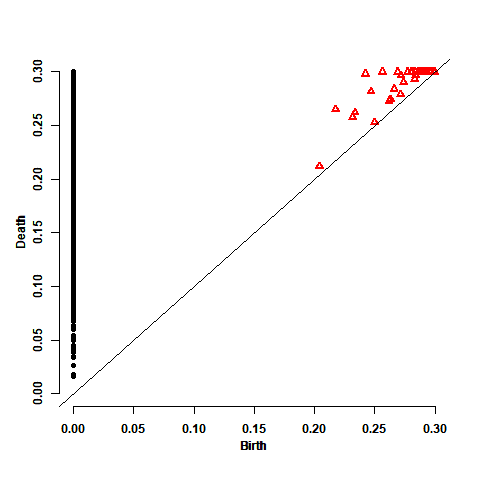} \hskip0.01truein
\includegraphics[width=1.45in, height=1.45in]{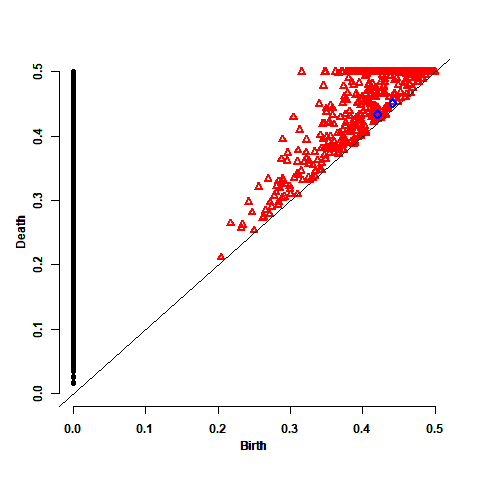} \hskip0.01truein
\includegraphics[width=1.45in, height=1.45in]{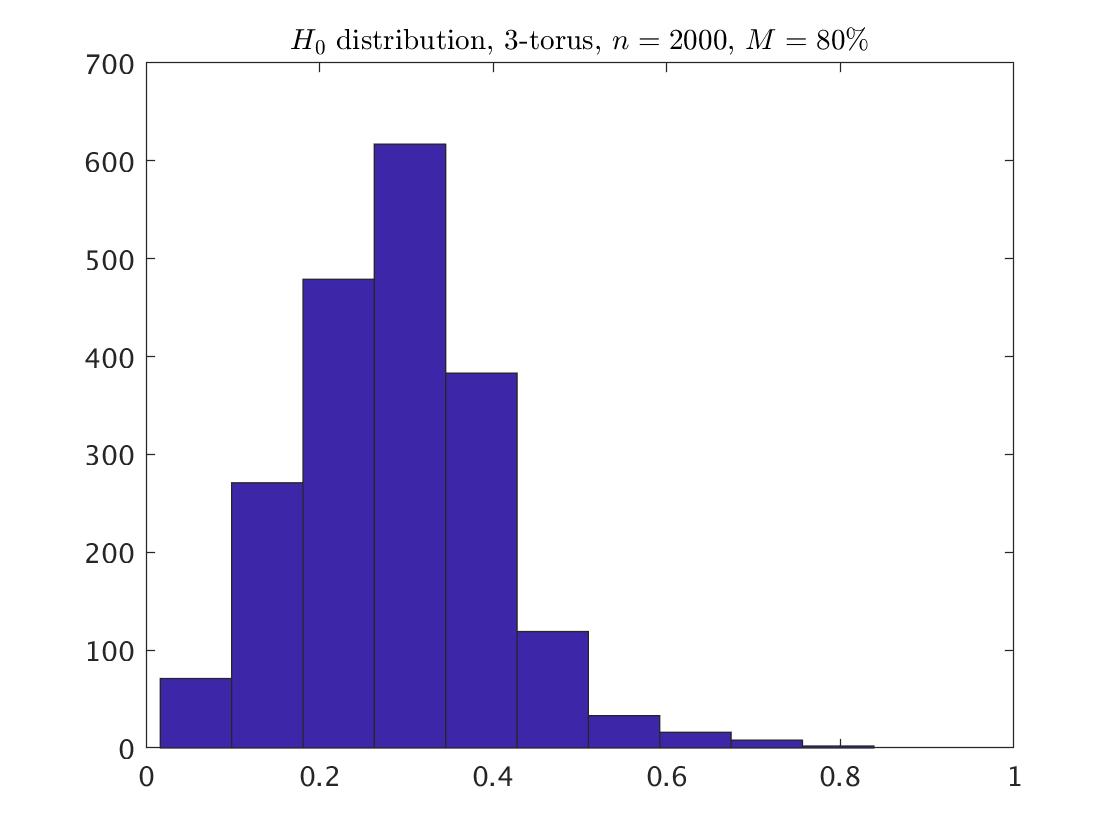} \hskip0.01truein
\includegraphics[width=1.45in, height=1.45in]{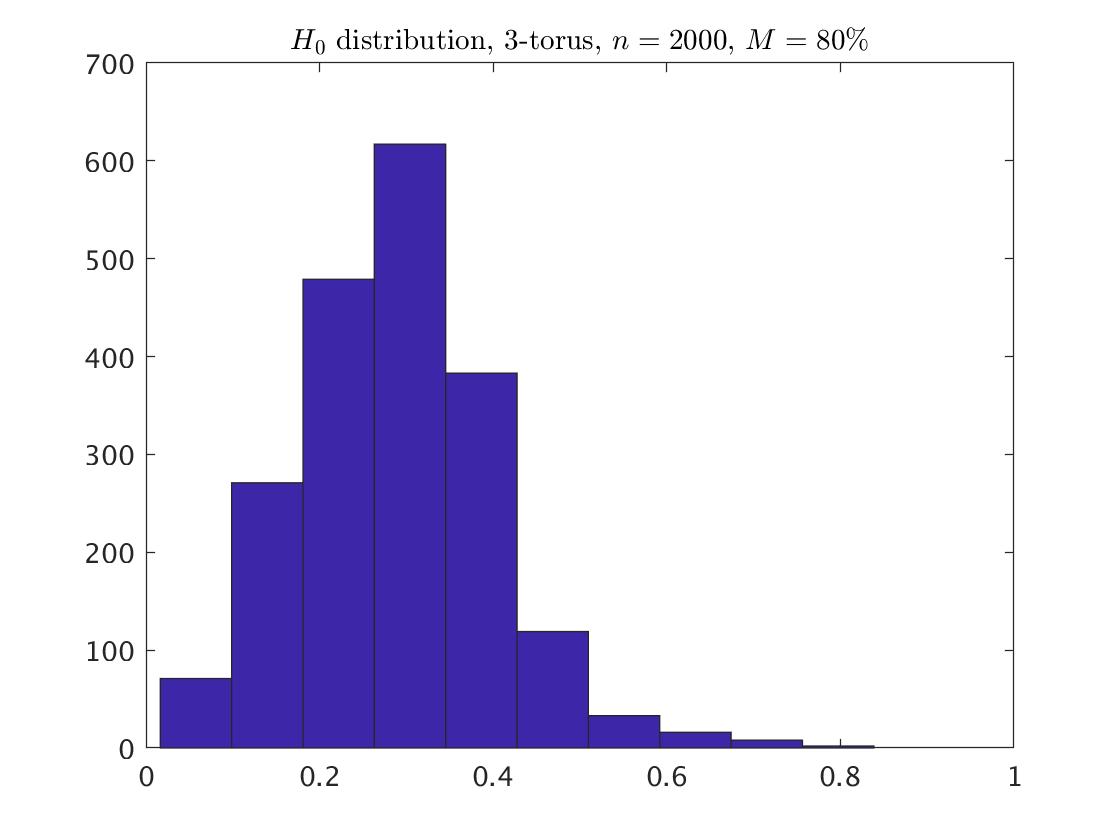} \hskip0.01truein

\includegraphics[width=1.45in, height=1.45in]{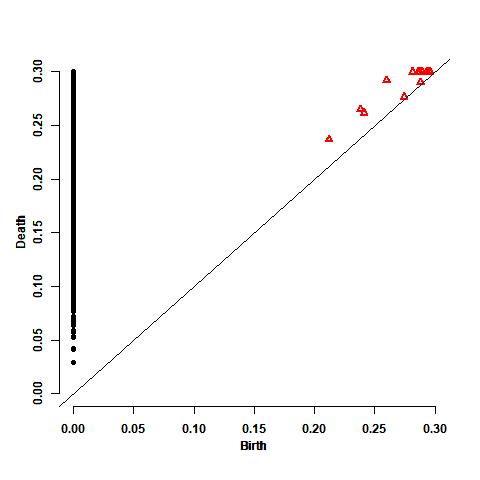} \hskip0.01truein
\includegraphics[width=1.45in, height=1.45in]{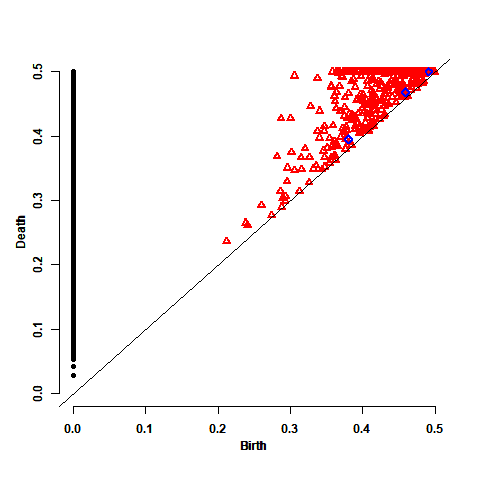} \hskip0.01truein
\includegraphics[width=1.45in, height=1.45in]{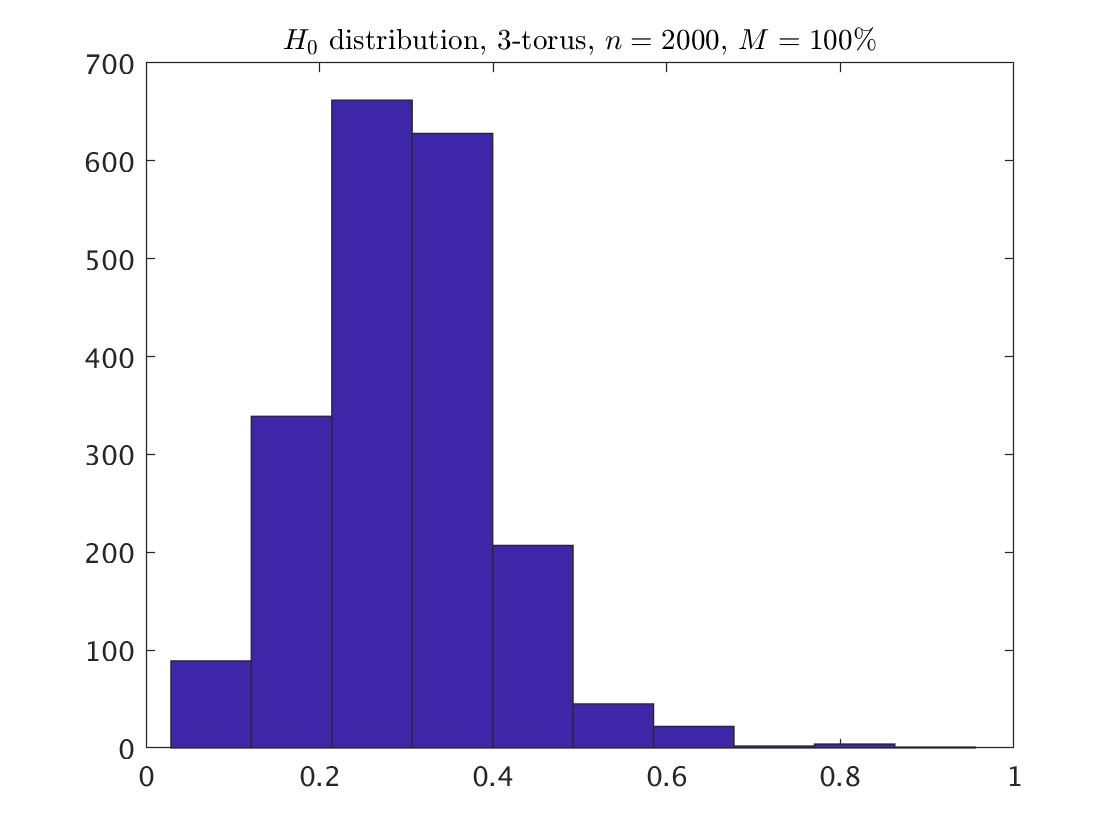} \hskip0.01truein
\includegraphics[width=1.45in, height=1.45in]{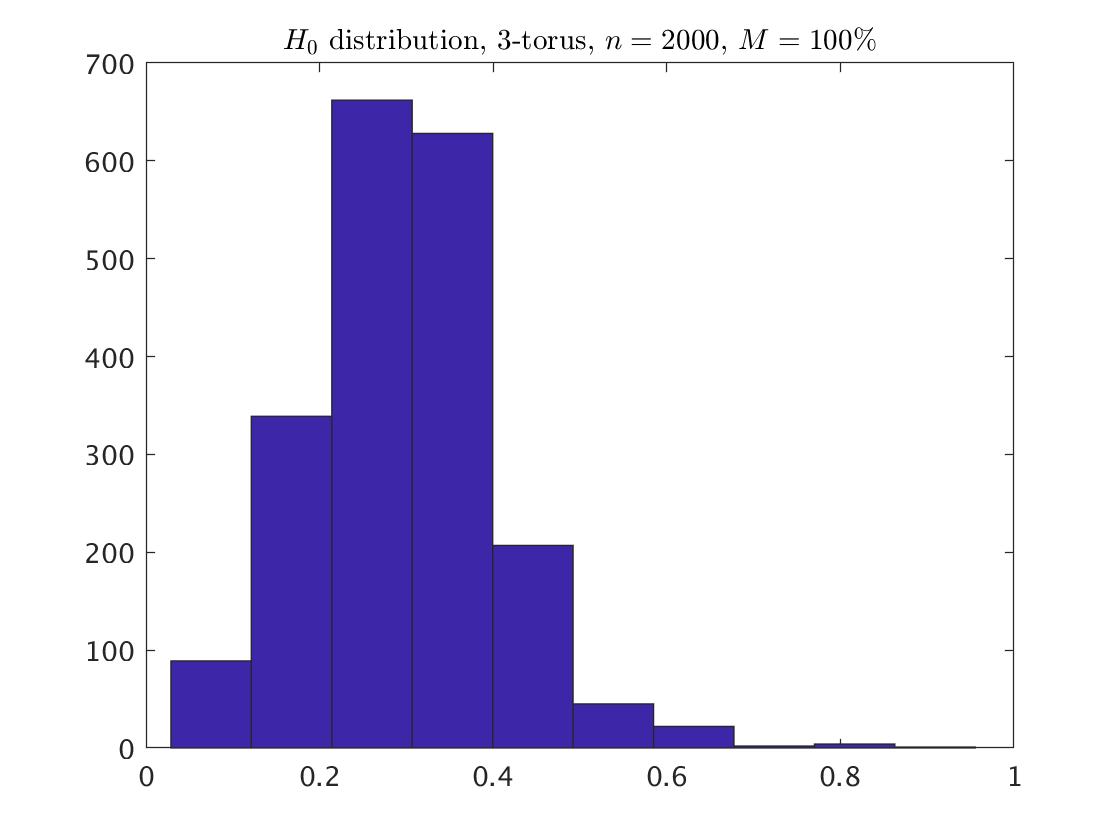} \hskip0.01truein
\ec
\caption{\footnotesize
 3-Torus with noise of $M\%$ of the sample size $n=2,000$. From top to bottom: $M = 30\%, 70\%, 80\%, 100\%$. Each row contains the persistence diagram based on maxscale=0.3, 0.5, and 1.
}
\label{fig:torus2000_M_noise}
\end{figure}
\end{landscape}

Here we examine the influence of noise on the example of the 3-torus, for the samples $n=1,500$ and $n=2,000$. The considered values of the maxscale are 0.3, 0.5, and 1. But, for some values of noise, the maxscale of 1 does not enough to cover the data of the 3-torus, and a larger value of 1.2 is needed.
Fig. \ref{fig:torus_M_noise} and Fig. \ref{fig:torus2000_M_noise} describe for each value of $M$ and the three (four) considered values of the maxscale, the 3-torus with $n=1,500$ and $n=2,000$ points, respectively, together with the distribution of the $H_0$ points (without the point at infinity). It is clear from the plots that the smaller values of the maxscale fail to capture the whole shape of the 3-torus (the $H_1$ points seem to be cut under the smaller values of the maxscale). The histogram at each row of Fig. \ref{fig:torus_M_noise} and Fig. \ref{fig:torus2000_M_noise} is the corresponded $H_0$ distribution based on the relevant larger value of the maxscale.

Again, as in the previous examples, we classified the bars into two groups of small and long bars. The results, as depending on noise and maxscale, as presented in Table \ref{table:torusLong}. For $n=2,000$, the proportion of long bars increases as the noise increases. However, for $n=1,000$, this proportion is similar for noise$\ge$ 70\%. The reason for this inconsistent result is the relative small sample as we already noted about above. Table \ref{table:torusLong} summarizes for each maxscale and a given noise, the $95$-th, $99$-th, and the $100$-th percentiles of the $H_0$ lengths. The different between the various values of maxscale for a given level of noise is at the $95$-th percentile of the $H_0$ lengths. Checking in more details the lengths under maxscale=0.3 and $M=30\%$ yields that the proportion of lengths that equal to 0.3 is 0.32 for $n=1,500$, and 0.205 for $n=2,000$. For $M=70\%$, again under maxscale=0.3, this proportion is 0.500 for $n=1,500$, and 0.382 for $n=2,000$. That is, maxscale of 0.3 at this example is relative very small.

\normalsize

\normalsize

\begin{center}
\fontsize{8.5}{0.9}\selectfont
\captionof{table}{3-Torus - classification of $H_0$ points}
\begin{tabular}{l|lc|cc}
\\
\\
$r=1$&$n=1,500$$^a$ &&$n=2,000$$^b$\\\hline
\\
Noise&short& long&short& long\\\hline
\\
\\
\\
\\
30\%& 0.971&0.029&0.962& 0.038\\
\\
\\
\\
\\
70\%& 0.959&0.041 &0.942&0.058\\
\\
\\
\\
\\
80\%& 0.961&0.039& 0.940&0.060 \\
\\
\\
\\
\\
100\%& 0.963&	0.037& 0.936&0.064\\
\label{table:torusLong}
\end{tabular}
\end{center}
\footnotesize{Classification of the $H_0$ lengths into short and long bars for 3-torus. The classification is based on one persistence diagram, with maxscale$\ge$ 1. $^a$ $c_{max}=0.527$, $^b$ $c_{max}=0.453$.}\\

\normalsize
\begin{center}
\fontsize{8.5}{0.9}\selectfont
\captionof{table}{3-torus - percentiles of $H_0$ points}
\begin{tabular}{l|l|c|c|c}
\\
\\
Noise&maxscale&percentile $\%$&$n=1,500$$^a$ &$n=2,000$$^b$\\\hline
\\
\\
\\
\\
30\%&0.3& 95&0.300&	0.300\\
\\
\\
\\
\\
&& 99&0.300&	0.300\\
\\
\\
\\
\\
&& 100&0.300&	0.300\\
\\
\\
\\
\\
& 0.5& 95&0.449&0.424\\
\\
\\
\\
\\
& & 99&0.500&0.500\\
\\
\\
\\
\\
& & 100&0.500&0.500\\
\\
\\
\\
\\
& 1& 99&0.653&0.586\\
\\
\\
\\
\\
& & 100&1&0.923\\
\\
\\
\\
\\
& 1.2& 100&1.084&0.923\\
& &&&\\\hline
\\
\\
\\
70\%&0.3& 95&0.300&	0.300\\
\\
\\
\\
\\
&& 99&0.300&0.300\\
\\
\\
\\
\\
&& 100&0.300&	0.300\\
\\
\\
\\
\\
& 0.5& 95&0.500&0.465\\
\\
\\
\\
\\
& & 99&0.500&0.500\\
\\
\\
\\
\\
&& 100&0.500&0.500\\
\\
\\
\\
\\
& 1& 95&0.512 &0.465\\
\\
\\
\\
\\
& & 99&0.687&	0.632\\
\\
\\
\\
\\
& & 100&1&	0.873\\
\\
\\
\\
\\
& 1.2& 100&1.126&	0.873\\
& &&&\\\hline
\\
\\
\\
80\%& 0.3& 95&0.300&	0.300\\
\\
\\
\\
\\
&& 99&0.300&0.300\\
\\
\\
\\
\\
&& 100&0.300&	0.300\\
\\
\\
\\
\\
& 0.5& 95&0.500&	0.473\\
\\
\\
\\
\\
& & 99&0.500&0.500\\
\\
\\
\\
\\

&& 100&0.500&	0.500\\
\\
\\
\\
\\
&1& 95&0.506&	0.473\\
\\
\\
\\
\\
&& 99&0.655&	0.623\\
\\
\\
\\
\\
&& 100&0.864&	0.839\\
& &&&\\\hline
\\
\\
\\
100\%& 0.3& 95&0.300&	0.300\\
\\
\\
\\
\\
&& 99&0.300&0.300\\
\\
\\
\\
\\
&& 100&0.300&	0.300\\
\\
\\
\\
\\
& 0.5& 95&0.500	&0.477\\
\\
\\
\\
\\
&1& 95&0.508&	0.477\\
\\
\\
\\
\\
&& 99&0.613&	0.611\\
\\
\\
\\
\\
&& 100&0.946&	0.955\\
\\
\\
\\
\\
\label{table:torusPercntiles}
\end{tabular}
\end{center}
\footnotesize{Percentiles of the $H_0$ lengths. The $H_0$ lengths are based on one persistence diagram, with maxscale$\ge$ 0.3. }\\
\normalsize

\section*{Acknowledgements}
I thank so much to Prof. Robert J. Adler for helpful conversations in various stages of this research.
This research was supported in part by the ISF-NSFC joint research program (grant No. 2539/17), and I gratefully acknowledge this support.

\newpage
\section*{Appendix}

\subsection*{A.1. \enspace Fitted distribution for clean data}
The following tables summarize the parameters of the fitted distributions for the $H_0$ points corresponded to each of the clean data examples. The fitted distribution for all the cases except those that are labeled with $^*$, is the beta distribution. The cases that are labeled with $^*$ have the generalized Pareto as the fitted distribution.

\begin{center}
\fontsize{8.5}{0.9}\selectfont
\captionof {table}{One circle}
\begin{tabular}{l|lcc|ccc}
&& $r=1$ &&&$r=3$ \\\hline
\\
\\
Noise=0\%& $n=500$ & $n=1,000$& $n=2,000$  & $n=500$ & $n=1,000$& $n=2,000$\\\hline
\\
\\
\\
\\
\\
\\
\\

maxscale=0.3& [0.976,77.302]&[1.032,164.246]& [0.996,316.773]&[0.957,24.658] & [1.022,53.565] &[0.991,104.391]\\
\\

\end{tabular}
\end{center}
\normalsize

\begin{center}
\fontsize{8.5}{0.9}\selectfont
\captionof{table}{Two concentric circles}
\begin{tabular}{l|l|c|c|cc}
Noise=0\%& $n=800$ &$n=1,200$& $n=2,400$\\\hline
\\
\\
\\
maxscale=0.3&[1.002,41.323]&[0.984,61.230] &[0.982,122.492] \\
\\
maxscale=0.5&[0.977,39.624]&[0.960,58.700]&[0.958,117.367] \\
\\
maxscale=1&[0.407,6.098] &[0.402,8.939]&[0.400,17.546] \\
\\
\\
\\
\\
\\
\\

\label{table:TwoCirclesPar}
\end{tabular}
\end{center}
\normalsize

\begin{center}
\fontsize{8.5}{0.9}\selectfont
\captionof{table}{Two distinct circles}
\begin{tabular}{l|l|ccc}
\\
\\
Noise=0\%& $n=600$ & $n=1,200$\\\hline
\\
\\
\\
\\
\\
\\
\\
\\
maxscale=0.3&[0.918,137.747]&[0.942,280.341]\\

\label{table:TwoDistinctPar}
\end{tabular}
\end{center}
\normalsize

\begin{center}
\fontsize{8.5}{0.9}\selectfont
\captionof{table}{Sphere $S^2$}
\begin{tabular}{l|l|ccc}
\\
\\
Noise=0\%& $n=1,000$ & $n=1,500$\\\hline
\\
\\
\\
\\
\\
\\
\\
\\
maxscale=0.3&[3.307,42.375]&[3.558,55.769]\\
\\
\\
\\
\\
\\
\\

\label{table:spherePar}
\end{tabular}
\end{center}
\normalsize

\begin{center}
\fontsize{8.5}{0.9}\selectfont
\captionof{table}{Torus $T^3$}
\begin{tabular}{l|l|ccc}
\\
\\
Noise=0\%& $n=1,500$ & $n=2,000$\\\hline
\\
\\
\\
\\
\\
\\
\\
\\
maxscale=0.3&[-1.326,0.389,0.007]$^*$ &[-1.173,0.347,0.004]\\
\\
maxscale=0.5&[3.211,12.935]&[3.413,16.392]\\
\\
maxscale=1&[3.209,12.922]\\
\label{table:TorusPar}
\end{tabular}
\end{center}
\normalsize

\subsection*{A.2. \enspace Fitted distribution for noisy data}
The following table summarizes the parameters of the fitted distribution for the $H_0$ points corresponded to the noisy circle example. The fitted distribution for all the cases except those that are labeled with $^*$, is the beta distribution. The cases that are labeled with $^*$ have the generalized Pareto as the fitted distribution.
\newpage
\begin{center}
\fontsize{8.5}{0.9}\selectfont
\captionof {table}{One circle}
\begin{tabular}{l|lcc|ccc}
&& $r=1$ &&&$r=3$ \\\hline
\\
\\
Noise& $n=500$ & $n=1,000$& $n=2,000$  & $n=500$ & $n=1,000$& $n=2,000$\\\hline
\\
\\
\\
\\
\\
\\
\\

30\%, maxscale=0.3& [0.607,11.120]&[0.522,13.586]&[0.454,17.439]&[0.854,8.857]&[0.693,10.181]&[0.571,12.204]\\
\\
30\%, maxscale=0.5&[0.597,10.686]&[0.513,12.963]&[0.452,17.217]&[0.787,7.529]&[0.660,9.185]&[0.563,11.787]  \\
\\
30\%, maxscale=1&&                               &[0.451,17.198] &[0.777,7.345]&[0.651,8.928]& \\
\\
\\
\\
\\
\\
\\
\\
\\
\\
70\%, maxscale=0.3 &[1.414,17.316]&[1.201,20.541]&[1.011,25.361]&[-0.831,0.252,0.001]$^*$&[1.509,14.321]&[1.289,17.717]\\
\\
70\%, maxscale=0.5&[1.390,16.872]&[1.170,19.717]& [1.001,24.922]&[1.668,10.724]&[1.440,13.346]&[1.266,17.219] \\
\\
70\%, maxscale=1&&                              & [1.000,24.886] &[1.665,10.694]&[1.436,13.284]&[1.266,17.219]\\
\\
\\
\\
\\
\\
\\
\\
\\
\\
80\%, maxscale=0.3&[1.784,20.858]&[1.508,23.897]& [1.398,32.000]&[-0.877,0.265,4.10e-05]$^*$&[1.853,16.509]&[1.664,21.315]\\
\\
80\%, maxscale=0.5 &[1.672,19.051]&[1.476,23.155]& [1.391,31.786]&[2.130,12.842]&[1.775,15.533]&[1.637,20.820]\\
\\
80\%, maxscale=1 &[1.665,18.935]&                   &&[2.090,12.525]&\\
\\
\\
\\
\\
\\
\\
100\%, maxscale=0.3&[2.926,30.497]&[2.804,41.598]&[2.644,57.753]&[-1.035,0.301,0.009]$^*$&[2.757,22.770]&[2.541,30.561] \\
\\
100\%, maxscale=0.5&[2.905,30.217]&[2.795,41.440]&[2.640,57.642]&[2.802,15.331]&[2.649,21.650]&[2.508,30.065] \\
\\
100\%, maxscale=1&&&&[2.795,15.284]&[2.632,21.473]& \\
\\
\end{tabular}
\end{center}
\normalsize

\end{document}